\documentclass[11pt]{article}%{ctexart}

\usepackage{graphicx}
\graphicspath{{figures/}}
\usepackage[margin=0.5in]{geometry}
\usepackage{bm}
\usepackage{amsthm}
\usepackage{amscd}
\usepackage{amsbsy}
\usepackage{appendix}
\usepackage{amsmath}
\usepackage{amsfonts}
\usepackage{amssymb}
\usepackage{calligra}
\usepackage{xcolor}
\usepackage{float}
\usepackage[T1]{fontenc}
\usepackage[utf8]{inputenc}
\usepackage{multirow}
\usepackage{mathrsfs}
\usepackage{mathtools}
\usepackage{epstopdf}
\usepackage{psfrag}
\usepackage{subfigure}
\usepackage{booktabs}
\usepackage{listings}
\usepackage{longtable}
\usepackage{latexsym}
\usepackage[shortlabels]{enumitem}
\usepackage{epsfig}
\usepackage{cite}
\usepackage{array}
\usepackage{algorithm}
\usepackage[noend]{algorithmic}
\usepackage{caption}
\usepackage{tikz}
\usetikzlibrary{positioning}
\IfFileExists{tikzlibraryarrows.meta.code.tex}{%
  \usetikzlibrary{arrows.meta}
  \tikzset{workflow forward/.style={-{Implies[length=3mm]}},
           workflow equivalence/.style={{Implies[length=3mm]}-{Implies[length=3mm]}}}
}{%
  \usetikzlibrary{arrows}
  \tikzset{workflow forward/.style={-implies},
           workflow equivalence/.style={implies-implies}}
}
\let\citep\cite
\usepackage{verbatim}
\usepackage{overpic}
\usepackage{abstract}
\usepackage{extarrows}
\usepackage[pagewise]{lineno}
\usepackage[none]{hyphenat}

\allowdisplaybreaks[4]

\date{}
\newtheorem {theorem*}{Theorem}
\newtheorem{theorem}{Theorem}[section]
\newtheorem{lemma}[theorem]{Lemma}
\newtheorem{corollary}[theorem]{Corollary}
\newtheorem{proposition}[theorem]{Proposition}
\theoremstyle{definition}
\newtheorem{definition}[theorem]{Definition}
\newtheorem{remark}[theorem]{Remark}

\theoremstyle{plain}

\newtheorem{open problem} {Open problem}

\numberwithin{equation}{section}

\usepackage[colorlinks,
            CJKbookmarks=true,
            linkcolor=blue,
            anchorcolor=red,
             urlcolor=blue,
            citecolor=blue
            ]{hyperref}
\hypersetup{
  pdftitle={Elimination Methods for Computing Invariant Algebraic Surfaces and Hopf Bifurcations},
  pdfauthor={Changbo Chen, Bo Huang, Valery Romanovski, and Linpeng Wang}
}

\DeclareMathOperator{\lt}{lt}

\DeclareMathOperator{\term}{term}
\DeclareMathOperator{\coef}{coef}

\DeclareMathOperator{\ini}{init}
\newcommand{\field}[1]{\mathbb{#1}}
\newcommand{\fk}{\field{K}}
\newcommand{\p}[1]{\mbox{{\em\bfseries #1}}}

\newcommand{\kx}{\fk[\p{x}]}
\newcommand{\zero}{\mathsf{Z}}

\newcommand{\pset}[1]{\mathcal{#1}}
\newcommand{\ide}[1]{\langle #1 \rangle}

\newcommand{\frI}{{\mathfrak{I}}}
\newcommand{\frJ}{{\mathfrak{J}}}

\newcommand{\caT}{{\mathcal{T}}}

\newcommand{\caP}{{\mathcal{P}}}
\newcommand{\caQ}{{\mathcal{Q}}}

\newcommand{\bases}[1]{\langle #1 \rangle}

\DeclareMathOperator{\lc}{lc}

\DeclareMathOperator{\sat}{sat}

\DeclareMathOperator{\prem}{{prem}}

\DeclareMathOperator{\nform}{nform}
\DeclareMathOperator{\cls}{cls}

\DeclareMathOperator{\mdeg}{mdeg}
\DeclareMathOperator{\mvar}{mvar}
\newcommand{\K}{{\mathbb{K}}}
\newcommand{\Q}{{\mathbb{Q}}}
\newcommand{\R}{{\mathbb{R}}}
\newcommand{\x}{{\boldsymbol{x}}}
\newcommand{\bfu}{{\boldsymbol{u}}}
\newcommand{\bfy}{{\boldsymbol{y}}}
\newcommand{\caR}{{\mathcal{R}}}
\newcommand{\caH}{{\mathcal{H}}}

\newcommand{\bfz}{{\boldsymbol{z}}}
\newcommand{\bmu}{{\boldsymbol{\mu}}}
\let\hat\widehat
\let\tilde\widetilde
\let\bar\overline

\title{Elimination Methods for Computing Invariant Algebraic Surfaces and Hopf Bifurcations\thanks{Changbo Chen is supported by the National Key Research Project of China under Grant No. 2023YFA1009402 and Chongqing Talents Plan Youth Top-Notch Project (2021000263). Bo Huang's work is supported by National Natural Science Foundation of China (NSFC 12671551, 12561026), Beijing Natural Science Foundation (1262015). Valery Romanovski acknowledges the support by the Slovenian Research and Innovation Agency (core research program P1-0306) and by the project 101183111-DSYREKI-HORIZON-MSCA-2023-SE-01 ``Dynamical Systems and Reaction Kinetics Networks''.}}
\author{Changbo Chen\textsuperscript{a}, Bo Huang\textsuperscript{b},\footnote{Corresponding author.}
  Valery Romanovski\textsuperscript{c,d,e}, and Linpeng Wang\textsuperscript{b}\\[0.4em]
  {\itshape\footnotesize \textsuperscript{a}Chongqing Key Laboratory of Automated Reasoning and Cognition,}\\
  {\itshape\footnotesize Chongqing Institute of Green and Intelligent Technology (CIGIT), Chinese Academy of Sciences,}\\
  {\itshape\footnotesize Chongqing 400714, China}\\
  {\itshape\footnotesize \href{mailto:chenchangbo@cigit.ac.cn}{chenchangbo@cigit.ac.cn}}\\[0.2em]
  {\itshape\footnotesize \textsuperscript{b}LMIB -- School of Mathematical Sciences,
    Beihang University, Beijing 100191, China}\\
  {\itshape\footnotesize \href{mailto:bohuang0407@buaa.edu.cn}{bohuang0407@buaa.edu.cn},
    \href{mailto:linpengwang@buaa.edu.cn}{linpengwang@buaa.edu.cn}}\\[0.2em]
  {\itshape\footnotesize \textsuperscript{c}Center for Applied Mathematics and Theoretical Physics, SI-2000 Maribor, Slovenia}\\
  {\itshape\footnotesize \textsuperscript{d}Faculty of Natural Science and Mathematics, University of Maribor, SI-2000 Maribor, Slovenia}\\
  {\itshape\footnotesize \textsuperscript{e}Faculty of Electrical Engineering and Computer Science,
    University of Maribor, SI-2000 Maribor, Slovenia}\\
  {\itshape\footnotesize \href{mailto:valerij.romanovskij@um.si}{valerij.romanovskij@um.si}}}

\begin{document}
\maketitle
\begin{center}
    \textit{Dedicated to Professor Dongming Wang on the occasion of his 65th birthday.}
\end{center}

\begingroup\small
\noindent \textbf{Abstract.} This paper presents efficient algorithmic methods for computing invariant algebraic surfaces and Hopf bifurcations in parametric polynomial differential systems. For the existence of a prescribed invariant algebraic surface of bounded degree, we establish two complementary elimination frameworks: combining Gr\"obner bases with saturation and prime decomposition yields a concise representation of the Zariski closure of the parameter existence locus (necessary conditions), whereas triangular decomposition into strong regular systems recovers the exact pointwise constructible existence locus (necessary and sufficient conditions) via the projection property. We further present a hybrid approach for identifying such surfaces by integrating prime decomposition with branch-wise triangular solving, thereby efficiently deriving necessary and sufficient conditions with enhanced computational flexibility and rectifying subtle classification errors in a class of Li\'enard systems.

For Hopf bifurcations, necessary conditions for an equilibrium point to be a Hopf bifurcation point were characterized in [SIAM J. Numer. Anal. {\bf 34} (1997)], while a class of sufficient conditions was independently established in [J. Math. Anal. Appl. {\bf 182} (1994)] and [J. Symbolic Comput. {\bf 30} (2000)]. In this paper, we provide a complete, necessary and sufficient algebraic characterization of Hopf bifurcation points via two equivalent formulations: one based on Hurwitz determinants and subresultant sequences, and the other derived from the greatest common divisor (GCD) of the real and imaginary parts of the characteristic polynomial. Using the polynomial GCD formulation, we establish a novel and, for the first time, complete semi-algebraic description of the transversality condition for Hopf bifurcation, leading to an algorithmic decision procedure via triangular decomposition of semi-algebraic systems and heuristic real quantifier elimination. We demonstrate that our approach can resolve cases that are out of reach for the state-of-the-art polynomial-based methods for Hopf bifurcations. The feasibility and performance of the proposed algorithmic frameworks are robustly validated by our experimental results.

\smallskip

\noindent {\bf Mathematics Subject Classification (2020).} 34C08; 37G10; 68W30

\smallskip

\noindent \textbf{Keywords.} Elimination methods, Gr\"obner basis, Hopf bifurcation, Hurwitz determinants, Invariant algebraic surface, Subresultant sequences, Triangular decomposition.

\par\endgroup
\clearpage
\tableofcontents

\section{Introduction}

Many real-world phenomena are modeled using autonomous systems of parametric differential equations. In this paper, we deal with polynomial differential systems of the form
\begin{equation}\label{L-dot_xx}
\dot{\boldsymbol{x}} = F(\boldsymbol{x}, \boldsymbol{\mu}),\quad F:\mathbb{K}^n\times\mathbb{K}^m\rightarrow\mathbb{K}^n,
\end{equation}
where $\mathbb{K}$ is an effective field of characteristic zero (typically, $\mathbb{K}=\mathbb{C}$ or $\mathbb{R}$), $\boldsymbol{x}=(x_1,\ldots,x_n)$ is a list of variables, $\boldsymbol{\mu}=(\mu_1,\ldots,\mu_m)$ is a list of real parameters, and $F(\boldsymbol{x}, \boldsymbol{\mu})=(F_1,\ldots,F_n)$ is a list of polynomials in $\mathbb{K}[\boldsymbol{x}, \boldsymbol{\mu}]$. We are concerned with two fundamental problems about differential systems of the form \eqref{L-dot_xx}. The first is to discuss the problem of computing invariant algebraic surfaces, while the second focuses on the algebraic description of Hopf bifurcations. More precisely, the present article considers the following two algorithmic problems:

($\mathcal{P}_1$): Given a prescribed form of an invariant algebraic surface of bounded degree for system \eqref{L-dot_xx}, we derive the conditions on the parameters $\mu$ under which the system admits this surface, or prove its nonexistence.

($\mathcal{P}_2$): Given a system of the form \eqref{L-dot_xx}: how to determine the critical parameter values of $\mu$ where a Hopf bifurcation occurs? What is the specific nature of this bifurcation?

The interest in the above two problems comes for two main reasons: one is the techniques used in their investigations are both based on the elimination theory, and the other is that there is a connection between the two problems, namely, the invariant algebraic surfaces of the system \eqref{L-dot_xx} may define its center manifold so that problem ($\mathcal{P}_1$) can be used to further study problem ($\mathcal{P}_2$) on the center manifold (see \cite{ayb2015inv,LiRomanovski2021} for instance). Throughout this paper, unless otherwise specified, problem ($\mathcal{P}_1$) is addressed over the complex field $\mathbb{K}=\mathbb{C}$, while problem ($\mathcal{P}_2$) is considered over the real field $\mathbb{K}=\mathbb{R}$.

\smallskip

\noindent{{\it Results related to problem ($\mathcal{P}_1$).}}

\smallskip

The study of the problem ($\mathcal{P}_1$) is one of the classical topics in the qualitative theory of differential equations. Given a planar differential system, the problem of finding a bound for the degree of irreducible invariant algebraic curves is known to be difficult: this is the so-called \textit{Poincar\'e problem}. It has been deeply studied in the literature and many partial results exist (see \cite{dar1878meo,poincare1886courbes,Cerveau1991,pre1983ele,Carnicer1994,Moulin-Ollagnier1995,Campillo1997MM,chr2000alg,Walcher2000,Pereira2002,Chavarriga2006,Yu2008Ser,Bostanetal2016,Demina2022} and references therein) though the problem has not yet been fully solved. Note that the uniform upper bound that depends entirely on the degree of the differential system under consideration may not exist. Indeed, there are polynomial differential systems that possess irreducible invariant algebraic curves with degrees depending on the coefficients of the system (see \cite{ChristopherJL2002,ZX:2017}).

Since 1983, several algorithms have been developed for finding invariant algebraic curves. Interest in invariant algebraic curves has been revived by the appearance of the Prelle--Singer’s method \cite{pre1983ele}. Substantial progress has been achieved in algorithmic design and computer algebra implementation. For example, Man and MacCallum \cite{Man1993,MacCallum1997} explain how we can use the Prelle-Singer method in practice. In particular, Man in \cite{Man1993} explains how to solve the involved polynomial system with a Gr\"obner basis. The resolution of this polynomial system is difficult because no particular structure of the considered system is known. Coutinho and Menasch\'e Schechter \cite{CMS2006,CMS2009} give algorithms that allow one
to test a given vector field (with rational coefficients) for the existence of invariant algebraic curves, and discuss the performance of their implementations in the computer algebra system \textsc{SINGULAR}. The bottleneck of Prelle--Singer’s method and all of its variants is the computation by undetermined coefficients of all irreducible invariant algebraic curves of bounded degree, which leads to solving a polynomial system. This yields an exponential number of reducible invariant algebraic curves (see \cite{Cheze2011}). More precisely, in \cite{Cheze2011}, Ch\'eze further shows that if a polynomial differential system (with integer coefficients) admits only finitely many irreducible invariant algebraic curves of degree at most $N$, then it is possible to compute all of them using the so-called ecstatic curve from \cite{Pereira2001}. To our knowledge, this is the first polynomial-time result on the computation of invariant algebraic curves. The importance of this result is mainly theoretical since the exponent in the polynomial complexity estimate is bigger than 10. More recently, Bostan and his co-workers present an efficient algorithm for computing the set of all irreducible invariant algebraic curves (of a given degree) of planar polynomial differential systems (without parameters). The algorithm compares well to the previous polynomial time algorithm in \cite{Cheze2011}.

%Other necessary conditions are contained in \cite{CGGL2003,CGG05} and also in works on inverse integrating factors \cite{CGGL2003,CFL10}.

% In addition, the method of Puiseux series can be generalized to the case of nonautonomous dynamical systems in the plane \cite{Demina:JPA2018} and higher-order autonomous ordinary differential equations \cite{DMV2020}.

Recently a novel method of finding invariant algebraic curves, which we shall call the method of Puiseux series, was introduced \cite{Demina2018}. This method was used to establish the general structure of irreducible invariant algebraic curves for a wide class of polynomial differential systems in the plane, see \cite{DMV20181,Demina2022,Chenetal2025}. Invariant algebraic curves also play a key role in the qualitative analysis of polynomial differential systems. For example, the inverse integrating factor forms a special type of invariant algebraic curve, and algebraic limit cycles appear as its factors (see \cite{GLV1996}). In \cite{GJ2007} the reader can find some open questions and some relations between the computation of invariant algebraic curves and Hilbert's 16th problem. For other results, references or applications of invariant algebraic curves the reader can consult for example \cite{DLA2006,Goriely2001,ZX:2017}.

%These pioneering ideas have stimulated extensive research over the past century. For example, Llibre and Zhang \cite{llibre2000invariant} described all invariant algebraic surfaces of the Rikitake system by solving linear partial differential equations along characteristic curves, obtaining parameter conditions for the existence of polynomial and rational first integrals. The same methodology was later applied to the classical Lorenz system \cite{llibre2002invariant}, leading to a complete classification of its Darboux polynomials, rational first integrals, and algebraic integrability. Invariant algebraic curves of large degree of the planar quadratic system were found in \cite{Christopher2005}. A novel method of Puiseux series, which works for planar systems, has been recently developed in \cite{Demina2022,Chenetal2025}.

Despite these advances, the algorithmic determination of invariant algebraic curves remains a challenging computational problem. One of the first systematic applications of computer algebra techniques to the computation of invariant algebraic surfaces in parametric polynomial differential systems appeared in \cite{ayb2015inv,adrt2016}. Their approach uses Gr\"obner basis elimination to derive algebraic conditions on the system parameters yielding existence of Darboux polynomials. This idea motivates us to formalize the problem from the perspective of computational algebraic geometry and to incorporate more tools from computer algebra for further analysis. This is exactly what we do in this paper: We develop computational algebra methods for determining invariant algebraic curves using elimination techniques and triangular decomposition. We show that the algebraic conditions obtained by Gr\"obner bases are necessary for the existence of a prescribed invariant algebraic surface of bounded degree. More concretely, we improve the elimination techniques used in \cite{ayb2015inv,adrt2016} and develop the method in a more general way. Furthermore, we demonstrate that triangular decomposition serves as a powerful alternative to Gr\"obner basis elimination, enabling the decomposition of the parameter space into irreducible components and yielding both necessary and sufficient conditions for the existence of invariant algebraic curves (see Section \ref{sec:dsys} for details).

%Recall that Gr\"obner bases \cite{buchberger1995g} and triangular decompositions \cite{wu1986o,lazard1991a,kalkbrener1993a,wang1993a,wang2001e} are two of the principal computational methods for studying polynomial ideals and algebraic varieties; see, for example, \cite{yang1994s,lazard1991a,kalkbrener1993a,faugere1999a,cox2005usi,shimoyama1996l}. While Gr\"obner basis elimination is well suited for deriving elimination ideals, triangular decomposition, in particular Wang's method \cite{wang1993a,wang2001e}, is capable of performing elimination over both equations and inequations, and has demonstrated high efficiency in solving positive-dimensional systems \cite{aubry1999t}.

\smallskip

\noindent{{\it Results related to problem ($\mathcal{P}_2$).}}

\smallskip

It should be noted that the literature on Hopf bifurcation and its applications is indeed quite extensive. A classical comprehensive source is Marsden and McCracken \cite{MJMM1976}. Among the monographs on differential equations and dynamical systems that discuss Hopf bifurcations, we mention only a few: Guckenheimer and Holmes \cite[Theorem 3.4.2]{GJHP1983}, Amann \cite[Theorem 26.21]{Amann:1990}, Han and Yu \cite[Chapter 4]{han2012yu}, and Chicone \cite[Section 8.3]{Chicone:2006}. Previous investigators have proposed a variety of approaches to the determination of Hopf bifurcation points, which may be divided into two broad classes: numerical detection of Hopf bifurcation points (see \cite{Reddien:1983,Roose:1985a,Roose:1985b,Garratt1991} and references therein) and symbolic computation of Hopf bifurcation points. Below, we review some known results from non-numeric studies of Hopf bifurcation points, which are particularly relevant to the problem ($\mathcal{P}_2$). An algebraic criterion for detecting a special type of Hopf bifurcation point in polynomial differential systems was introduced by Liu \cite{JMAA:1994}. At such a point, a pair of complex conjugate eigenvalues of the Jacobian matrix crosses the imaginary axis, while all other eigenvalues have negative real parts. Such type of Hopf bifurcation is called a simple Hopf bifurcation in order to distinguish it from the Hopf bifurcations with some other eigenvalues on the right half plane. Guckenheimer and his co-workers \cite{Guckenheimer:1996,Guckenheimer:1997} presented algorithms for detecting Hopf bifurcation points based on classical algebraic results regarding the roots of polynomial equations and properties of Kronecker products for matrices. Their algorithms give the necessary conditions for an isolated equilibrium point of a differential system to be a Hopf bifurcation point (see Remark \ref{Hopf-Lcomp}). On the other hand, being unaware of the result of Liu \cite{JMAA:1994}, El Kahoui and Weber \cite{Weber2000} rediscovered such a criterion for simple Hopf bifurcations independently and obtained a fully algebraic criterion for traditional Hopf bifurcation points. Moreover, they showed that the problem ($\mathcal{P}_2$) may be formalized as an existentially quantified formula in the language of the first-order theory of the ordered field of the real numbers. The technique has already been successfully applied to systems with mass action kinetics \cite{MCS:2009,JCP:2015}. In particular, using a combination of this algebraic criterion with techniques arising from stoichiometric network analysis, Errami and co-workers \cite{JCP:2015} devised an efficient algorithmic methods to detect Hopf bifurcation points in chemical reaction networks with symbolic rate constant. It is worth noting that the fully algebraic criterion for identifying Hopf bifurcation points in parametric polynomial vector fields was originally proposed by El Kahoui and Weber \cite{Weber2000}. Subsequent research by Chen and Moreno Maza \cite{ChenMM-2011} revealed that the condition presented in \cite[Theorem 3.5]{Weber2000} is merely a sufficient, rather than necessary and sufficient condition for Hopf bifurcation points. Chen and Moreno Maza attempted to derive a revised, complete algebraic characterization for this bifurcation scenario, as briefly outlined in \cite[Remark 3]{ChenMM-2011}. Nevertheless, we observe that this correction still leaves a critical gap in the rigorous algebraic description of Hopf bifurcation points for general parametric polynomial systems. This unresolved issue, along with our systematic remedy to establish a necessary and sufficient algebraic criterion, will be formally addressed in full detail in Theorem \ref{Hopf-th3} of the present work.

We also remark that all the above-mentioned works did not consider the nature of the Hopf bifurcations. A common way to determine the nature of the bifurcations presented in most textbooks is via computations of normal forms (see for example  \cite{Kuznetsov2004}). However, the computation of normal form is a very laborious procedure. Over the years, normal forms, and the ensuing reductions, were discussed in coordinate-independent settings (including an algorithm for their computation), see \cite{Scheurle2002,ZAMM2004}. In principle, these results may be used for the necessary computations to determine the nature of a Hopf bifurcation, but this general approach has the drawback of high computational expense. See \cite{NKSW:2020,TPVR:2026} for some new improvements on the computation of Poincar\'e--Dulac normal forms.

%The main contributions regarding the problem ($\mathcal{P}_2$) are as follows:

%We present a complete, constructive algebraic characterization of Hopf bifurcation points that admits two distinct yet mutually consistent representations: one is formulated in terms of Hurwitz determinants and subresultant sequences (building on the ideas of \cite{Weber2000,ChenMM-2011}), while the other relies on the greatest common divisor (GCD) of the real and imaginary parts of a characteristic polynomial (offering a novel perspective on this problem). We then give a semi-algebraic description for deciding whether a given differential system undergoes a Hopf bifurcation and devise an algorithmic procedure for solving such semi-algebraic systems using triangular decomposition. The proposed algorithmic approach not only derives the conditions for a Hopf bifurcation, but also characterizes the structure of the equilibrium equations under those conditions. For the detection of the nature of the Hopf bifurcations, we propose a symbolic-numeric approach based on making use of the projection method to compute Lyapunov coefficients and the routine {\sf SamplePoints} (within the {\sf RegularChains} library in \textsc{Maple}), which appears much simpler and more efficient.

The main theoretical and algorithmic contributions of this work are summarized as follows:

\begin{itemize}
    \item \textbf{Exact constructible locus for invariant algebraic surfaces ($\mathcal{P}_1$):} We formalize the distinction between the Zariski closure of the existence locus, obtained via Gr\"obner basis elimination and saturation (Proposition \ref{proposition:nontrivial-condition}), and the exact pointwise constructible existence locus, recovered via strong regular system decomposition and constructible-set difference (Proposition \ref{proposition:dsys-trimethod}). Their hybrid integration (Proposition \ref{proposition:GB+dsys}) achieves both global conciseness and pointwise exactness, enabling the classification of higher-degree surfaces and identifying hand-computation errors in the published literature (Remark \ref{hand-error}).

    \item \textbf{Complete algebraic characterization of Hopf bifurcation points ($\mathcal{P}_2$):} By extending Hurwitz determinant definitions to degenerate indices ($\Delta_{-1}=\Delta_0=a_0$, Theorem \ref{HSres}), we close a subtle gap in prior algebraic criteria \cite{Weber2000,ChenMM-2011} concerning the corner case $i=m$. This yields a complete necessary and sufficient characterization of Hopf bifurcation points (Theorem \ref{Hopf-th3} and Corollary \ref{coro-th-hp3}), complemented by a novel, mutually consistent characterization via polynomial GCDs (Theorem \ref{Hopf-gcd1}).

    \item \textbf{Semi-algebraic Hopf decision procedure and symbolic-numeric stability analysis ($\mathcal{P}_2$):} Using the polynomial GCD formulation, we establish a novel and, for the first time, complete semi-algebraic description of the Hopf transversality condition (Lemma \ref{Hopf-lemma-1} and Proposition \ref{Hopf-prop1}). Pairing triangular decomposition of semi-algebraic systems with a heuristic real quantifier elimination procedure via {\sf RealTriangularize} (Algorithm \ref{algorithm:rtd}), we obtain a decision procedure that projects exact parameter conditions while preserving the algebraic structure of equilibrium curves. To further determine the specific nature of the bifurcation, we integrate sample-point isolation with a symbolic-numeric projection algorithm (Algorithm \ref{alg:lyapunov}) to evaluate the first Lyapunov coefficient.

\end{itemize}

The rest of the paper is organized as follows. In Section \ref{sec:preli}, we explain the basic theory and concepts of elimination ideal, Gr\"obner basis and triangular decomposition together with some improvements of Hurwitz determinants and subresultant sequences. We analyze in Section \ref{sec:dsys} the connection between invariant algebraic surfaces of parametric differential systems and elimination methods. We first propose a Gr\"obner basis algorithm solving problem ($\mathcal{P}_1$), from which necessary conditions for the existence of a prescribed invariant algebraic surface are derived. A second algorithm, based on triangular decomposition, then yields necessary and sufficient conditions for such existence. As a byproduct, a hybrid method for computing invariant algebraic surfaces by combining Gr\"obner bases with triangular decomposition is introduced. Experimental results on the Belousov--Zhabotinsky reaction model and a class of Li\'enard differential systems, together with the connections of the two algorithms are also presented in Section \ref{sec:dsys}. In Section \ref{Hopf-sect4}, we study the necessary and sufficient conditions for Hopf bifurcation points, and develop efficient algorithm to compute the first Lyapunov coefficient. Then, in Section \ref{sect-5} we present algebraic criteria for Hopf bifurcations, and propose an algorithmic decision procedure for detecting Hopf bifurcations by using method of triangular decomposition of semi-algebraic systems. We demonstrate the efficiency of our algorithmic procedure on some examples. Finally, Section \ref{sect-6} offers concluding remarks.

\section{Preliminary Results}
\label{sec:preli}
% {\color{red}TODO for CC: This section has to be rewritten. There are some wrong statements and the notations are not unified. I plan to start with general definitions
% of regular chains, regular systems and state systems with stronger properties. }

%We begin by recalling some necessary definitions from elimination methods, including the theories of Gr\"obner basis and triangular decomposition.
% We then introduce the basic concepts of invariant algebraic surfaces and first integrals of polynomial differential systems.

% {\color{red}Comments by CC:1)should term include the coefficient or not?2)Is it necessary to define polynomial ring?}

In this section, we review some fundamental concepts of elimination ideal, Gr\"obner basis and triangular decomposition together with their properties that will be used throughout this paper. For more detailed descriptions, readers are referred to \citep{Wu2001, wan2001eli, Hubert-2003, eis1995com, Sturmfels2002, cox2005usi, cox2013ide}. We also discuss some improvements about the connections of Hurwitz determinants and subresultant sequences.

Let $\fk$ be a field and let $x_1, \ldots, x_n$ be indeterminates over $\fk$. A formal power product $x_1^{\alpha_1} \cdots x_n^{\alpha_n}$ ($\alpha_i \in \mathbb{Z}_{\geq 0}$) is called a \emph{term} in the variables $x_1, \ldots, x_n$, denoted by $\boldsymbol{x}^{\boldsymbol{\alpha}}$, where $\boldsymbol{x} = x_1, \ldots, x_n$ and $\boldsymbol{\alpha} = \alpha_1, \ldots, \alpha_n$. Furthermore, for $i = 1, \ldots, n$, we denote $\boldsymbol{x}_i = x_1, \ldots, x_i$ (the sequence of the first $i$ indeterminates).

A finite sum
\begin{equation*}
P = \sum_{\boldsymbol{\alpha}} c_{\boldsymbol{\alpha}} \boldsymbol{x}^{\boldsymbol{\alpha}},
\; \text{where}\,  c_{\boldsymbol{\alpha}} \in \mathbb{K} \setminus\{0\}, \; \boldsymbol{\alpha} \in \mathbb{Z}^n_{\ge 0}
\end{equation*}
is called a \emph{polynomial} in $x_1, \ldots, x_n$ over $\mathbb{K}$. We call $c_{\boldsymbol{\alpha}}$ the \emph{coefficient} of the term $\boldsymbol{x}^{\boldsymbol{\alpha}}$ in $P$, denoted by $\coef(P, \boldsymbol{x}^{\boldsymbol{\alpha}})$. We call $\boldsymbol{x}^{\boldsymbol{\alpha}}$ a term of $P$, and denote the set of all terms in $P$ by $\term(P)$. Let $\mathbb{F}$ be an extension field of $\mathbb{K}$, i.e., $\mathbb{K} \subseteq \mathbb{F}$. For any $\boldsymbol{a} \in \mathbb{F}^n$, define the evaluation homomorphism
\begin{equation*}
\operatorname{ev}_{\boldsymbol{a}}: \mathbb{F}[\boldsymbol{x}] \to \mathbb{F},  \;\,  \operatorname{ev}_{\boldsymbol{a}}(P):= \sum_{\boldsymbol{\alpha}}c_{\boldsymbol{\alpha}} \boldsymbol{a}^{\boldsymbol{\alpha}},
\end{equation*}
and denote the result of applying it to $P$ by $P(\boldsymbol{a})$. If $P(\boldsymbol{a}) = 0$, then $\boldsymbol{a}$ is called a \emph{zero} of the polynomial $P$ in $\mathbb{F}^n$.

The set of all polynomials in $x_1, \ldots, x_n$ over $\mathbb{K}$ with the usual addition and multiplication forms a commutative ring with unity, called the \emph{$n$-variate polynomial ring} over $\mathbb{K}$, denoted by $\mathbb{K}[\boldsymbol{x}]$. Unless otherwise stated, we fix the variable ordering as $x_1<\cdots<x_n$.

%Let $\mathbb{K}$ be a field. We write $\boldsymbol{x}$ for $x_1,\ldots,x_n$ or $(x_1,\ldots,x_n)$ and denote by $\mathbb{K}[\boldsymbol{x}] = \mathbb{K}[x_1, \ldots, x_n]$ the ring of polynomials in $n$ variables over $\mathbb{K}$. Our notational conventions are as follows: Uppercase letters (e.g., $G$, $T$) denote polynomial sets unless otherwise specified and gothic letters (e.g., $\mathfrak{I}$, $\mathfrak{J}$) are used to denote polynomial ideals. Throughout the paper, we fix the variable ordering as $x_1 < x_2 < \cdots < x_n$.

\subsection{Elimination ideal and Gr\"obner basis}

%For each nonzero polynomial $f \in \mathbb{K}[\boldsymbol{x}]$, $f$ can be written as $f = \sum a \boldsymbol{x}^{\boldsymbol{p}}$ where $\boldsymbol{p} = (p_1,\ldots,p_n)$ and $\boldsymbol{x^p}$ is a term, $a \in \mathbb{K}$ is called the \emph{coefficient} of $f$ in $\boldsymbol{x^p}$, denoted by $\coef(f, \boldsymbol{x^p})$. The greatest term of $f$ with respect to a fixed term ordering is called the \emph{leading term} of $f$, denoted by $\lt(f)$. We write $\operatorname{term}(f)$ for the set of all terms appearing in $f$. For a proper ideal $\mathfrak{I} \subseteq \mathbb{K}[\boldsymbol{x}]$, we denote by $\lt(\mathfrak{I})$ the set of leading terms of all nonzero elements of $\mathfrak{I}$.

\begin{definition}
Let $\mathfrak{I} \subseteq \kx$ be an ideal and $<$ be a term ordering. A finite set of polynomials $\{G_1, \ldots, G_s\}$ in $\mathfrak{I}$ is called a Gr\"obner basis of $\mathfrak{I}$ with respect to $<$ if $\langle \lt(G_1), \ldots, \lt(G_s) \rangle = \langle \lt(\mathfrak{I}) \rangle$, where $\lt(\mathfrak{I})$ denotes the set of leading terms of all polynomials in $\mathfrak{I}$. If a polynomial set $\mathcal{P}$ is a set of generators for the ideal $\mathfrak{I}$, we also say $\mathcal{G}$ is a Gr\"obner basis of $\mathcal{P}$ with respect to $<$.
\end{definition}

\begin{proposition}[{\cite[Section 2.5, Corollary 2]{cox2013ide}}]
\label{proposition:GBproperty}
Let $F$ be a polynomial in $\mathbb{K}[\boldsymbol{x}]$ and $\frI$ be an ideal. Let $\mathcal{G}$ be a Gr\"obner basis of $\frI$. Then the normal form of $F$ modulo $\mathcal{G}$ is unique, denoted by $\nform(F, \mathcal{G})$, and $F \in \mathfrak{I}$ if and only if $\nform(F, \mathcal{G}) = 0$.
\end{proposition}

\begin{definition}
Let $\mathfrak{I}$ be an ideal in $\mathbb{K}[\boldsymbol{x}]$. The \emph{$\ell$-th elimination ideal} of $\mathfrak{I}$ ($\ell = 1, \ldots,n-1$) is
\begin{equation*}
\mathfrak{I}_{\ell} := \mathfrak{I} \cap \mathbb{K}[x_1,  \ldots,  x_{n-\ell}].
\end{equation*}
\end{definition}

Here the variables are ordered so that the variables to be retained come first and the variables to be eliminated come last. Thus $\mathfrak{I}_{\ell}$ denotes the ideal obtained by eliminating the last $\ell$ variables.

Let $F \in \kx$ be a polynomial and $<$ be a term ordering. If for any integer $\ell = 1, \ldots, n-1$, the leading term $\lt_{<}(F)$ of $F$ with respect to $<$ satisfies the following property: if $\lt_{<}(F)$ belongs to the subring $\mathbb{K}[x_1, \ldots, x_{n-\ell}]$ then $F$ also belongs to this subring, then we say the term ordering $<$ is an \emph{elimination ordering}.

\begin{theorem}[{\cite[Section 3.1, Theorem 2]{cox2013ide}}]
\label{theorem:poly-alge-elim}
Let $\mathfrak{I}$ be an ideal in $\mathbb{K}[\boldsymbol{x}]$ and $\mathcal{G}$ be a Gr\"obner basis of $\mathfrak{I}$ with respect to an elimination ordering $<$. Then the polynomial set $\mathcal{G} \cap \mathbb{K}[x_1, \ldots, x_{n-\ell}]$ is a Gr\"obner basis of the elimination ideal $\mathfrak{I}_{\ell}$ with respect to $<$.
\end{theorem}

It is easy to verify that the lexicographic order is an elimination ordering. Let $\mathcal{G}=\{G_1, \cdots, G_s\}$ be a Gr\"obner basis. If for $i=1, \ldots, s$, each $G_i$ is monic and no term of $G_i$ is divisible by any $\lt(G_j)$ ($j \neq i$), then we say $\mathcal{G}$ is \emph{reduced}. When no ambiguity arises, we also say $\mathcal{G}$ is a reduced Gr\"obner basis.

For a given ideal $\mathfrak{I}$ and a fixed term ordering, the reduced Gr\"obner basis of $\mathfrak{I}$ with respect to this ordering is unique (see Proposition 6 in Section 2.9 of \cite{cox2013ide}).

\begin{definition}
\label{definition:ideal-quotient-sat}
Let $\frI$ and $\frJ$ be ideals in $\mathbb{K}[\boldsymbol{x}]$. The \emph{ideal quotient} and the \emph{saturation} of $\frI$ with respect to $\frJ$ are defined respectively as
\begin{equation*}
\frI:\frJ:=\{F \in \kx \mid FG \in \mathfrak{I}\ \text{for all}\ G\in \frJ\},
\end{equation*}
\begin{equation*}
\frI:\frJ^{\infty}:=\{F\in \kx \mid FG^N \in \mathfrak{I}\ \text{for some}\ N \geq 0\ \text{and all}\ G\in \frJ\}.
\end{equation*}
There exists a sufficiently large positive integer $N$ such that $\frI:\frJ^{\infty}=\frI:\frJ^{N}$. A polynomial $F$ is called \emph{regular modulo $\frI$} if $\frI :\langle F \rangle = \frI$.
\end{definition}

Let $\mathcal{F}$ be any subset of the ring $\mathbb{K}[\boldsymbol{x}]$ and let $\overline{\mathbb{K}}$ denote the algebraic closure of $\mathbb{K}$. In this paper, we call the set of common zeros of all polynomials in $\mathcal{F}$ in $\overline{\mathbb{K}}^n$ the \emph{zero set} of $\mathcal{F}$, denoted by
\begin{equation*}
\zero(\mathcal{F}) := \{ \boldsymbol{a} \in \overline{\mathbb{K}}^n \, |\,  F(\boldsymbol{a}) = 0,  \forall F \in \mathcal{F}  \}.
\end{equation*}
For an ideal $\frI$ in $\mathbb{K}[\boldsymbol{x}]$, regarded as a subset of $\mathbb{K}[\boldsymbol{x}]$, we can similarly define $\zero(\frI)$. In particular, if $\mathbb{K}$ itself is an algebraically closed field, then $\zero(\frI)$ coincides exactly with the (not necessarily irreducible) algebraic variety $\mathsf{V}(\frI)$ determined by the ideal $\frI$.

The zero sets of polynomial sets in $\mathbb{K}[\boldsymbol{x}]$ are the closed sets of the \emph{Zariski topology} on $\overline{\mathbb{K}}^n$; their complements are the \emph{Zariski-open sets}. The \emph{Zariski closure} of a set $W\subseteq\overline{\mathbb{K}}^n$, denoted by $\overline{W}$, is the smallest Zariski-closed set containing $W$ (see, e.g., \cite[Chapter~4]{cox2013ide}). A set is called \emph{locally closed} if it is the intersection of a Zariski-closed set and a Zariski-open set, and \emph{constructible} if it is a finite union of locally closed sets.

\begin{definition}
Let $\frI$ be an ideal in $\mathbb{K}[\boldsymbol{x}]$. Define the \emph{radical} of $\frI$ as
\begin{equation*}
\sqrt{\frI} := \{ F \in \kx \, |\,  F^m \in \frI ,  \exists\,  m \in \mathbb{N} \}.
\end{equation*}
If $F^m \in  \frI \implies F \in \frI$, or equivalently $\frI = \sqrt{\frI}$, then we say $\frI$ is a \emph{radical ideal}.
\end{definition}

\begin{definition}
An ideal $\mathfrak{p}$ in $\mathbb{K}[\boldsymbol{x}]$ is called a \emph{prime ideal} if for any $F$, $G \in \kx$,
\begin{equation*}
F G \in \mathfrak{p} \implies F \in \mathfrak{p} \;\text{or}\; G \in \mathfrak{p}.
\end{equation*}
A zero set is called \emph{irreducible} if it cannot be written as the union of two strictly smaller zero sets. If
$$\zero(\frI) = \zero(\frI_1)\cup \cdots \cup \zero(\frI_s)$$
with each $\zero(\frI_i)$ irreducible and $\zero(\frI_i) \nsubseteq \zero(\frI_j)$ for $i\neq j$, then this expression is called a \emph{minimal irreducible decomposition} of $\zero(\frI)$.
\end{definition}

Since all zero sets are defined over the algebraic closure $\overline{\mathbb{K}}$, Hilbert's Nullstellensatz implies that this decomposition is unique up to the order of components.

\subsection{Triangular decomposition}

Let $F$ be a nonzero polynomial in $\mathbb{K}[\boldsymbol{x}]$.
We assume that the variable ordering is given by $x_1<\cdots<x_n$.
For any $i=1, \ldots, n$, $F$ can be written as
\begin{equation*}
F = c_kx_i^k + c_{k-1} x_i^{k-1} + \cdots+ c_1x_i + c_0,
\end{equation*}
where $c_j \in \mathbb{K}[x_1, \ldots, x_{i-1}, x_{i+1}, \ldots, x_n]\;(j=0, \ldots, k)$ and $c_{k}\neq 0$.
We call the integer $k$ the degree of $F$ with respect to variable $x_i$, denoted by $\deg(F, x_i)$; we call the coefficient $c_k$ the \emph{leading coefficient} of $F$ with respect to variable $x_i$, denoted by $\lc(F, x_i)$.

Furthermore, define the \emph{class} of $F$ as the largest index of variables appearing in $F$:
\begin{equation*}
\cls(F) := \max\{ i \, |\,  \deg(F, x_i) >0,  i=1, \ldots, n \};
\end{equation*}
if $F\in \mathbb{K}$, we define $\cls(F) = 0$. When $\cls(F) > 0$, define the \emph{main variable} and \emph{main degree} of $F$ respectively as
% \begin{equation*}
% \lv(F) := x_{\cls(F)},  \quad \ldeg(F) := \deg(F,  \lv(F)).
% \end{equation*}
\begin{equation*}
\mvar(F) := x_{\cls(F)},  \quad \mdeg(F) := \deg(F,  \mvar(F)).
\end{equation*}
Define the \emph{initial} of $F$ as
% \begin{equation*}
% \ini(F) := \lc(F, \lv(F)),  \quad \tail(F) := F-\ini(F)\cdot \lv(F)^{\ldeg(F)}.
% \end{equation*}
\begin{equation*}
\ini(F) := \lc(F, \mvar(F)).
\end{equation*}
For a polynomial set $\mathcal{F}$ in $\mathbb{K}[\boldsymbol{x}]$, denote the set of all initials of polynomials in $\mathcal{F}$ by $\ini(\mathcal{F}) := \{ \ini(F) \, |\,  F \in \mathcal{F} \}$.

Let $F, G$ be two polynomials in $\mathbb{K}[\boldsymbol{x}]$, $x_i$ be a variable, and let $\ell = \deg(G, x_i)$, $m=\deg(F, x_i)$. If $\ell >0$, then there exist two polynomials $Q$ and $R$ in $\mathbb{K}[\boldsymbol{x}]$ and an integer $0 \leq s \leq m-\ell+1$ such that
\begin{equation*}
\lc(G, x_i)^s F = QG+R, \,  \text{and } \deg(R, x_i) < \ell.
\end{equation*}
If $s$ is fixed, then $Q$ and $R$ are uniquely determined. We call $R$ the \emph{pseudo-remainder} of $F$ by $G$ with respect to $x_i$, denoted by $\prem(F, G, x_i)$. If $x_i = \mvar(G)$, we also simply denote it by $\prem(F, G)$. In this case, if $m < \ell$, we say \emph{$F$ is reduced with respect to polynomial $G$}. It is easy to see that $\prem(F,  G)$ is always reduced with respect to $G$.

\begin{definition}
Let $\mathcal{T} =  [T_1, \ldots, T_r]$ be a finite nonempty ordered set in $\mathbb{K}[\boldsymbol{x}]$. If $0<\cls(T_1)<\cdots<\cls(T_r)$, then we call $\mathcal{T}$ a \emph{triangular set}.
\end{definition}

It is easy to see that for a triangular set $\mathcal{T} = [T_1, \ldots, T_r]$, the ordered set consisting of the first $i$ polynomials ($i=1, \ldots, r$), denoted by $[T_1, \ldots, T_i]$, is also a triangular set, denoted by $\mathcal{T}^{[i]}$.
Let $P$ be a polynomial in $\mathbb{K}[\boldsymbol{x}]$. Define the \emph{pseudo-remainder of $P$ with respect to triangular set $\mathcal{T}$} as:
\begin{equation*}
\prem(P, \mathcal{T})=\prem(\, \ldots \,\prem(\, \prem(P, T_r), T_{r-1}), \ldots, T_1).
\end{equation*}
If $P$ is reduced with respect to $T_1, \ldots, T_r$, then we say $P$ is reduced with respect to triangular set $\mathcal{T}$. It is easy to see that the polynomial $\prem(P, \mathcal{T})$ is always reduced with respect to $\mathcal{T}$.

Let $\mathcal{F}$ and $\mathcal{H}$ be two polynomial sets in $\mathbb{K}[\boldsymbol{x}]$. We introduce the following notation to describe the set of points in the zero set of $\mathcal{F}$ where all polynomials in $\mathcal{H}$ are nonzero:
\begin{equation*}
    \zero(\mathcal{F}/ \mathcal{H})
    := \{ \boldsymbol{a} \in \overline{\mathbb{K}}^n \, |\,  \text{for any}\, F \in \mathcal{F},\, F(\boldsymbol{a}) = 0, \,\text{and for any}\, H \in \mathcal{H},\, H(\boldsymbol{a}) \ne 0 \}.
\end{equation*}

The set $\ini(\mathcal{T}) := \{\ini(T) \, |\,  T\in \mathcal{T} \}$ consists of all initials of polynomials in triangular set $\mathcal{T}$. From the pseudo-division formula of polynomials with respect to triangular sets, we can easily obtain the following zero set inclusion relationship.

\begin{proposition}
\label{proposition:prem-zero}
For any triangular set $\mathcal{T}$ and polynomial $P$ in $\mathbb{K}[\boldsymbol{x}]$, if $\prem(P, \mathcal{T})=0$, then
$\zero(\mathcal{T}/\ini(\mathcal{T})) \subseteq \zero(P).$
\end{proposition}

% \begin{proposition}[{\cite[Lemma 2.1.1]{wan2001eli}}]
% \label{proposition:prem-zero}
% For any triangular set $\mathcal{T}$ and polynomial $P$ in $\mathbb{K}[\boldsymbol{x}]$, if $\prem(P, \mathcal{T})=0$, then
% $\zero(\mathcal{T}/\ini(\mathcal{T})) \subseteq \zero(P).$
% \end{proposition}

For any two polynomial sets $\mathcal{F}$ and $\mathcal{H}$ in $\mathbb{K}[\boldsymbol{x}]$, define the projection of $\zero(\mathcal{F}/\mathcal{H})$ to the first $i$ coordinates as
\begin{equation*}
    \operatorname{Proj}_i \zero(\mathcal{F}/\mathcal{H})
    := \{ (a_1,  \ldots,  a_i) \in \overline{\mathbb{K}}^i \, |\,
        \exists\,  a_{i+1},  \ldots,  a_n \in \overline{\mathbb{K}}, \
        (a_1,  \ldots,  a_n) \in \zero(\mathcal{F}/\mathcal{H})
    \}.
\end{equation*}
For any polynomial set $\mathcal{F}$, denote by $\mathcal{F}^{(\le i)} := \mathcal{F}\cap\K[x_1,\ldots,x_{i}]$ its truncation to the first $i$ variables.
We say $[\mathcal{F},\mathcal{H}]$ has the \emph{projection property} if
\begin{equation}
\label{equation:projection}
\zero(\mathcal{F}^{(\le i)}/\mathcal{H}^{(\le i)}) = \operatorname{Proj}_i \zero(\mathcal{F} / \mathcal{H})
\end{equation}
holds for every $i = 1, \ldots, n-1$.

% Let $\mathcal{T}$ be a triangular set in $\mathbb{K}[\boldsymbol{x}]$.
% % For $i=1, \ldots, n-1$ denote $\mathcal{T}^{[i]} := \mathcal{T} \cap \field{K} [x_1, \ldots, x_i]$, i.e., the subset of $\mathcal{T}$ containing only variables $x_1,  \ldots,  x_i$.
% If for any $i = 1,  \ldots,  n-1$, we have
% % \begin{equation}
% % \label{equation:projection}
% % \zero(\mathcal{T}^{[i]}/ \ini(\mathcal{T}^{[i]})) \subseteq \operatorname{Proj}_i \zero(\mathcal{T} / \ini(\mathcal{T})),
% % \end{equation}
% %Comments by CC: using $=$ is better than $\subseteq$.
% \begin{equation}
% \label{equation:projection}
% \zero(\mathcal{T}^{[i]}/ \ini(\mathcal{T}^{[i]})) = \operatorname{Proj}_i \zero(\mathcal{T} / \ini(\mathcal{T})),
% \end{equation}
% then we say $\mathcal{T}$ has the \emph{projection property}.

\begin{definition}
Let $\pset{T} = [T_1, \ldots, T_r]$ be a triangular set in $\mathbb{K}[\boldsymbol{x}]$.
If $\pset{T}$ is empty, define $\sat(\pset{T})=\langle 0\rangle$. Otherwise, define its \emph{saturated ideal} $$\sat(\pset{T}) := \bases{\pset{T}}: \langle J\rangle^{\infty}, $$ where $J=\ini(T_1)\cdots \ini(T_r)$. If for each $i=1, \ldots, r$, the initial $\ini(T_i)$ is regular modulo the saturated ideal $\sat(\pset{T}^{[i-1]})$, then we call the triangular set $\mathcal{T}$ a \emph{regular chain}.
\end{definition}

\begin{proposition}
%[{\cite[Proposition 5.1.8]{wan2001eli}}]
\label{proposition:tri-projection}
Let $\mathcal{T}$ be a regular chain in $\mathbb{K}[\boldsymbol{x}]$. Then $\overline{\zero(\mathcal{T}\big/\ini(\mathcal{T}))} = \zero(\sat(\mathcal{T})) \neq \emptyset$.
\end{proposition}

\begin{theorem}[{\cite[Theorem 6.1]{aub1999ont}}]
\label{proposition:regular-tri-regular=sat}
A triangular set $\mathcal{T}$ in the polynomial ring $\mathbb{K}[\boldsymbol{x}]$ is a regular chain if and only if $\sat(\mathcal{T}) = \{F \, |\,  \prem(F, \mathcal{T}) = 0\}$.
\end{theorem}

\begin{definition}
Let  $\mathcal{T}$ be a triangular set in $\mathbb{K}[\boldsymbol{x}]$  and $\mathcal{H}\subseteq \kx$.
The pair $[\mathcal{T},\mathcal{H}]$  is called a \emph{triangular system}. For any $i=1,\ldots, n$, let $\mathcal{T}^{(i)}$ (resp. $\mathcal{T}^{(<i)}$) and $\mathcal{H}^{(i)}$ (resp. $\mathcal{H}^{(<i)}$) denote respectively the set of polynomials in $\mathcal{T}$ and $\mathcal{H}$ whose main variables are $x_i$ (resp. less than $x_i$).
The pair $[\mathcal{T},\mathcal{H}]$ is called a \emph{(weak) regular system} if $\mathcal{T}$ is a regular chain and each polynomial in $\mathcal{H}$ is regular modulo ${\rm sat}(\mathcal{T})$.
The pair $[\mathcal{T},\mathcal{H}]$ is called a \emph{strong regular system} if the following properties hold for any $i=1,\ldots, n$:
$(1)$ either $\mathcal{T}^{(i)}=\emptyset$ or $\mathcal{H}^{(i)}=\emptyset$;
$(2)$ at any point ${\bf a}^{(<i)}=(a_1,\ldots, a_{i-1})$ of $\zero(\mathcal{T}^{(<i)}/\mathcal{H}^{(<i)})$ in $\overline{\mathbb{K}}^{i-1}$,
for any $T\in \mathcal{T}^{(i)}$ and any $H\in \mathcal{H}^{(i)}$, we have $\ini(T)({\bf a}^{(<i)})\neq 0$ and $\ini(H)({\bf a}^{(<i)})\neq 0$.
\end{definition}

Strong regular system and weak regular system were respectively proposed in~\cite{wan2000com} and~\cite{DBLP:conf/casc/ChenGLMP07}.
It has already been remarked in~\cite{DBLP:conf/casc/ChenGLMP07} that a weak regular system may not be a strong one.
For instance, the system $[\{(x_1+x_2)x_3-1, x_2\},\{x_1+x_2\}]$
is a weak but not a strong regular system under the variable ordering $x_1<x_2<x_3$.
The system $[\{(x_1+x_2)x_3-1, x_2\},\{x_1\}]$ is a strong regular system.
%and a simple example $[\{x_1x_4-x_2\}, \{x_2x_3\}]$ is also provided.
Next we show that a strong regular system is also a weak one.

\begin{lemma}
\label{strong-weak}
Let $R :=[\mathcal{T},\mathcal{H}]$ be a strong regular system.
Then $R$ is also a weak regular system.
\end{lemma}
\begin{proof}
   By Theorem 5.1 of~\cite{wan2000com}, ${\cal T}$ is a regular chain.
   By Corollary 5.1 of \cite{wan2000com}, each polynomial $P$ in ${\cal H}$
   is regular modulo ${\rm sat}(\mathcal{T})$.
\end{proof}

\begin{proposition}[{\cite[Proposition 1]{DBLP:conf/casc/ChenGLMP07}}]
\label{prop:regular}
For any regular system $R :=[\mathcal{T},\mathcal{H}]$,
we have $\mathsf{Z}(R)\neq \emptyset$.
\end{proposition}

\begin{definition}
Let $\mathcal{F}\subseteq\kx$ and $\mathcal{G}\subseteq\kx$ be two polynomial sets. If there exist finitely many triangular systems such that
\begin{equation}
\label{equation:tri-dec}
\mathsf{Z}(\mathcal{F}/\mathcal{G})
=
\bigcup_{i=1}^{s}\mathsf{Z}(\mathcal{T}_i/\mathcal{H}_i),
\end{equation}
then we call the set $\{[\mathcal{T}_1,\mathcal{H}_1],\ldots,[\mathcal{T}_s,\mathcal{H}_s]\}$ a \emph{triangular decomposition} of $\mathsf{Z}(\mathcal{F}/\mathcal{G})$. If every $[\mathcal{T}_i,\mathcal{H}_i]$ is a (strong) regular system,
the decomposition is said to be (strong) regular. When $\zero(\mathcal{G})=\emptyset$ and $\mathcal{H}_i=\ini({\cal T}_i)$ for $i=1,\ldots,s$, we also call $\{\mathcal{T}_1,\ldots,\mathcal{T}_s\}$ a \emph{triangular decomposition} of $\mathsf{Z}(\mathcal{F})$.
\end{definition}

Many algorithms have been proposed to compute different types of triangular decompositions. In this work, we are primarily interested in the case that
triangular systems in the decomposition are regular systems,
which have the nice property that their zero sets are always nonempty by Proposition \ref{prop:regular}.

\begin{remark}\label{remark:cs}
In practice, we may encounter systems more general than $[\mathcal{F},\mathcal{G}]$. We call a finite union of differences of algebraic varieties in $\overline{\mathbb{K}}^n$ a {\em constructible set} of $\overline{\mathbb{K}}^n$. In Lemma 5 of \cite{DBLP:conf/casc/ChenGLMP07}, it is proved that every constructible set can be written as a finite union of zero sets of regular systems.
Moreover, an algorithm is provided to compute such a representation based
on a dedicated operation {\sf Difference} for computing the difference of the zero sets of two regular systems (see {\cite[Algorithm 1]{DBLP:conf/casc/ChenGLMP07}}). The concept of ``specializing well'' (see Definition \ref{def:spec}) is also proposed in the same paper.
This notion allows one to compute a refined triangular decomposition of a constructible set $cs$
into the union of zero sets of finitely many regular system $\{R_i\}$, called a pre-comprehensive triangular decomposition,
from which the projection of $cs$ onto the parameter space
can be obtained by taking the union of the defining sets of $R_i's$.

In \cite{wan2000com}, Algorithm {\sf RegSer} is proposed to decompose systems in the form of $[\mathcal{F},\mathcal{G}]$ into strongly regular systems.
As a result, any constructible set can be written as a union of zero sets of strong regular systems. From the definition, it is immediate that strong regular systems enjoy projection property. Note that weak regular systems may not enjoy projection property, as shown by the example $[\{(x_1+x_2)x_3-1, x_2\},\{x_1+x_2\}]$ given earlier.
\end{remark}

% \begin{definition}
% Let $\mathcal{F}\subseteq \kx$ be a polynomial set. If there exist finitely many triangular sets $\mathcal{T}_1, \ldots, \mathcal{T}_s \subseteq \kx$ satisfying the following zero relation:
% \begin{equation}
% \label{equation:tri-dec}
% \mathsf{Z}(\mathcal{F})= \bigcup_{i=1}^{s}\mathsf{Z}(\mathcal{T}_i/\ini(\mathcal{T}_i)) = \bigcup_{i=1}^{s}\zero(\sat(\mathcal{T}_i)),
% \end{equation}
% then we call the set $\{\mathcal{T}_1, \ldots, \mathcal{T}_s \}$ a \emph{triangular decomposition} of $\mathcal{F}$. In particular, if each $\mathcal{T}_i$ is a regular chain, normal chain, simple set, or irreducible triangular set, then we call $\{\mathcal{T}_1, \ldots, \mathcal{T}_s \}$ a regular, normal, simple, or irreducible triangular decomposition, respectively.
% \end{definition}

%###################################

\subsection{Triangular decompositions of semi-algebraic systems}

% XXXXXXXXXXXXXXXXXXXXXXXXXXXXX

In this section, we fix the base field $\K$ as
the rational number field $\Q$,
and recall related results on triangular decompositions
of semi-algebraic systems introduced in~\cite{JSCChen:2013}.

\begin{definition}
\label{def:spec}
A  regular system $[\caT,\caH]$ of $\Q[\x]$, is \textit{squarefree} if for all $t\in \caT$, the derivative of $t$ w.r.t. $\mvar(t)$, denoted by ${\rm der}(t)$, is regular modulo ${\rm sat}(\caT)$. Given $u\in\mathbb{R}^d$, we say that a (squarefree) regular system $[\caT,\caH]$ \textit{specializes well} at $u$ if $h_T(u)\neq0$ and $[\caT(u),\caH(u)]$ is a (squarefree) regular system.
% A regular chain is called $d$-dimensional if it has exactly $d$ free variables, that is, if its saturated ideal has dimension $d$.
\end{definition}

\begin{definition}[Semi-algebraic system]
Given finite polynomial sets $F=\{f_1,\ldots,f_s\}$, $N=\{n_1,\ldots,n_k\}$, $P=\{p_1,\ldots,p_e\}$, and $H=\{h_1,\ldots,h_{\ell}\}$ of $\mathbb{Q}[\x]$,
let  $F_{=}$ denote the set of equations $f_1=0,\ldots,f_s=0$,
$N_{\geq}$ denote the set of nonnegative inequalities $\{n_1\geq0,\ldots,n_k\geq0\}$, $P_>$ denote the set of positive inequalities $\{p_1>0,\ldots,p_e>0\}$, and $H_{\neq}$ denote the set of inequations $\{h_1\neq0,\ldots,h_{\ell}\neq0\}$. We denote by $S=[F_{=},N_{\geq},P_>,H_{\neq}]$ the \textit{semi-algebraic system} (SAS) defined as the conjunction of the constraints in $F_{=}$, $N_{\geq}$, $P_>$, and $H_{\neq}$.
\end{definition}

% \begin{definition}[Border polynomial]
% Let $\mathfrak{R}$ be either a squarefree regular chain $T$, or a squarefree regular system $[T,P]$ in $\mathbb{Q}[\boldsymbol{x}]$. We denote by $B_{{\rm sep}}(T)$, $B_{{\rm ini}}(T)$, $B_{{\rm ineqs}}([T,P])$ the set of irreducible factors of $\prod_{t\in T}{\rm res}({\rm discrim}(t,{\rm mvar}(t)),T)$, $\prod_{t\in T}{\rm res}({\rm init}(t),T)$, and $\prod_{f\in P}{\rm res}(f,T)$, respectively. Denote by $BP(\mathfrak{R})$ the set $B_{{\rm sep}}(T)\cup B_{{\rm ini}}(T)\cup B_{{\rm ineqs}}([T,P])$. Then $BP(\mathfrak{R})$ (resp. the polynomial $\prod_{f\in BP(\mathfrak{R})}f$) is called the \textit{border polynomial set} (resp. \textit{border polynomial}) of $\mathfrak{R}$.
% \end{definition}

% More properties of border polynomials can be found in \cite{YHX:2001,YX:2005,JSCChen:2013}.

% \begin{definition}[Pre-regular semi-algebraic system]
% Let $[T,P]$ be a squarefree regular system of $\mathbb{Q}[\boldsymbol{u},\boldsymbol{y}]$. Let $bp$ be the border polynomial of $[T,P]$.
% Let $B\subseteq\mathbb{Q}[\boldsymbol{u}]$ be a polynomial set such that $bp$ divides the product of polynomials in $B$. We call the triple $[B_{\neq},T,P_>]$ a \textit{pre-regular semi-algebraic system} of $\mathbb{Q}[\boldsymbol{x}]$. Its zero set, writtrn as $\zero_{\mathbb{R}}(B_{\neq},T,P_>)$, is the set $(u,y)\in\mathbb{R}^n$ such that $b(u)\neq0$, $t(u,y)=0$, $p(u,y)>0$, for all $b\in B$, $\in T$, $p\in P$.
% \end{definition}

\begin{definition}[Regular semi-algebraic system]
Let $\caR=[\caT,\caP]$ be a squarefree regular system of $\Q[\x]$.
Denote by $\bfu=u_1,\ldots,u_d$ and $\bfy=y_1,\ldots,y_r$ $(r+d=n)$ respectively the free and the main variables of $\caT$.
Let $\mathcal{Q}$ be a quantifier-free formula over $\Q[\bfu]$,
which defines a semi-algebraic set $S$ of $\R^d$.
We say that $rsa:=[\caQ,\caT,{\cal P}_{>}]$ is a \textit{regular semi-algebraic system} if:
$(i)$ $\mathcal{Q}$ defines a nonempty open subset $S$ in $\mathbb{R}^d$;
$(ii)$ the regular system $[\caT,\caP]$ specializes well at every point $u$ of $S$;
$(iii)$ at each point $u$ of $S$, the specialized system $[\caT(u),P(u)_>]$ admits real solutions.
\end{definition}

\begin{theorem}\label{Theorem-rsas}
Let $\mathfrak{S}$ be a semi-algebraic system of $\mathbb{Q}[\boldsymbol{x}]$. Then one can compute a  triangular decomposition of $\mathfrak{S}$, that is, finitely many regular semi-algebraic systems such that the union of their zero sets is the zero set of $\mathfrak{S}$.
\end{theorem}

\begin{remark}
\label{remark:rsas}
An algorithm for computing such a triangular decomposition
is called {\sf RealTriangularize},
which has been implemented as part of the RegularChains library in \textsc{Maple}.
\end{remark}

\subsection{Hurwitz determinants and subresultant sequences}

%M_i=\overbrace{}^{p+q-i}

The terminology of subresultants here is from \cite{Loos1982}. The main results about the Hurwitz determinants and subresultant
sequences are quoted from \cite{Weber2000,ChenMM-2011}. Let $p\leq q$ be two positive integers and $M$ be a $p\times q$ matrix with coefficients in $\mathbb{R}$. Let $M_i$ be the square submatrix of $M$ consisting of the first $p-1$ columns of $M$ and the $i$-th column of $M$, for $i=p,\ldots,q$. Let ${\rm det}M_i$ be the determinant of $M_i$. We denote by ${\rm dpol}(M)$ the
element of $\mathbb{R}[y]$, called the \textit{determinant polynomial} of $M$, given by
\[{\rm det}M_py^{q-p}+{\rm det}M_{p+1}y^{q-p-1}+\cdots+{\rm det}M_q.\]

Let $f_1(y),\ldots,f_p(y)$ be a set of polynomials in $\mathbb{R}[y]$. Let
\[q=1+\max({\rm deg}f_1(y),\ldots,{\rm deg}f_p(y)).\]
The matrix $M$ of $f_1,\ldots,f_p$ is defined by $M_{ij}={\rm coeff}(f_i,y^{q-j})$.

We now introduce subresultants. Let $f=a_my^m+\cdots+a_0$, $g=b_ny^n+\cdots+b_0$ be two polynomials in $\mathbb{R}[y]$. For $i=0,1,\ldots,\min(m,n)-1$, the $i$-th \textit{subresultant} of $f$ and $g$ is defined as
\begin{equation}\label{subequ1}
\begin{split}
S_i(f,m,g,n,y):={\rm dpol}(\underbrace{y^{n-1-i}f,\ldots,yf,f}_{n-i},\underbrace{y^{m-1-i}g,\ldots,yg,g}_{m-i}).
\end{split}
\end{equation}
The defining matrix on the right-hand side of \eqref{subequ1} has $m+n-2i$ rows and $m+n-i$ columns. Hence $S_i(f,g,y)$ is a polynomial in $\mathbb{R}[y]$ with degree at most $i$. Let $s_i(f,g,y)={\rm coeff}(S_i(f,g,y),y^i)$ and call it the \textit{principle subresultant coefficient} of $S_i(f,g,y)$. Let $\delta=\min\{m,n\}$. We extend the above definitions of subresultants and principal subresultant coefficients in such a way that
$f$ and $g$ are themselves subresultants:
\begin{equation}\label{S-defi-1}
\left\{\begin{array}{l}
S_{\delta+1}=f,\quad S_{\delta}=g,\quad s_{\delta+1}=a_m,\quad s_{\delta}=b_n,\quad \text{when $m\geq n$}\\
S_{\delta+1}=g,\quad
S_{\delta}=f,\quad
s_{\delta+1}=b_n,\quad
s_{\delta}=a_m,\quad \text{when $m<n$}\\
\end{array}\right..
\end{equation}

In \eqref{subequ1} and \eqref{S-defi-1}, $m$ and $n$ are treated as \textit{formal degrees}, allowing degree drops ($\deg(f) < m$ with $a_m = 0$, or $\deg(g) < n$ with $b_n = 0$). In particular, when $g = 0$ with $f \neq 0$, or symmetrically when $f = 0$ with $g \neq 0$, the zero polynomial retains its prescribed formal degree with formal leading coefficient $b_n = 0$ or $a_m = 0$, respectively.

Let $f(y)=a_0y^n+a_1y^{n-1}+\cdots+a_n\in\mathbb{R}[y]$ be a polynomial of degree $n$. The $i$-th Hurwitz determinant $\Delta_i$, for $i=1,2\ldots,n$, are defined by the following determinants:
\begin{equation*}% \label{H-matrix}
\begin{split}
&\Delta_1=a_1,\quad \Delta_2=\left|
 \begin{matrix}
    a_1& a_3 \\
   a_0 & a_2 \\
  \end{matrix}
  \right|,\quad \Delta_3=\left|
 \begin{matrix}
    a_1& a_3 & a_5\\
   a_0 & a_2 & a_4 \\
   0 & a_1 & a_3 \\
  \end{matrix}
  \right|,\\
  &\Delta_j=\left|
 \begin{matrix}
    a_1& a_3 & \cdots & \cdots\\
   a_0 & a_2 & a_4 & \cdots\\
   0 & a_1 & a_3 & \cdots\\
   0 & a_0 &a_2 & \cdots\\
   \cdots & \cdots & \cdots & \cdots\\
   0 & 0 & \cdots & a_j\\
  \end{matrix}
  \right|,\ldots,j=4,\ldots,n.
\end{split}
\end{equation*}
We write $f(y)=f_1(y^2)+yf_2(y^2)$. If $n$ is even, $n=2m$, then
\[f_1(y)=a_0y^m+a_2y^{m-1}+\cdots+a_{2m},\quad f_2(y)=a_1y^{m-1}+a_3y^{m-2}+\cdots+a_{2m-1}.\]
On the other hand, if $n$ is odd, $n=2m+1$, then
\[f_1(y)=a_1y^m+a_3y^{m-1}+\cdots+a_{2m+1},\quad f_2(y)=a_0y^{m}+a_2y^{m-1}+\cdots+a_{2m}.\]

The next result is a refined version of \cite[Theorem 5]{ChenMM-2011} as well as a corrected version of \cite[Theorem 3.1]{Weber2000}, which demonstrates the close relationship between the Hurwitz determinants and principle subresultant sequences.

\begin{theorem}\label{HSres}
Let $f(y)=a_0y^n+a_1y^{n-1}+\cdots+a_n\in\mathbb{R}[y]$ be a polynomial of degree $n$, and we write $f(y)=f_1(y^2)+yf_2(y^2)$. Let $\Delta_{-1}=\Delta_0=a_0$ and $\{\Delta_1,\ldots,\Delta_n\}$ be the Hurwitz determinants sequence of $f$. Then the following conclusion holds:
\begin{itemize}
\item[(i)] If $n=2m+1$, we have
$\Delta_{n-1-2i}=(-1)^{\frac{(m-i)(m-i-1)}{2}}s_i(f_1,m,f_2,m,y)$ hold, for $i=0,1,\ldots,m$.

\item[(ii)] If $n=2m$, we have $\Delta_{n-1-2i}=(-1)^{\frac{(m-i)(m-i-1)}{2}}s_i(f_1,m,f_2,m-1,y)$ hold, for $i=0,1,\ldots,m$.

\item[(iii)] If $n=2m+1$, we have
$\Delta_{n-2i}=(-1)^{\frac{(m-i)(m-i+1)}{2}}s_i(f_1,m,yf_2,m+1,y)$ hold, for $i=0,1,\ldots,m$.

\item[(iv)] If $n=2m$, we have $\Delta_{n-2i}=(-1)^{\frac{(m-i)(m-i-1)}{2}}s_i(yf_2,m,f_1,m,y)$ hold, for $i=0,1,\ldots,m$.

\end{itemize}
\end{theorem}

\begin{remark}
The above statements have used the extension \eqref{S-defi-1}. Note that this result originally comes from \cite[Theorem 3.1]{Weber2000} but has sign errors, which were later fixed in \cite[Theorem 5]{ChenMM-2011}. However, the ranges of $i$ in previous studies cannot cover the case $i=m$ except for $(iii)$. This situation is crucial for deriving the complete algebraic conditions for Hopf bifurcation points (see Theorem \ref{Hopf-th3}). Here we extend the Hurwitz determinants by defining $\Delta_{-1}=\Delta_0=a_0$. In such a way, the ranges of $i$ in Theorem \ref{HSres} can be extended to $i=m$. In fact, such extension is important to fully take advantage the power of subresultants, which would be clear in the proof of Theorem \ref{Hopf-th3}. It is worth noting that we switch the position of $f_1$ and $yf_2$ in $(iv)$ since we agree that $\Delta_0=a_0$ and the extension \eqref{S-defi-1}.
\end{remark}

The following two corollaries are direct consequences of Theorem \ref{HSres}.
%{\color{red} CC: specify the exact range of $k$ in the following two corollaries. The range can be pushed to the limit if we extend the definition of $\Delta_0$.}

\begin{corollary}\label{lea-HS1}
Let the polynomial $f$ be as above, where $a_n\neq0$. Then $f$ has a nonzero root pair
$\{z,-z\}$ in $\mathbb{C}$ if and only if $z^2$ is a common zero of $f_1(y)$ and $f_2(y)$. Moreover, $f$ has exactly $k$ pairs of symmetric roots $\{z_i,-z_i\}$ in $\mathbb{C}$ (counted with multiplicity) if and only if
\[
 \Delta_{n-1-2j}=0\quad(0\le j<k),\qquad \Delta_{n-2k-1}\ne0,
 \qquad 0\le k\le\lfloor n/2\rfloor.
\]
\end{corollary}

\begin{proof}
Since $a_n\neq 0$, $f(x)$ has no zero roots. Thus $f$ has a nonzero root pair
$\{z,-z\}$ in $\mathbb{C}$ if and only if $z^2$ is a common zero of $f_1(y)$ and $f_2(y)$. For $t=z^2\ne0$, $(\lambda^2-t)^j\mid f(\lambda)$ if and only if $(y-t)^j$ divides both $f_1(y)$ and $f_2(y)$, by taking even and odd parts. The largest such $j$ is the multiplicity of the symmetric pair $\{z,-z\}$. Therefore, $f$ has exactly $k$ pairs of symmetric roots $\{z_i,-z_i\}$ if and only if $f_1(y)$ and $f_2(y)$ has $k$ common roots counted with multiplicity. Then the result follows from the above Theorem \ref{HSres}, which builds the relation between Hurwitz determinants and principle subresultant coefficients, and  \cite[Theorem 5]{DBLP:journals/jsc/ChenM12},
which establishes the relation between the GCD of $f_1$ and $f_2$ and the subresultants covering the corner case $k=\lfloor \frac{n}{2} \rfloor$.
\end{proof}

\begin{corollary}\label{lea-HS2}
Notation as above, if $a_n\neq0$, then $f$ has a nonzero root pair
$\{z,-z\}$ in $\mathbb{C}$ if and only if $z^2$ is a common zero of $f_1(y)$ and $yf_2(y)$. Moreover, $f$ has exactly $k$ pairs of symmetric roots $\{z_i,-z_i\}$ in $\mathbb{C}$ (counted with multiplicity) if and only if
\[
 \Delta_{n-2j}=0\quad(0\le j<k),\qquad \Delta_{n-2k}\ne0,
 \qquad 0\le k\le\lfloor n/2\rfloor.
\]

\end{corollary}

\begin{remark}\label{rem-cor-HS2}
Note that when $n=2m$, neither $f_2(y)$ nor $yf_2(y)$ can have $m$ pairs of nonzero symmetric roots $\{z_i,-z_i\}$ since $f_2(y)$ is a polynomial of degree $m-1$. As a result, to make $f(y)$ have $m$ pairs of symmetric roots, $f_2(y)$ and $yf_2(y)$ must be identically zero. This situation was not considered in \cite{Weber2000,KWeber-2002,ChenMM-2011}, and is covered now by our extension in Theorem \ref{HSres}. A conclusion about this fact is placed in Corollary \ref{coro-th-hp3}.

\end{remark}

% The following result shows the behavior of the $\Delta_i$'s when we add to the roots of the polynomial $f$ some pairs of symmetric points $z_j$ and $-z_j$ of the plane. Its proof can be found in \cite{Weber2000}.

% \begin{lemma}\label{lem-HS2}
% Let $f(y)\in\mathbb{R}[y]$ be a polynomial of degree $n$ and $z_1,\ldots,z_k$ be arbitrary
% complex numbers. Let
% \[f^*(y)=f(y)(y^2-z_1^2)\cdots(y^2-z_k^2).\]
% If $\Delta_i^*$ is the Hurwitz determinants of order $i$ of the polynomial $f^*(y)$, then $\Delta_i=\Delta_i^*$, for $i=1,\ldots,n$. Similarly, let $f^*(y)=f(y)y^k$, then we also have $\Delta_i=\Delta_i^*$ hold.
% \end{lemma}

\section{Invariant Algebraic Surfaces via Elimination Theory}
\label{sec:dsys}

Computational algebra methods have been widely used in the qualitative analysis of differential equations, see \citep{chr1999lio,chr2000alg,wan2005sta,wan2005dif}. In this section, we review algorithms based on elimination theory for computing invariant algebraic surfaces of parametric dynamical systems. We first introduce the Gr\"obner basis method via elimination ideal computation, following the methodological framework established by Antonov, Aybar and co-authors \citep{adrt2016,ayb2015inv}. As shown in subsequent subsections, this approach yields only necessary conditions for the existence of invariant algebraic surfaces, whose sufficiency requires separate verification.

To circumvent this limitation, triangular decomposition offers a more comprehensive elimination alternative: it accepts nonvanishing conditions on the parameters directly as part of the input, and it can either complement or be used independently of Gr\"obner basis methods. However, the existing literature lacks algorithms for this problem that work directly with triangular sets. To address this gap, we propose an approach based on triangular decomposition for this task. To our knowledge, this is the first application of triangular set methods to directly compute invariant algebraic surfaces for parametric dynamical systems. Crucially, we prove that our method yields both necessary and sufficient conditions for the existence of such surfaces. Furthermore, following an explicit constructible-set postprocessing step on each branch, the approach can recover the nontrivial invariant surfaces together with their exact parameter conditions.

Consider a parametric polynomial differential system of the form \eqref{L-dot_xx}. Let
\begin{equation}\label{equation:X}
\mathcal{X}=F_1\frac{\partial}{\partial x_1}+\cdots+F_n\frac{\partial}{\partial x_n}
\end{equation}
denote the polynomial vector field associated with system \eqref{L-dot_xx}. A polynomial $L\in \mathbb{K}[\boldsymbol{x},\boldsymbol{\mu}]$ is called a \emph{Darboux polynomial} of system \eqref{L-dot_xx} if there exists a polynomial $K \in \mathbb{K}[\boldsymbol{x},\boldsymbol{\mu}]$ such that
\begin{equation}\label{equation:dsys-dlkl}
\mathcal{X}(L)=KL,
\end{equation}
where $K$ is referred to as the \emph{cofactor} of $L$. Clearly, if $L$ is a Darboux polynomial of system \eqref{L-dot_xx}, then for each fixed $\boldsymbol{\mu}$, the algebraic variety $\boldsymbol{V}(L)=\{\boldsymbol{x}\in\mathbb{K}^n\mid L(\boldsymbol{x},\boldsymbol{\mu})=0\}$ is \emph{invariant} under the flow of system \eqref{L-dot_xx}; that is, any orbit of system \eqref{L-dot_xx} starting in $\boldsymbol{V}(L)$ remains in $\boldsymbol{V}(L)$ for all time within its maximal interval of existence.

\begin{remark}\label{rem-Ias-def}
For a real polynomial vector field $\mathcal{X}$, its Darboux polynomials and corresponding cofactors (if any) may be complex. In this case, they appear in complex conjugate pairs. If $L$ is a Darboux polynomial of system \eqref{L-dot_xx}, then $L=0$ is called an \emph{invariant algebraic curve} for $n=2$, an \emph{invariant algebraic surface} for $n=3$, and an \emph{invariant algebraic hypersurface} for $n>3$. Throughout the remainder of this paper, when the dimension $n$ of the system is unspecified, we simply refer to $L=0$ as an invariant algebraic surface regardless of the dimension.
\end{remark}

Recall that a \emph{rational first integral} of system \eqref{L-dot_xx} is a nonconstant rational function $H=\frac{F}{G}\in\mathbb{K}(\boldsymbol{x},\boldsymbol{\mu})\setminus\mathbb{K}(\boldsymbol{\mu})$ (with $F,G\in\mathbb{K}[\boldsymbol{x},\boldsymbol{\mu}]$ and $G\neq0$) satisfying $\mathcal{X}(H)=0$. Darboux polynomials are closely linked to rational first integrals; in particular, the following theorem establishes that the existence of sufficiently many Darboux polynomials guarantees a rational first integral.

\begin{theorem}[Darboux--Jouanolou]\label{D-J-integrability}
Assume that the $n$-dimensional polynomial differential system \eqref{L-dot_xx} of degree $d$ has $p$ irreducible Darboux polynomials. Set
\[N=\binom{d+n-1}{n}.\]
\begin{itemize}
    \item[(a)] If $p\geq N+1$, system \eqref{L-dot_xx} has a Darboux first integral.
    \item[(b)] System \eqref{L-dot_xx} has a rational first integral if and only if $p\geq N+n$.
\end{itemize}
\end{theorem}

Statement (a) of Theorem \ref{D-J-integrability} was originally established by Darboux \cite{dar1878meo} in 1878, whereas statement (b) was proved by Jouanolou \cite{Jou1979} in 1979 using sophisticated tools from algebraic geometry (see \cite{CCJL2000,JLZX2012} for an elementary proof of Jouanolou's result).

In light of Theorem \ref{D-J-integrability}, finding Darboux polynomials (or equivalently, invariant algebraic surfaces) is a central step in constructing rational first integrals. Motivated by this, we focus in this paper on the problem of determining whether the parametric differential system \eqref{L-dot_xx} admits an invariant algebraic surface $L=0$ associated with a cofactor $K$ of prescribed forms:
\begin{equation}\label{eq:preL}
L = \sum_{i=0}^{\ell_1} f_{i} \boldsymbol{x}^{\boldsymbol{p}_i}=0, \quad K = \sum_{j=0}^{\ell_2} a_{j} \boldsymbol{x}^{\boldsymbol{q}_j} ,
\end{equation}
where $\boldsymbol{p}_0=\boldsymbol{0}$, so that $f_0$ is the constant term of $L$; the monomial sets
$\{\boldsymbol{x}^{\boldsymbol{p}_i}\}_{i=0}^{\ell_1}$ and
$\{\boldsymbol{x}^{\boldsymbol{q}_j}\}_{j=0}^{\ell_2}$
are prescribed, and all $f_i$ and $a_j$ are undetermined coefficients.

Recall that for a polynomial $P \in \mathbb{K}[\boldsymbol{x}]$, $\operatorname{term}(P)$ denotes the set of all terms in $P$. For a term $\boldsymbol{t}$ of $P$, $\operatorname{coef}(\boldsymbol{t})$ denotes the coefficient of $\boldsymbol{t}$ in $P$.
If equation~\eqref{equation:dsys-dlkl} holds,
then $\mathcal{X}L - KL \equiv 0$ may be regarded as a polynomial identity in the variables $\boldsymbol{x}$.
This identity is equivalent to the following system of equations:
\begin{equation}\label{equation:dsys-termeq}
\{ \operatorname{coef}(\mathcal{X}L - KL, \boldsymbol{t}) = 0 \;|\;  \boldsymbol{t} \in \operatorname{term}(\mathcal{X}L - KL)\},
\end{equation}
where $\mathcal{X}L-KL$ is treated as a polynomial in $\boldsymbol{x}$ with coefficients in $\mathbb{K}[\boldsymbol{\mu},\boldsymbol{f},\boldsymbol{a}]$, i.e., as an element of the polynomial ring $\mathbb{K}[\boldsymbol{\mu}, \boldsymbol{f}, \boldsymbol{a}][\boldsymbol{x}]$, where $\boldsymbol{f} = (f_0,\ldots,f_{\ell_1})$ and $\boldsymbol{a} = (a_0,\ldots,a_{\ell_2})$. Based on the above observation, we establish the following lemma.

\begin{lemma}
\label{lemma:dsys-equivalence}
For fixed values of $(\boldsymbol{\mu},\boldsymbol{f},\boldsymbol{a})$, the polynomials $L$ and $K$ of the prescribed forms \eqref{eq:preL} satisfy $\mathcal{X}L=KL$ if and only if the coefficient equations \eqref{equation:dsys-termeq} hold. In particular, for a fixed $\overline{\boldsymbol{\mu}}$, a nonconstant $L$ with cofactor $K$ exists if and only if there are coefficients $(\boldsymbol{f},\boldsymbol{a})\in\overline{\mathbb{K}}^{\ell_1+\ell_2+2}$ satisfying \eqref{equation:dsys-termeq} at $\overline{\boldsymbol{\mu}}$ and $(f_1,\ldots,f_{\ell_1})\neq\boldsymbol{0}$.
\end{lemma}
Unless otherwise stated, the coefficient witnesses $\boldsymbol{f},\boldsymbol{a}$ in problem ($\mathcal{P}_1$) range over $\overline{\mathbb{K}}$. Thus, for a real parameter point, membership in the constructible existence locus guarantees a witness over the algebraic closure, not necessarily a real invariant surface.

The polynomial set
$\mathcal{P} := \{\operatorname{coef}(\mathcal{X}L - KL, \boldsymbol{t}) \mid \boldsymbol{t} \in \operatorname{term}(\mathcal{X}L - KL)\} \subseteq \mathbb{K}[\boldsymbol{\mu},\boldsymbol{f},\boldsymbol{a}]$
consists entirely of polynomials in the parameters, including system parameters $\boldsymbol{\mu}$ and undetermined coefficients $\boldsymbol{f}, \boldsymbol{a}$. Consequently, we may analyze $\mathcal{P}$ in the parameter ring $\mathbb{K}[\boldsymbol{\mu},\boldsymbol{f},\boldsymbol{a}]$ to investigate the existence of invariant algebraic surfaces. We call $\mathcal{P}$ the \emph{associated parametric polynomial set} of differential system \eqref{L-dot_xx} with respect to the ansatz \eqref{eq:preL}, or simply the parametric polynomial set. Once $\ell_1$ and $\ell_2$ are fixed, the set $\mathcal{P}$ is finite.

In the following, we present general theoretical results and complete computational methods based on Gr\"obner basis computation and triangular decomposition.

\subsection{The Gr\"obner basis approach}
\label{sec:GBmethod}

Lemma \ref{lemma:dsys-equivalence} provides necessary and sufficient conditions for determining whether a general $n$-dimensional polynomial differential system admits a prescribed invariant algebraic surface of a bounded degree.
Let $\mathcal{P} \subseteq \mathbb{K}[\boldsymbol{\mu},\boldsymbol{f},\boldsymbol{a}]$ be the parametric polynomial set of system \eqref{L-dot_xx} with respect to the ansatz \eqref{eq:preL}. To determine the parameter conditions under which this criterion is satisfied, we must investigate the projection of the zero set of $\mathcal{P}$ onto the system parameter space. A feasible approach using elimination theory consists in computing the elimination ideal
$\langle \mathcal{P} \rangle \cap \mathbb{K}[\boldsymbol{\mu}].$
According to Theorem \ref{theorem:poly-alge-elim},
we may obtain the $(\ell_1+\ell_2+2)$-th elimination ideal of $\langle \mathcal{P} \rangle$ via Gr\"obner basis computation. Under the variable ordering $\boldsymbol{\mu}<\{\boldsymbol{f},\boldsymbol{a}\}$, this elimination ideal is precisely $\langle \mathcal{P} \rangle \cap \mathbb{K}[\boldsymbol{\mu}].$

Note that equation \eqref{equation:dsys-dlkl} always admits the solutions with constant $L$, i.e., with $f_1=\cdots=f_{\ell_1}=0$ (then $K=0$ if $f_0\neq0$, and $K$ is arbitrary if $f_0=0$), for all system parameter values, which implies that the elimination ideal mentioned above is always the zero ideal.
A constant $L$ does not define an algebraic surface: its zero set is empty if $f_0\neq0$ and is the whole space if $f_0=0$. Hence the solutions with constant $L$ must be excluded before the elimination ideal is computed; accordingly, an invariant algebraic surface $L=0$ is called \emph{nontrivial} if $L$ is nonconstant.

To exclude these solutions, we work with the saturated ideal
$\mathfrak{I} = \langle \mathcal{P} \rangle : \langle f_1,\ldots,f_{\ell_1} \rangle^{\infty}$.
Let
\[
cs:=\zero(\mathcal{P})\setminus\zero(f_1,\ldots,f_{\ell_1})
\]
be the constructible set of nontrivial coefficient solutions.
By Proposition~\ref{proposition:sat-zeros} below, the zero set of $\mathfrak{I}$ is the Zariski closure of $cs$,
\[\zero(\mathfrak{I}) = \overline{cs}.\]
Let $m$ be the number of system parameters $\boldsymbol{\mu}$. The set
\[S := \operatorname{Proj}_{m}(cs)\]
of system parameter values for which system \eqref{equation:dsys-termeq}
admits a nontrivial solution is called the \emph{existence locus} of the
ansatz \eqref{eq:preL}.

% It should be noted that saturation does not guarantee that the resulting elimination ideal is nonzero.
\begin{remark}
\label{remark:sat-nonzero}
The elimination ideal of $\mathfrak{I}$ being \emph{nonzero} cannot be
proved a priori from the saturation step alone: it is the zero ideal if and
only if the Zariski closure of the existence locus is the entire parameter
space. Whether this
degenerate situation occurs is part of the classification itself; for a
concrete system it is certified by the computation (see
Section~\ref{sec:result-1}, where the elimination ideal of the BZ system
turns out to be generated by twelve explicit nonzero polynomials).
\end{remark}

\begin{proposition}[{
\cite[Theorem 10, Section 4.4]{cox2013ide}}]
\label{proposition:sat-zeros}
Let $\mathfrak{I}$ and $\mathfrak{J}$ be ideals in $\mathbb{K}[\boldsymbol{x}]$,
then
$\zero(\mathfrak{I}:\mathfrak{J}^{\infty}) = \overline{\zero(\mathfrak{I}) \setminus \zero(\mathfrak{J})}$,
i.e., it coincides with the Zariski closure of the difference between $\zero(\mathfrak{I})$ and $\zero(\mathfrak{J})$.
\end{proposition}

The following result provides necessary conditions for a differential system of the form \eqref{L-dot_xx} to admit an invariant algebraic surface of bounded degree.

% \begin{proposition}\label{proposition:nontrivial-condition}
% Let $\mathcal{F}$ denote the associated parametric polynomial set of system \eqref{L-dot_xx} with respect to the prescribed form \eqref{eq:preL}, and let $\mathfrak{I} = \langle \mathcal{F} \rangle : \langle f_1,\ldots,f_m \rangle^{\infty}$. If system \eqref{L-dot_xx} admits a nontrivial invariant algebraic surface for some system parameter value $\overline{\boldsymbol{\mu}}$, then
% \begin{equation*}
% \overline{\boldsymbol{\mu}} \in \zero(\frI_{m+h+2}).
% \end{equation*}
% That is, the parameter value must lie in the zero set of the $(m+h+2)$-th elimination ideal of $\mathfrak{I}$.
% \end{proposition}

\begin{proposition}\label{proposition:nontrivial-condition}
Let $\mathcal{P}$ denote the associated parametric polynomial set of system \eqref{L-dot_xx} with respect to the prescribed form \eqref{eq:preL}. With $cs$, $m$, and $S$ as above, set
\[
\mathfrak{I}
=
\langle\mathcal{P}\rangle:
\langle f_1,\ldots,f_{\ell_1}\rangle^{\infty},
\qquad
\mathfrak{I}_{\ell_1+\ell_2+2}
=
\mathfrak{I}\cap\mathbb{K}[\boldsymbol{\mu}].
\]
Then
\[
\zero(\mathfrak{I})=\overline{cs},
\qquad
\zero(\mathfrak{I}_{\ell_1+\ell_2+2})=\overline{S}.
\]
In particular, if system \eqref{L-dot_xx} admits a nontrivial invariant algebraic surface for a parameter value $\overline{\boldsymbol{\mu}}$, then $\overline{\boldsymbol{\mu}}\in\zero(\mathfrak{I}_{\ell_1+\ell_2+2})$. Conversely, every point of $\zero(\mathfrak{I}_{\ell_1+\ell_2+2})$ lies in the Zariski closure of the existence locus, although a boundary point of this closure need not itself lift to a nontrivial solution of the prescribed ansatz.
\end{proposition}

\begin{proof}
Proposition~\ref{proposition:sat-zeros} gives
$\zero(\mathfrak{I})=\overline{cs}$.
By part~(i) of the Closure Theorem in elimination theory, applied over $\overline{\mathbb{K}}$
(see, e.g., \cite[Chapter~3, \S 2, Theorem~3]{cox2013ide}),
\[
\zero(\mathfrak{I}_{\ell_1+\ell_2+2})
=
\overline{\operatorname{Proj}_{m}(\zero(\mathfrak{I}))}.
\]
Since the coordinate projection onto the system-parameter space is Zariski continuous (the preimage of a polynomial zero set is again a polynomial zero set), we have
\[
\overline{\operatorname{Proj}_{m}(\zero(\mathfrak{I}))}
=
\overline{\operatorname{Proj}_{m}(\overline{cs})}
=
\overline{\operatorname{Proj}_{m}(cs)}
=
\overline{S}.
\]
Hence $\zero(\mathfrak{I}_{\ell_1+\ell_2+2})=\overline{S}$.
\end{proof}

% \begin{proof}
% By Lemma \ref{lemma:dsys-equivalence} and the nontriviality requirement that $L$ be nonconstant,
% for the system parameter value $\overline{\boldsymbol{\mu}}$,
% there exist undetermined coefficients
% $\boldsymbol{f},\boldsymbol{a}$ with $(f_1,\ldots,f_m)\neq\boldsymbol{0}$ such that
% $(\overline{\boldsymbol{\mu}},\boldsymbol{f},\boldsymbol{a})\in \zero(\mathcal{F})\setminus\zero(f_1,\ldots,f_m)\subseteq\zero(\mathfrak{I})$. Since every polynomial $F \in \mathfrak{I}_{m+h+2}$ belongs to $\mathfrak{I}$ and $(\overline{\boldsymbol{\mu}},\boldsymbol{f},\boldsymbol{a})\in \zero(\mathfrak{I})$,
% we have $F(\overline{\boldsymbol{\mu}}) = 0$ for all $F \in \mathfrak{I}_{m+h+2}$. Hence $\overline{\boldsymbol{\mu}} \in \zero(\mathfrak{I}_{m+h+2})$.
% \end{proof}

% The above proposition demonstrates that the parameter conditions derived by computing elimination ideals via Gr\"obner bases are necessary for the system to admit invariant algebraic surfaces of the specified form.
To determine whether these conditions are sufficient, the following two approaches may be adopted:
\begin{enumerate}[label=\textup{(\arabic*)}, ref=\arabic*]
    \item \textit{Theoretical approach:} Apply the Extension Theorem from elimination theory to prove that, under appropriate nondegeneracy conditions, parameter values in the zero set of the elimination ideal may indeed be extended to solutions of the original parametric polynomial set, thereby ensuring the existence of the corresponding invariant algebraic surfaces.
     \item \textit{Constructive approach:} For each computed parameter condition, impose the condition on the original parametric polynomial system and solve the resulting restricted system for the remaining undetermined coefficients $\boldsymbol{f}$ and $\boldsymbol{a}$. Then verify that the resulting polynomial $L$ satisfies equation \eqref{equation:dsys-dlkl} and is nonconstant.
\end{enumerate}

\begin{algorithm}
\renewcommand{\algorithmicrequire}{\textbf{Input:}}
\renewcommand{\algorithmicensure}{\textbf{Output:}}
\caption{Computing Parameter Conditions via Gr\"obner Basis Elimination
\\ $\Psi:={\rm InvSurGB}((F_1,\ldots,F_n),L,K)$}
\label{algorithm:invsurgb}
\begin{algorithmic}[1]
\REQUIRE
{$\left\{
        \begin{array}{l}
        (F_1,\ldots,F_n),\; \text{vector field of dynamical system } \dot{\boldsymbol{x}} = (F_1,\ldots,F_n); \\
        L = \sum\limits_{i=0}^{\ell_1} f_i \boldsymbol{x}^{\boldsymbol{p}_i},\; \text{prescribed invariant algebraic surface};\\
        K = \sum\limits_{j=0}^{\ell_2} a_j \boldsymbol{x}^{\boldsymbol{q}_j},\; \text{prescribed cofactor}.\\
        \end{array}
        \right.$}

\ENSURE {$\Psi = \big\{ \mathfrak{p}_1, \ldots, \mathfrak{p}_s \big\}$, the minimal prime decomposition of $\sqrt{\frI_{\ell_1+\ell_2+2}}$ in $\mathbb{K}[\boldsymbol{\mu}]$; equivalently, the varieties $\zero(\mathfrak{p}_1),\ldots,\zero(\mathfrak{p}_s)$ are the irreducible components of the Zariski closure of the existence locus for which a nontrivial invariant algebraic surface of the prescribed form exists}
\STATE Set $\Psi := \emptyset$  and compute the associated parametric polynomial set $\mathcal{P}$ of vector field $(F_1,\ldots,F_n)$\;
\STATE $\frI := \ide{\mathcal{P}} : \ide{f_1,\ldots,f_{\ell_1}}^{\infty}$ \label{line:invsurgbl2}\;
\STATE Compute the $(\ell_1+\ell_2+2)$-th elimination ideal $\frI_{\ell_1+\ell_2+2}$ of $\frI$\;
\STATE Compute the minimal prime decomposition $\mathfrak{p}_1,\ldots,\mathfrak{p}_s$ of $\sqrt{\frI_{\ell_1+\ell_2+2}}$ \label{line:invsurgbl4}\;
\FOR{$i = 1$ \TO $s$}
\STATE $\Psi := \Psi \cup \{ \mathfrak{p}_i \}$\;
\ENDFOR
\RETURN $\Psi$\;
\end{algorithmic}
\end{algorithm}

\begin{remark}
The saturation step has a natural interpretation in terms of primary decomposition. Let $\mathfrak{J}=\bigcap_r\mathfrak{Q}_r$ be a primary decomposition and set $\mathfrak{I}_f=\langle f_1,\ldots,f_{\ell_1}\rangle$. Then
\[
\mathfrak{J}:\mathfrak{I}_f^{\infty}
=
\bigcap_{\mathfrak{I}_f\nsubseteq\sqrt{\mathfrak{Q}_r}}\mathfrak{Q}_r;
\]
see \cite[Chapter~4, \S4 and \S8]{cox2013ide}. Thus, saturation discards precisely those primary components whose associated primes contain $\mathfrak{I}_f$, or geometrically, those components entirely contained in $\zero(f_1,\ldots,f_{\ell_1})$. In the present problem, these are the components on which all the undetermined coefficients of $L$ vanish identically, so that $L$ is constant. Trivial points lying inside a retained mixed component are not removed by this step and must be treated separately when pointwise parameter conditions are required.

Although we compute the saturation directly, the underlying idea follows \cite{adrt2016,ayb2015inv}, where the Rabinowitsch trick is used to treat the conditions $f_i\neq0$ separately for $i=1,\ldots,\ell_1$. Collectively, their computations amount to saturation with respect to $\langle f_1,\ldots,f_{\ell_1}\rangle$.
\end{remark}

When addressing relatively complex systems or cases involving numerous parameters, the Gr\"obner basis method may encounter computational bottlenecks during elimination ideal and saturation calculations. Triangular decomposition provides an alternative: it represents the solution set branch by branch and accepts nonvanishing conditions as part of the input. These features are particularly useful for the present problem, where trivial branches and prescribed parameter constraints must be excluded.

\subsection{The triangular decomposition method}
\label{sec:trisysmethod}

As analyzed above, determining invariant algebraic surfaces of a parametric differential system \eqref{L-dot_xx} requires not only eliminating variables, but also preserving certain nonvanishing conditions on the parameters. Consequently, triangular decomposition of polynomial systems is particularly suited to this task.

For $\boldsymbol{r} = (r_1, \ldots, r_n) \in \overline{\mathbb{K}}^n$ and $i = 1,\ldots,n$,
recall that $\operatorname{Proj}_i(\boldsymbol{r}) = (r_1,\ldots,r_i)$ denotes the projection of $\boldsymbol{r}$ onto the first $i$ coordinates.
Let $\mathcal{F}, \mathcal{H} $ be two polynomial sets and  $\mathcal{T}$ be a triangular set in $\mathbb{K}[\boldsymbol{x}]$.
Recall that $\operatorname{Proj}_i \zero(\mathcal{F}/\mathcal{H})$ denotes the projection of the zero set $\zero(\mathcal{F}/\mathcal{H})$ onto the first $i$ coordinates and $\mathcal{T}^{(\le i)}$ denotes the truncated triangular set $\mathcal{T} \cap \mathbb{K}[x_1,\ldots,x_{i}]$.

\begin{proposition}
\label{proposition:dsys-trimethod}
Let $\mathcal{P}$ denote the parametric polynomial set of system \eqref{L-dot_xx} with respect to the prescribed form \eqref{eq:preL}
and let $cs := \zero(\mathcal{P})\setminus\zero(f_1,\ldots,f_{\ell_1})$.
Let $m$ be the number of system parameters $\boldsymbol{\mu}$.
Then system \eqref{L-dot_xx} admits a nontrivial invariant algebraic surface of the form \eqref{eq:preL} for some system parameter value $\overline{\boldsymbol{\mu}} \in \overline{\mathbb{K}}^{m}$ if and only if
$\overline{\boldsymbol{\mu}} \in \operatorname{Proj}_{m}(cs)$.
In particular, let $\{[\mathcal{T}_1, \mathcal{H}_1], \ldots, [\mathcal{T}_s, \mathcal{H}_s]\}$ be a finite set of
strong regular systems such that
$\zero(\mathcal{P})\setminus\zero(f_1,\ldots,f_{\ell_1}) = \bigcup_{j=1}^s \zero(\mathcal{T}_j/\mathcal{H}_j)$.
Then we have
\begin{equation}\label{equation:dsys-tri-existence}
\operatorname{Proj}_{m}(cs) = \bigcup_{j=1}^s\zero\big(\mathcal{T}_j^{(\le m)}/\mathcal{H}_j^{(\le m)}\big).
\end{equation}
\end{proposition}

\begin{proof}
Recall that the variable ordering is chosen such that the system parameters $\boldsymbol{\mu}$ are smaller than the undetermined coefficients $\boldsymbol{f}$ and $\boldsymbol{a}$.
By Lemma \ref{lemma:dsys-equivalence} and the nontriviality requirement that $L$ be nonconstant, the system admits a nontrivial invariant algebraic surface for the parameter value $\overline{\boldsymbol{\mu}}$ if and only if there exists a point $(\overline{\boldsymbol{\mu}}, \boldsymbol{f}, \boldsymbol{a}) \in \zero(\mathcal{P})$ with $(f_1,\ldots,f_{\ell_1})\neq\boldsymbol{0}$.
The claim follows from Remark~\ref{remark:cs} and the fact that  every strong regular system possesses the projection property \eqref{equation:projection}.

\end{proof}

\begin{remark}
By Remark \ref{remark:cs}, $\operatorname{Proj}_m(cs)$ can also be computed by decomposing $cs$ firstly into the zero sets of (weak) regular systems, and then refining this decomposition into a pre-comprehensive triangular decomposition. Alternatively, in practice, one could also firstly compute a triangular decomposition of $\mathcal{P}$ into (strong or weak) regular systems $R_1,\ldots,R_t$.
As $R := [\{f_1,\ldots,f_{\ell_1}\},\emptyset]$ is trivially a regular system with $\zero(R)=\zero(f_1,\ldots,f_{\ell_1})$,
one could compute $\operatorname{Proj}_{m}(\zero(R_i)\setminus\zero(R))$
instead to obtain the conditions on the parameters. This approach is formalized as Algorithm \ref{algorithm:invsurtri}. In our experiments, the difference and projection of constructible sets were performed with the \textsf{RegularChains} library in \textsc{Maple} \cite{ChenMorenoMazaPanXie2008}.

\end{remark}

\begin{remark}
In practice, the {TriSer} algorithm in \textsc{Epsilon} may also be used to decompose a polynomial system into triangular systems. When its output branches satisfy the required projection property (in particular, if they are strong regular systems), the same projection argument applies. Since {TriSer} avoids some of the refinement steps needed to enforce strong regularity, it often returns substantially fewer branches than {RegSer}. Branches without the required projection property can still be handled by a general constructible-set projection. As before, the trivial part is removed before projection. For the linear ansatz \eqref{equation:dsys-BZ-Ldeg1} of the BZ system, all $176$ {TriSer} branches passed the projection check, whereas {RegSer} returned $632$ branches.
\end{remark}

\begin{algorithm}
\renewcommand{\algorithmicrequire}{\textbf{Input:}}
\renewcommand{\algorithmicensure}{\textbf{Output:}}
\caption{Computing Parameter Conditions via Triangular Decomposition\\
$\Psi:={\rm InvSurTri}((F_1,\ldots,F_n),L,K)$}
\label{algorithm:invsurtri}
\begin{algorithmic}[1]
\REQUIRE {
        $\left\{
        \begin{array}{l}
        (F_1,\ldots,F_n),\; \text{vector field of dynamical system } \dot{\boldsymbol{x}} = (F_1,\ldots,F_n); \\
        L = \sum\limits_{i=0}^{\ell_1} f_i \boldsymbol{x}^{\boldsymbol{p}_i},\; \text{prescribed invariant algebraic surface};\\
        K = \sum\limits_{j=0}^{\ell_2} a_j \boldsymbol{x}^{\boldsymbol{q}_j},\; \text{prescribed cofactor}.\\
        \end{array}
        \right.$
}
\ENSURE
{A finite collection $\Psi=\{\mathcal{Q}_1,\ldots,\mathcal{Q}_t\}$ of constructible subsets of the parameter space such that the system admits a nontrivial invariant algebraic surface of the prescribed form~\eqref{eq:preL} for $\overline{\boldsymbol{\mu}}$ if and only if $\overline{\boldsymbol{\mu}}\in\bigcup_{j=1}^{t}\mathcal{Q}_j$. Each $\mathcal{Q}_j$ is obtained by removing all constant-$L$ points from a triangular branch before projection.}
\STATE Compute the associated parametric polynomial set $\mathcal{P}$ of vector field $(F_1,\ldots,F_n)$\;
\STATE Compute a decomposition $\{R_1,\ldots,R_s\}$ of $[\mathcal{P},\emptyset]$ into strong regular systems; set $\Psi:=\emptyset$ \label{line:invsurtril2}\;
\STATE Set $R:=[\{f_1,\ldots,f_{\ell_1}\},\emptyset]$\;
\FOR{$i = 1$ \TO $s$}
\STATE Compute a constructible-set representation of $\mathcal{N}_i:=\zero(R_i)\setminus\zero(R)$, refining it into regular systems when needed \;
    \IF{$\mathcal{N}_i=\emptyset$}
        \STATE Discard $\mathcal{N}_i$: this branch has no nontrivial solution\;
    \ELSE
        \STATE $\mathcal{Q}_i:=\operatorname{Proj}_{m}(\mathcal{N}_i)$\;
        \STATE  $\Psi:=\Psi\cup \{\mathcal{Q}_i\}$\;
    \ENDIF
\ENDFOR
\RETURN $\Psi$\;
\end{algorithmic}
\end{algorithm}

As a summary of this section, we remark that the two elimination methods are complementary. The Gr\"obner basis method, together with saturation and prime decomposition, gives a concise description of the Zariski closure of the existence locus and thus supplies necessary conditions. Triangular decomposition then refines each prime component branch by branch. When strong regular systems are used, their projection property yields pointwise necessary and sufficient conditions, while inequality constraints can be imposed directly during the decomposition. The following proposition formalizes the combination. It is based on Propositions \ref{proposition:nontrivial-condition} and \ref{proposition:dsys-trimethod}.

\begin{proposition}
\label{proposition:GB+dsys}
Let $\mathcal{P}$ denote the parametric polynomial set of system \eqref{L-dot_xx} with respect to the prescribed form \eqref{eq:preL}. Let $\mathfrak{I} = \langle \mathcal{P} \rangle : \langle f_1,\ldots,f_{\ell_1} \rangle^{\infty}$ and $\frI_{\ell_1+\ell_2+2}$ be the $(\ell_1+\ell_2+2)$-th elimination ideal of $\mathfrak{I}$. Suppose that the minimal prime decomposition of $\frI_{\ell_1+\ell_2+2}$ is
\[
\sqrt{\frI_{\ell_1+\ell_2+2}}=\bigcap_{i=1}^{r}\mathfrak{p}_i.
\]
Let $\mathcal{G}_i$ be a generating set of $\mathfrak{p}_i$ for $i\in\{1,\ldots,r\}$.
Then system \eqref{L-dot_xx} admits a nontrivial invariant algebraic surface of the form \eqref{eq:preL} if and only if there exists $i\in\{1,\ldots,r\}$ such that the polynomial system $[\mathcal{P}\cup\mathcal{G}_i,\emptyset]$ admits a solution with $(f_1,\ldots,f_{\ell_1})\neq\boldsymbol{0}$.
\end{proposition}

\begin{proof}
The zero set of the elimination ideal decomposes into the irreducible varieties of its minimal prime components:
\[
\zero(\frI_{\ell_1+\ell_2+2})=\zero(\mathfrak{p}_1)\cup\cdots\cup\zero(\mathfrak{p}_{r}).
\]
\emph{Sufficiency.} If $[\mathcal{P}\cup\mathcal{G}_i,\emptyset]$ admits a solution $(\overline{\boldsymbol{\mu}},\boldsymbol{f},\boldsymbol{a})$ with $(f_1,\ldots,f_{\ell_1})\neq\boldsymbol{0}$ for some $i$, then $(\overline{\boldsymbol{\mu}},\boldsymbol{f},\boldsymbol{a})\in\zero(\mathcal{P})$, and the claim follows from Lemma \ref{lemma:dsys-equivalence}.

\noindent\emph{Necessity.} If the system admits a nontrivial invariant algebraic surface for the parameter value $\overline{\boldsymbol{\mu}}$, then there exist $\boldsymbol{f}$ with $(f_1,\ldots,f_{\ell_1})\neq\boldsymbol{0}$ and $\boldsymbol{a}$ such that $(\overline{\boldsymbol{\mu}},\boldsymbol{f},\boldsymbol{a})\in\zero(\mathcal{P})$; by Proposition~\ref{proposition:nontrivial-condition}, $\overline{\boldsymbol{\mu}}\in\zero(\frI_{\ell_1+\ell_2+2})$, hence $\overline{\boldsymbol{\mu}}\in\zero(\mathfrak{p}_i)$ for some $i\in\{1,\ldots,r\}$, which means that $(\overline{\boldsymbol{\mu}},\boldsymbol{f},\boldsymbol{a})$ is a nontrivial solution of $[\mathcal{P}\cup\mathcal{G}_i,\emptyset]$.
\end{proof}

\begin{figure}[H]
\centering
\begin{tikzpicture}[
 font=\small,
 box/.style={draw, rounded corners=2pt, align=center, text width=7.7cm, inner sep=5pt},
 lab/.style={font=\footnotesize\itshape, align=left},
 sta/.style={font=\footnotesize\itshape, align=left, text width=4.0cm},
 arr/.style={workflow forward, double, double distance=1.1pt, thick},
 equ/.style={workflow equivalence, double, double distance=1.1pt, thick},
 node distance=9mm]
\node[box] (geo) {Prescribed invariant algebraic\\surface $L$ and cofactor $K$};
\node[box, below=of geo] (alg) {Parametric polynomial system $\mathcal{P}=0$ with $(f_1,\ldots,f_{\ell_1})\neq(0,\ldots,0)$};
\node[box, below=of alg] (clo) {Elimination ideal with saturation\\ and its minimal prime decomposition};
\node[box, below=of clo] (tri) {Strong regular system decomposition\\ of $\mathcal{P}\cup\mathcal{G}_i$ on each prime component};
\node[box, below=of tri, fill=yellow!15] (iff) {Finite parameter conditions\\with explicit $L$ and $K$\\(equations and inequations)};
\draw[equ] (geo) -- node[lab, right, xshift=2mm]{comparing coefficients\\ (Lemma~\ref{lemma:dsys-equivalence})} (alg);
\draw[arr] (alg) -- node[lab, right, xshift=2mm]{Gr\"obner basis elimination,\\ saturation and prime decomposition} (clo);
\draw[arr] (clo) -- node[lab, right, xshift=2mm]{branch-wise triangular decomposition\\ (projection property)} (tri);
\draw[arr] (tri) -- node[lab, right, xshift=2mm]{remove trivial points by\\ constructible-set difference before projection} (iff);
\node[sta, right=9mm of clo] {necessary conditions\\ (Zariski closure)};
\node[sta, right=9mm of tri] {branch refinement\\with nonvanishing\\constraints};
\node[sta, right=9mm of iff] {necessary and sufficient\\ (pointwise)};
\end{tikzpicture}
\caption{The combined Gr\"obner basis and triangular decomposition routine.}
\label{figure:workflow}
\end{figure}
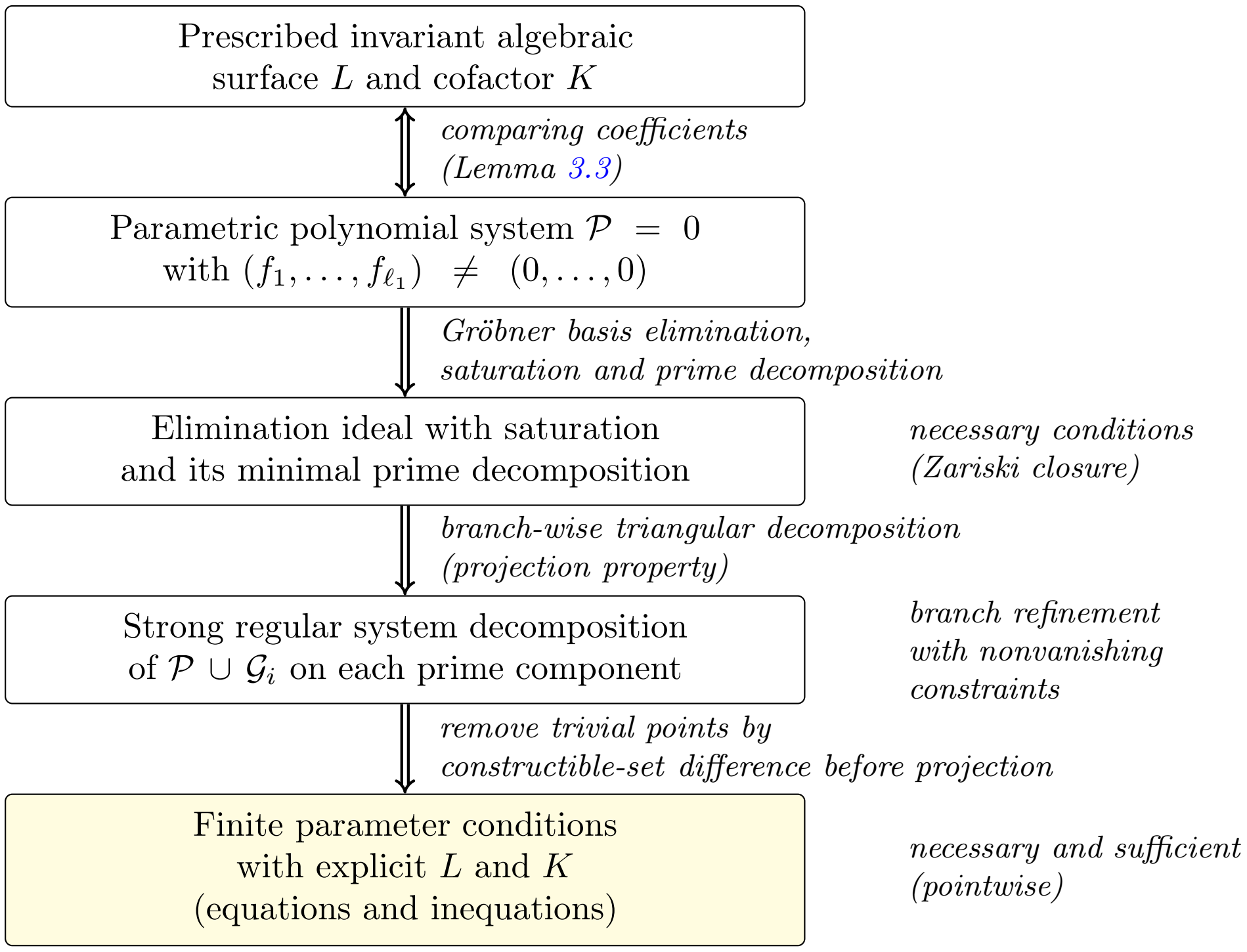

In practice, the prime components obtained in the first step reduce the global problem to several smaller systems $[\mathcal{P}\cup\mathcal{G}_i,\emptyset]$. A subsequent decomposition into strong regular systems preserves the projection property on each resulting branch. While prescribed nonvanishing conditions can be incorporated directly into the input as inequations, the nontriviality requirement on $L$ is handled by Algorithm \ref{algorithm:invsurtri} through the constructible-set difference $\zero(R_i)\setminus\zero(R)$ on every branch before projection. This combination retains the compact global description given by prime decomposition and the pointwise exactness given by strong regular systems and constructible-set operations. Figure \ref{figure:workflow} summarizes this combined routine.

\subsection{Computational results}
\label{sec:result}

\subsubsection{Belousov--Zhabotinsky reaction model}
\label{sec:result-1}

We employ the Belousov--Zhabotinsky (BZ) reaction model to illustrate the two proposed methods. First, we demonstrate that combining these methods yields the exact existence conditions for invariant algebraic planes, together with representative plane-cofactor pairs (see Theorem \ref{theorem:BZ-invariant-planes}). Next, we apply the triangular decomposition method independently under nonzero parameter constraints. Finally, we compare both methods in terms of computational efficiency and the characteristics of their respective results.

Consider the following three-dimensional dynamical system describing a modified BZ reaction model (see \cite{tot2018rea}):
\begin{equation}
\label{equation:dsys-BZ-system}
\left\{
\begin{aligned}
\dot{x}&=-k_1xy+k_2ay-2k_3x^2+k_4ax-\frac{1}{2}k_5xz,\\
\dot{y}&=-k_1xy-k_2ay+k_6mz,\\
\dot{z}&=2k_4ax-k_5xz-2k_6mz,
\end{aligned}
\right.
\end{equation}
where $\boldsymbol{\mu}=(k_1,k_2,k_3,k_4,k_5,k_6,a,m)$ are system parameters, and $\boldsymbol{x}=(x,y,z)$ represent the concentrations of chemical species, respectively. While the BZ model has been extensively studied in terms of reaction kinetics, bifurcation analysis, and numerical simulation (see \cite{fie1974osc,sco1994che}), a systematic investigation of its invariant algebraic surfaces using symbolic computation remains lacking.

First, taking the linear ansatz as an example, assume that system \eqref{equation:dsys-BZ-system} admits an invariant algebraic plane of the form
\begin{equation}
\label{equation:dsys-BZ-Ldeg1}
L = 1+f_1x+f_2y+f_3z=0,
\end{equation}
where $f_1,f_2,f_3$ are undetermined coefficients. Note that for a general linear ansatz $L = f_0+f_1x+f_2y+f_3z$, we may divide it into two cases for separate computation: (1) $f_0 \neq 0$, by scaling we may set $f_0 = 1$, i.e., $L = 1+f_1x +f_2y +f_3z$; (2) $f_0 = 0$, i.e., $L = f_1x+f_2y+f_3z=0$ represents a plane passing through the origin. In this section, we focus primarily on the case \eqref{equation:dsys-BZ-Ldeg1} since the case for $f_0=0$ may be handled similarly.

Since the maximum degree of the polynomials on the right-hand side of system \eqref{equation:dsys-BZ-system} is 2, we assume the cofactor has the form
$K = s_0+s_1x+s_2y+s_3z$. Expanding the equation $\mathcal{X}L - KL$ and collecting the coefficients of the monomials in $x,y,z$,
we obtain the following polynomials:
\begin{equation}
\label{equation:d-system-BZF}
\setlength{\jot}{2pt}
\left\{\;
\begin{aligned}
P_1 &= s_0, \\
P_2 &= -af_1k_4 - 2af_3k_4 + f_1s_0 + s_1, \\
P_3 &= -f_2k_6m + 2f_3k_6m + f_3s_0 + s_3, \\
P_4 &= f_3s_3, \\
P_5 &= f_2s_2, \\
P_6 &= 2f_1k_3 + f_1s_1, \\
P_7 &= f_1k_5 + 2f_1s_3 + 2f_3k_5 + 2f_3s_1, \\
P_8 &= f_1k_1 + f_1s_2 + f_2k_1 + f_2s_1, \\
P_9 &= f_2s_3 + f_3s_2, \\
P_{10} &= -af_1k_2 + af_2k_2 + f_2s_0 + s_2.
\end{aligned}
\right.
\end{equation}

The identity $\mathcal{X}L - KL\equiv0$ holds if and only if $P_1=\cdots=P_{10}=0$. Consequently, by Lemma \ref{lemma:dsys-equivalence}, we immediately obtain the following corollary.

\begin{corollary}
\label{corollary:d-system-BZF}
System~\eqref{equation:dsys-BZ-system} admits a nontrivial invariant algebraic plane of the form~\eqref{equation:dsys-BZ-Ldeg1} if and only if the polynomial system~\eqref{equation:d-system-BZF} has a solution with $(f_1,f_2,f_3) \neq (0,0,0)$.
\end{corollary}

% \begin{equation}
% \boxed{
% \zero(\mathcal F\cup\mathcal G_i)
% =
% \bigcup_{j=1}^{s_i}
% \zero\!\left(\mathcal T_{i,j}/\mathcal U_{i,j}\right),
% \qquad i=1,\ldots,8.
% }
% \end{equation}

% \begin{equation}
% \boxed{
% \zero(\mathcal F\cup\mathcal G_i)
% \setminus\zero(f_1,f_2,f_3)
% =
% \bigcup_{j=1}^{s_i}
% \left[
% \zero\!\left(\mathcal T_{i,j}/\mathcal U_{i,j}\right)
% \setminus\zero(f_1,f_2,f_3)
% \right],
% \qquad i=1,\ldots,8.
% }
% \end{equation}

\begin{theorem}\label{theorem:BZ-invariant-planes}
System \eqref{equation:dsys-BZ-system} admits a nontrivial invariant algebraic plane of the form \eqref{equation:dsys-BZ-Ldeg1} if and only if the system parameters satisfy at least one of the following conditions:
\begin{enumerate}[label=\textup{(\arabic*)}, ref=\arabic*]
\item $\mathcal{C}_{1,1}=\{k_6m=0,\quad ak_4k_5\neq0\}$ \, $\vee$ \,
$\mathcal{C}_{1,2}=\{k_6m=0,\quad ak_4=k_5=0\}$;

\item $\mathcal{C}_{2,1}=\{k_1=k_3=k_5,\quad ak_4k_5\neq0\}$;

\item $\mathcal{C}_{3,1}=\{k_2=0,\quad k_1=k_5,\quad ak_4k_5\neq0\}$;

\item $\mathcal{C}_{4,1}=\{k_5=0,\quad k_1k_4=2k_2k_3,\quad ak_3k_4\neq0\}$;

\item $\mathcal{C}_{5,1}=\{k_3=k_4=k_5=0,\quad ak_1k_2\neq0\}$;

\item $\mathcal{C}_{6,1}=\{k_5=0,\quad k_6m=0,\quad k_2=0,\quad ak_3k_4\neq0,\quad k_1(k_1-2k_3)\neq0\}$ \, $\vee$ \, \\
$\mathcal{C}_{6,2}=\{k_5=0,\quad k_6m=0,\quad k_1=k_3,\quad ak_3k_4\neq0,\quad k_1(k_1-2k_3)\neq0\}$;

\item $\mathcal{C}_{7,1}=\{k_1=k_4=k_5=k_2k_3=0,\quad a\neq0\}$;

\item $\mathcal{C}_{8,1}=\{a=k_1=k_5=0\}$;

\item $\mathcal{C}_{9,1}=\{k_3=ak_2=0\}$;

\item $\mathcal{C}_{10,1}=\{k_1=k_6m=ak_2k_3=0\}$.

\end{enumerate}
In each case, representative invariant algebraic planes and their corresponding cofactors are presented explicitly in Table \ref{table:BZ-deg1}.
% Specifically, the first column lists explicit planes realizing each parameter condition. Together with their cofactors, these examples establish sufficiency by direct substitution, whereas necessity follows from the completeness of the parameter space projection.
% We remark that, on each branch, all solutions can be recovered by successive back-substitution on the computed regular systems; since this paper focuses on Problem ($\mathcal{P}_1$), that is, the complete classification of the parameter conditions, Table \ref{table:BZ-deg1} lists only representative planes and cofactors.

\end{theorem}

\begin{table}[H]
\centering
\caption{Representative invariant algebraic planes and cofactors of the BZ system \eqref{equation:dsys-BZ-system} under the linear ansatz \eqref{equation:dsys-BZ-Ldeg1}, where $c\in\mathbb{C}\setminus\{0\}$ is an arbitrary constant.}
\vspace{5pt}
\footnotesize
\tabcolsep 4pt
\begin{tabular}{>{\raggedright\arraybackslash}p{4.7cm}>{\raggedright\arraybackslash}p{3.6cm}>{\raggedright\arraybackslash}p{4.4cm}}
\toprule
Invariant algebraic planes $L$ & Cofactors $K$ & Parameter conditions \\
\midrule
$1-\dfrac{k_5}{2ak_4}z$ \newline $1+c\,z$  & $-k_5x$ \newline $0$ & $\mathcal{C}_{1,1}$ \newline $\mathcal{C}_{1,2}$ \\[10pt]
$1-\dfrac{k_5}{ak_4}x-\dfrac{k_5}{ak_4}y-\dfrac{k_5}{2ak_4}z$ & $-2k_5x$ & $\mathcal{C}_{2,1}$ \\[8pt]
{$\begin{aligned}&1+\tfrac{k_5-2k_3}{ak_4}x-\tfrac{k_5}{ak_4}y-\tfrac{k_5}{2ak_4}z\\[2pt] &1-\tfrac{k_5}{ak_4}y-\tfrac{k_5}{2ak_4}z\end{aligned}$} & {$\begin{aligned}&-2k_3x\\[2pt]&-k_5x\end{aligned}$} & $\mathcal{C}_{3,1}$ \\[16pt]
$1-\dfrac{2k_3}{ak_4}x$ & $-2k_3x-k_1y$ & $\mathcal{C}_{4,1}$ \\[8pt]
$1-\dfrac{k_1}{ak_2}x$ & $-k_1y$ & $\mathcal{C}_{5,1}$ \\[8pt]
$1-\dfrac{2k_3}{ak_4}x+\dfrac{2k_3k_1}{ak_4(k_1-2k_3)}y$ & $-2k_3x$ & {$\begin{array}{@{}l@{}}\mathcal{C}_{6,1}\\[2pt]\mathcal{C}_{6,2}\end{array}$}\\[12pt]
$1+c\,(2y+z)$ \ if $k_2=0$;\newline $1+c\,(2x+2y+z)$ \ if $k_3=0$ & $0$ & $\mathcal{C}_{7,1}$ \\[8pt]
$1+c\,(2y+z)$ & $0$ & $\mathcal{C}_{8,1}$ \\[4pt]
$1+c\,(-2x+2y+z)$ & $0$ & $\mathcal{C}_{9,1}$\\[4pt]
$1+c\,y$ \ if $ak_2=0$;\newline $1+c\,(-2x-2y+z)$ \ if $k_3=0$ & $0$ & $\mathcal{C}_{10,1}$\\
\bottomrule
\end{tabular}
\label{table:BZ-deg1}
\end{table}

\begin{proof}
The algebraic computations are carried out over $\mathbb{C}$; conclusions for the physically relevant real parameters are then obtained by restricting to the real part of the existence locus.

\emph{Sufficiency.} Substituting each plane and its corresponding cofactor from Table \ref{table:BZ-deg1} directly into $\mathcal{X}L=KL$ verifies that the identity holds under the corresponding condition.

\emph{Necessity.} By Corollary \ref{corollary:d-system-BZF}, it suffices to determine the parameter values under which system \eqref{equation:d-system-BZF} admits a nontrivial solution $(f_1,f_2,f_3)\neq(0,0,0)$. This computation proceeds in two steps.

\smallskip
\noindent\emph{Step 1: Elimination and prime decomposition}. Compute the saturated ideal $\mathfrak{I}=\langle\mathcal{P}\rangle:\langle f_1,f_2,f_3\rangle^{\infty}$ in $\mathbb{C}[\boldsymbol{\mu},\boldsymbol{s},\boldsymbol{f}]$, and its elimination ideal $\mathfrak{I}_7=\mathfrak{I}\cap\mathbb{C}[\boldsymbol{\mu}]$, where $\boldsymbol{\mu}=(k_1,\ldots,k_6,a,m)$, $\boldsymbol{s}=(s_0,\ldots,s_3)$ and $\boldsymbol{f}=(f_1,f_2,f_3)$. The subscript is $7$ because the seven undetermined quantities $s_0,\ldots,s_3,\allowbreak f_1,f_2,f_3$ are eliminated; this is the normalized count $\ell_1+\ell_2+1$ for $\ell_1=\ell_2=3$, since the constant coefficient of $L$ has been fixed to $1$. The minimal prime decomposition is
\[
\sqrt{\mathfrak{I}_7}=\mathfrak{p}_1\cap\cdots\cap\mathfrak{p}_8,
\]
where $\mathfrak{p}_1=\langle m\rangle$, $\mathfrak{p}_2=\langle k_6\rangle$, $\mathfrak{p}_3=\langle a,k_5,k_1\rangle$, $\mathfrak{p}_4=\langle k_5,k_1k_4-2k_2k_3\rangle$, $\mathfrak{p}_5=\langle k_1-k_5,k_3-k_5\rangle$, $\mathfrak{p}_6=\langle k_2,k_1-k_5\rangle$, $\mathfrak{p}_7=\langle a,k_3\rangle$ and $\mathfrak{p}_8=\langle k_2,k_3\rangle$. Consequently, the existence locus is contained in $\bigcup_{i=1}^{8}\zero(\mathfrak{p}_i)$, and its Zariski closure equals this union.

\smallskip
\noindent\emph{Step 2: Branch-wise triangular solving}. For each $i\in\{1,\ldots,8\}$, let $\mathcal{G}_i$ be a generating set of $\mathfrak{p}_i$ and decompose the polynomial system $[\mathcal{P}\cup\mathcal{G}_i,\emptyset]$ into finitely many strong regular systems $R_{ij}$ via \textsf{RegSer}, using an elimination ordering in which the parameters $\boldsymbol{\mu}$ are smaller than the unknowns $s_0,\ldots,s_3,f_1,f_2,f_3$.
Set $R:=[\{f_1,f_2,f_3\},\emptyset]$. For each resulting system $R_{ij}$, set
\[
\mathcal{N}_{ij}=\zero(R_{ij})\setminus\zero(R),
\]
and compute $\operatorname{Proj}_{8}(\mathcal{N}_{ij})$, rather than projecting $\zero(R_{ij})$ before its trivial part is removed. Successive back-substitution on the sets $\mathcal{N}_{ij}$, with the disjunction $f_1\neq0\vee f_2\neq0\vee f_3\neq0$, also yields the corresponding nonconstant planes and cofactors. The union of the projected nontrivial parts over all $i$ and $j$, after absorbing redundant sub-branches, is exactly the disjunction of conditions (1)--(10).
\end{proof}

The branch-wise computation above provides more than the parameter conditions themselves: each branch retains all the information needed
to recover the coefficient solutions (indeed, for the present linear ansatz, the full solution families can be computed explicitly by
back-substitution). Nevertheless, Table \ref{table:BZ-deg1} lists only witnesses (i.e., representative plane--cofactor pairs) for each condition, as this paper focuses on problem ($\mathcal{P}_1$), that is, the complete classification of the parameter conditions.

With the complete classification at hand, it is natural to single out some nondegenerate cases where the classification simplifies considerably. For example, imposing the generic condition
$ak_4k_5\neq0$ yields a substantial simplification.

\begin{corollary}
\label{corollary:BZ-invariant-planes-generic}
Suppose $ak_4k_5\neq0$. Then system \eqref{equation:dsys-BZ-system} admits a nontrivial invariant algebraic plane of the form \eqref{equation:dsys-BZ-Ldeg1} if and only if the system parameters satisfy at least one of the following five conditions:
\begin{enumerate}[label=\textup{(\arabic*)}, ref=\arabic*]
\item $m=0$;
\item $k_6=0$;
\item $k_1=k_3=k_5$;
\item $k_2=0$, $k_1=k_5$;
\item $k_3=k_2=0$.
\end{enumerate}
\end{corollary}

\begin{proof}
Under the assumption that $ak_4k_5\neq0$, conditions (4)--(8) of Theorem \ref{theorem:BZ-invariant-planes} are ruled out because each requires $ak_4k_5=0$. Meanwhile, condition (10) implies $k_6m=0$ and is thus absorbed by condition (1). Since $a\neq0$, condition (9) reduces to $k_3=k_2=0$, while condition (1) itself simplifies to $m=0$ or $k_6=0$.
\end{proof}

\begin{remark}\label{remark:BZ-closure}
We present some observations about our computational methods:
\smallskip

(i) The full existence locus of Theorem \ref{theorem:BZ-invariant-planes} is a constructible set but not an algebraic variety in the parameter space: the inequations in conditions (1)--(7) are indispensable. For instance, the parameter point $(k_1,\ldots,k_6,a,m)=(1,1,1,1,1,1,0,0)$ satisfies $m=0$ and hence lies in the Zariski closure of the existence locus, but it violates each of conditions (1)--(10) (condition (1) fails since $ak_4k_5=0$ while $k_5\neq0$); accordingly, the specialized system~\eqref{equation:d-system-BZF} at this point has reduced Gr\"obner basis $\{s_0,s_1,s_2,s_3,f_1,f_2,f_3\}$ and admits only the trivial solution.
\smallskip

(ii) The Zariski closure of the existence locus is described exactly by the elimination ideal: for the prime decomposition $\sqrt{\mathfrak{I}_7}=\bigcap_{i=1}^{8}\mathfrak{p}_i$ computed in Step~1 of the proof, the closure of the locus cut out by conditions (1)--(10) is precisely $\bigcup_{i=1}^{8}\zero(\mathfrak{p}_i)$. More precisely, the sets cut out by conditions (1), (2), (3) and (4) have Zariski closures $\zero(\mathfrak{p}_1)\cup\zero(\mathfrak{p}_2)$, $\zero(\mathfrak{p}_5)$, $\zero(\mathfrak{p}_6)$ and $\zero(\mathfrak{p}_4)$, respectively, condition (8) coincides with $\zero(\mathfrak{p}_3)$, condition (9) coincides with $\zero(\mathfrak{p}_7)\cup\zero(\mathfrak{p}_8)$, and the remaining conditions are contained in these components. Moreover, these branches, together with the corresponding planes in Table~\ref{table:BZ-deg1}, can be found systematically by computing reduced Gr\"obner bases of $\langle\mathcal{P}\rangle$ over the function fields of the prime components: all denominators and leading coefficients then involve only the factors appearing in the inequations of conditions (1)--(6), and the bases specialize well precisely off the corresponding hypersurfaces. The branch-wise triangular solving in Step~2 of the proof refines this Zariski closure into the exact pointwise classification, detecting automatically the degenerate strata where specialization fails.
\smallskip

(iii) Under conditions (7)--(10), as well as in the degenerate case $ak_4=k_5=0$ of condition (1), the cofactor vanishes identically; consequently, the corresponding polynomials $L$ are, in fact, first integrals of the system.

\end{remark}

%These non-zero or positivity constraints from the model background make the feasible parameter space of the system constitute a semi-algebraic set.

In the BZ reaction model, $k_1, \ldots, k_6$ denote the rate constants of elementary reactions, which are generally positive real numbers, while $a$ and $m$ typically represent the initial concentrations of bromide and cerium ions (or related parameters), which are also inherently positive \cite{fie1974osc, tys1976bel}. It is therefore natural to investigate invariant algebraic surfaces under nonzero (in fact, strictly positive) parameter constraints. In theory, cylindrical algebraic decomposition (CAD) provides a universal framework for addressing problems with inequality constraints (see \cite{col1976qua, che2016qua}). However, when applied to computing invariant algebraic surfaces of dynamical systems, the CAD algorithm often generates a prohibitively large number of semi-algebraic cells. Consequently, CAD is better suited for verifying specific parameter points than for a complete classification. In contrast, the triangular decomposition method under inequality constraints operates directly on a polynomial system $[\mathcal{F},\mathcal{H}]$, which comprises equations $\mathcal{F}=0$ and inequations $\mathcal{H}\neq0$. It returns finitely many constrained triangular systems $[\mathcal{T}_i,\mathcal{H}_i]$ whose zero sets satisfy $\zero(\mathcal{F} / \mathcal{H}) = \bigcup_{i} \zero(\mathcal{T}_i / \mathcal{H}_i)$ (see Section \ref{sec:preli}; see also \cite{wan1998dec, wan2001eli}). In this manner, nonvanishing requirements on the parameters are imposed as part of the input, without any additional post-processing for these parameter constraints.

%In contrast, triangular decomposition with inequality constraints applies directly to a polynomial system $[\mathcal{F},\mathcal{H}]$, with equality constraints $\mathcal{F}=0$ and inequations $\mathcal{H}\neq0$, and returns finitely many triangular systems with inequality constraints $[\mathcal{T}_i,\mathcal{H}_i]$ whose zero sets satisfy $\zero(\mathcal{F} / \mathcal{H}) = \bigcup_{i} \zero(\mathcal{T}_i / \mathcal{H}_i)$ (see Section~\ref{sec:preli}, and also \cite{wan1998dec,wan2001eli}). In this way, non-vanishing requirements on the system parameters are incorporated as part of the input, without any additional processing.

Hereafter, we employ the RegSer algorithm to perform strong regular system decomposition under inequality constraints. Specifically, we consider the linear invariant algebraic planes of the BZ system, where $P_1,\ldots,P_{10}$ are defined in \eqref{equation:d-system-BZF}. By executing the computation under the constraint that all system parameters are nonzero:
\begin{equation*}
\operatorname{RegSer}([P_1, \dots, P_{10}], \{k_1,k_2,k_3,k_4,k_5,k_6,a,m\}),
\end{equation*}
we obtain the desired decomposition:
\begin{equation*}
\begin{split}
[\mathcal{T}_1,\mathcal{U}_1] &= [[s_0,s_1,s_2,s_3,f_1,f_2,f_3],\{a,k_1,k_2,k_3,k_4,k_5,k_6,m\}], \\
[\mathcal{T}_2,\mathcal{U}_2] &=[[k_1 - k_3, k_1 - k_5, s_0, 2k_1 + s_1, s_2, s_3, af_1k_4 + k_1, f_1 - f_2, f_1 - 2f_3], \{a, k_1, k_2, k_4, k_6, m\}],
\end{split}
\end{equation*}
In the first branch, $f_1=f_2=f_3=0$ corresponds to the trivial solution.
Solving the second branch: from $k_1-k_3=0$, $k_1-k_5=0$ we get $k_3=k_5=k_1$; from $af_1k_4+k_1=0$ we solve $f_1=-\frac{k_1}{ak_4}$;
from $f_1-f_2=0$, $f_1-2f_3=0$ we get $f_2=f_1$, $f_3=\frac{1}{2}f_1$.
Substituting into equation~\eqref{equation:dsys-BZ-Ldeg1}~and simplifying yields
$$L = 1-\frac{k_1}{ak_4}x-\frac{k_1}{ak_4}y-\frac{k_1}{2ak_4}z.$$

The above computational results lead directly to the following proposition.

\begin{proposition}
\label{proposition:d-system-T2}
When all system parameters are nonzero, the BZ
system \eqref{equation:dsys-BZ-system} admits a nontrivial invariant algebraic plane of the form \eqref{equation:dsys-BZ-Ldeg1} if and only if $k_3=k_5=k_1$.
\end{proposition}

The efficiency of triangular decomposition under inequality constraints allows us to treat the full quadratic ansatz \eqref{equation:bzdeg2} of the BZ system, for which the Gr\"obner basis method did not terminate within 3600 seconds in \textsc{Maple}:
\begin{equation}
\label{equation:bzdeg2}
L = 1+f_1x+f_2y+f_3z+f_4x^2+f_5y^2+f_6z^2+f_7xy+f_8xz+f_9yz=0.
\end{equation}
Computational results show that, when all system parameters are nonzero, the associated parametric polynomial set of the quadratic ansatz \eqref{equation:bzdeg2} consists of $20$ polynomials, and its triangular decomposition yields four branches: the trivial branch on which all undetermined coefficients vanish, two branches recovering the invariant plane of Proposition~\ref{proposition:d-system-T2} and its square, respectively, and a branch whose parameter conditions involve $k_1^2-4k_1k_3+8k_3^2=0$, which admits no real solutions since $k_1^2-4k_1k_3+8k_3^2=(k_1-2k_3)^2+4k_3^2>0$ whenever $k_3\neq0$. Hence, under the normalization $L(0)=1$ used in \eqref{equation:bzdeg2}, the BZ reaction model admits no irreducible quadratic invariant algebraic surfaces in real space when all system parameters are nonzero.

We evaluated the Gr\"obner basis and triangular decomposition methods on the same BZ ansatzes and on the same machine. The Gr\"obner basis timings $t_g^{\mathrm{MS}}$ use a hybrid implementation: \textsc{Maple} performs saturation and elimination, while \textsc{Singular} computes the minimal associated primes with \textsf{minAssGTZ}. The triangular decomposition timings $t_t$ use \textsf{TriSer} in \textsc{Epsilon}, followed by constructible set difference, projection, and union. The barred entries $\overline{t_t}$ and $\overline{n_t}$ refer to the computations with all system parameters constrained to be nonzero.

% \begin{table}
% \centering
% \caption{Experimental data for computing different invariant algebraic surface ansatzes using Gr\"obner basis and triangular set methods (unit: seconds).}
% \vspace{5pt}
% \small
% \tabcolsep 6pt
% % \setlength{\extrarowheight}{2pt}
% \begin{tabular}{>{\raggedright\arraybackslash}p{6.5cm}cccccccc}
% \toprule
% \textbf{Invariant Algebraic Surfaces} & $n_f$ & $t_g$ &$\overline{t_g}$ & $t_t$ & $\overline{t_t}$ & $n_g$ & $n_t$ & $\overline{n_t}$ \\
% \midrule
%  ~$\overline{L} = 1+f_1x+f_2y+f_3z+f_5y^2=0$ & 13 & 1.019& 0.983 & 0.673 & 0.090 & 8 & 57 & 2 \\
%  ~$\tilde{L} = 1+f_1x+f_2y+f_3z+f_7xy=0$ & 13 & 0.793& 0.801 & 0.729  & 0.079  & 8  & 61  & 2 \\
%  ~$\hat{L} = 1+f_1x+f_2y+f_3z+f_4x^2+f_7xy=0$ & 15 & 17.933 & 16.746 & 1.633 & 0.035 & 8 & 68 & 2\\[6pt]
%  ~{$\begin{aligned}
% &L = 1 + f_1x + f_2y + f_3z + f_4x^2 + f_5y^2 + f_6z^2 \\
% &\quad \quad + f_7xy + f_8xz + f_9yz=0
% \end{aligned}$} & 20 & $>3600$ & $>3600$ & 254.615  & 1.137 & -- & 222 & 4\\
% \bottomrule
% \end{tabular}
% \label{table:dsys-exp}
% \end{table}

\begin{table}[H]
\centering
\caption{Experimental data for computing different invariant algebraic surface ansatzes using a \textsc{Maple}--\textsc{Singular} Gr\"obner basis route and triangular decomposition (in seconds).}
\vspace{5pt}
\scriptsize
\tabcolsep 3pt
\resizebox{\textwidth}{!}{%
\begin{tabular}{>{\raggedright\arraybackslash}p{6.2cm}ccccccc}
\toprule
\textbf{Invariant Algebraic Surfaces} & $n_f$ & $t_g^{\mathrm{MS}}$ & $t_t$ & $\overline{t_t}$ & $n_g$ & $n_t$ & $\overline{n_t}$ \\
\midrule
~$\overline{L}=1+f_1x+f_2y+f_3z+f_5y^2=0$
& 13 & 0.425 & 7.328 & 0.219 & 8 & 57 & 2 \\[3pt]
~$\tilde{L}=1+f_1x+f_2y+f_3z+f_7xy=0$
& 13 & 0.452 & 7.500 & 0.219 & 8 & 61 & 2 \\[3pt]
~$\hat{L}=1+f_1x+f_2y+f_3z+f_4x^2+f_7xy=0$
& 15 & 5.155 & 10.844 & 0.235 & 8 & 68 & 2 \\[6pt]
~{$\begin{aligned}
&L=1+f_1x+f_2y+f_3z+f_4x^2+f_5y^2+f_6z^2\\
&\quad+f_7xy+f_8xz+f_9yz=0
\end{aligned}$}
& 20 & $>3600$ & 317.766 & 2.282 & -- & 222 & 4 \\
\bottomrule
\end{tabular}}
\label{table:dsys-exp}
\end{table}

In Table \ref{table:dsys-exp}, the \textbf{Invariant Algebraic Surfaces} column records the assumed invariant algebraic surface forms, while $n_f$ records the number of polynomials in the associated parametric polynomial set. The $t_g^{\mathrm{MS}}$ and $t_t$ columns record the computation times for the Gr\"obner basis and triangular decomposition methods, respectively, using a hybrid \textsc{Maple}--\textsc{Singular} implementation for the former and \textsf{TriSer} in \textsc{Epsilon} for the latter. The barred entries $\overline{t_t}$ and $\overline{n_t}$ refer to the corresponding computation under the nonvanishing constraints on all system parameters $(k_1,k_2,k_3,k_4,k_5,k_6,a,m)$. The columns $n_g$, $n_t$, and $\overline{n_t}$ report the numbers of branches obtained: $n_g$ counts the minimal prime components of the elimination ideal, $n_t$ counts the triangular systems returned by \textsf{TriSer}, and $\overline{n_t}$ gives the corresponding count under the nonvanishing constraints. These values measure the sizes of the decompositions. A dash -- indicates that no prime decomposition was obtained within the time limit. The Gr\"obner basis method usually produces fewer branches because it uses the minimal prime components of the elimination ideal, which provide concise necessary conditions; however, further branch-wise computation is needed to obtain pointwise necessary and sufficient conditions. The triangular decomposition method is particularly effective for higher degree prescribed forms and inequality constraints, and its computational efficiency further improves when nonzero constraints are imposed on the system parameters.
Experiments were performed using Maple 2026, Singular 4.0.2, and Epsilon 0.618 on an Intel Core Ultra 5 336H with 32 GB RAM. Table~\ref{table:dsys-exp} reports CPU times. The hybrid Gr\"obner basis timings include Maple's saturation and elimination steps and Singular's \textsf{facstd} and \textsf{minAssGTZ} computations. The triangular decomposition timings include \textsf{TriSer} and the subsequent constructible-set difference, projection, and union. A time limit of 3600 seconds was imposed; the entry $>3600$ indicates that the Gr\"obner-basis route did not return a prime decomposition within this limit.

\subsubsection{A class of Li\'enard differential systems}
\label{sec:result-2}

We consider the following Li\'enard differential system:
\begin{equation}\label{LD-mod1}
\begin{split}
\dot{x}=y,\quad \dot{y}=-(\zeta x^2+\beta x+\alpha)y-(\varepsilon x^3+ex^2+\sigma x+\delta);\quad \zeta\neq0,\quad \varepsilon\neq0.
\end{split}
\end{equation}
The changes of variables $x\mapsto A(x+B)$, $y\mapsto Ay$, $A\neq0$ relates system \eqref{LD-mod1} with its simplified version at $\zeta=3,e=0$. Thus without loss of generality we set $\zeta=3,e=0$. All the parameters are supposed to be from the field $\mathbb{C}$ with $\varepsilon\neq0$. The classification of irreducible invariant algebraic curves of system \eqref{LD-mod1} was studied by Demina \cite{Demina2018}. Here using our computational approach, we provide some new findings of irreducible invariant algebraic curves of system \eqref{LD-mod1}.

Introducing the invertible change of variables
\begin{equation}\label{LD-change}
x=s, \quad y=z-h(s)\quad \leftrightarrow \quad s=x, \quad z=y+h(x),
\end{equation}
where $h(x)=x^3+\beta x^2/2+(\alpha-\varepsilon)x$, system \eqref{LD-mod1} becomes the following dynamical system
\begin{equation}\label{LD-mod2}
\begin{split}
\dot{s}=z-s^3-\frac{\beta}{2}s^2+(\varepsilon-\alpha)s,\quad \dot{z}=\frac{\beta\varepsilon}{2}s^2+(\varepsilon(\alpha-\varepsilon)-\sigma)s-\delta-\varepsilon z.
\end{split}
\end{equation}
There exists the one-to-one correspondence between irreducible invariant algebraic curves $F(x,y)$ of Li\'enard system \eqref{LD-mod1} and irreducible invariant algebraic curves $G(s,z)$ of system \eqref{LD-mod2}.

The following proposition is a minor variant of \cite[Theorem 1.2]{Demina2018}.

\begin{proposition}\label{LD-Prop1}
If $G(s,z)$ is an irreducible invariant algebraic curve of system \eqref{LD-mod2}, then
\begin{equation}
G(s,z)=\mu_0z^N+\sum_{i=0}^{N-1}c_i(s)z^i,\quad N\in\mathbb{N},
\end{equation}
and its cofactor $K(s,z)=A_2s^2+A_1s+A_0$, where $A_2=-M$ with $M\in\{0,3\}$ being the degree of $G(s,z)$ with respect to $s$.
\end{proposition}

\begin{proof}
The fact $M\in\{0,3\}$ is a direct consequence of \cite[Theorem 1.2]{Demina2018}. We prove the remaining part. Note that the irreducible invariant algebraic curves of system \eqref{LD-mod2} satisfy the equation
\begin{equation}\label{LD2-iac1}
\left(z-s^3-\frac{\beta}{2}s^2+(\varepsilon-\alpha)s\right)G_s+\left(\frac{\beta\varepsilon}{2}s^2+(\varepsilon(\alpha-\varepsilon)-\sigma)s-\delta-\varepsilon z\right)G_z=KG.
\end{equation}
By direct substitution we verify that there are no invariant algebraic curves that do not depend on $z$. Let $G$ and $K$ have degrees $N$ and $\ell$ with respect to $z$ accordingly. Equation \eqref{LD2-iac1} implies ${\rm deg}_zK\leq1$ when $K\neq0$. Substituting relations $G=\mu(s)z^N+\cdots$, $\mu(s)\not\equiv0$ and $K=k_1(s)z+k_0(s)$ with $\mu(s),k_i(s)\in\mathbb{C}[s]$ into equation \eqref{LD2-iac1} and balancing higher-order terms in $z$ yields $\mu_{s}=k_1(s)\mu(s)$. Since $\mu(s)$ is a polynomial, we get $k_1(s)=0$, $\mu(s)=\mu_0$, where $\mu_0$ is a constant. In addition, we see that the cofactor $K$ does not depend on $z$.

Now suppose that $G$ and $K$ have degrees $M$ and $m$ with respect to $s$ accordingly. Taking expressions $G=\nu(z)s^M+\cdots$, $\nu(z)\not\equiv0$ and $K=A_ms^m$ with $\nu(z)\in\mathbb{C}[z]$, $A_m\in\mathbb{C}$. The bound ${\rm deg}_sK\leq2$ follows from \eqref{LD2-iac1}. Arguing as above and comparing $s^{M+2}$ gives
\[\frac{\beta\varepsilon}{2}\nu'(z)=(A_2+M)\nu(z).\]
Thus $A_2=-M$. This holds also when $\beta=0$. This completes the proof.

\end{proof}

Without loss of generality, we shall set $\mu_0=1$. In \cite{Demina2018}, all irreducible invariant algebraic curves for the Li\'enard system \eqref{LD-mod1} are obtained, showing that the degree of $G(s,z)$ with respect to $z$ is at most 3 (see the classification and argument in \cite[Theorem 1.3]{Demina2018}). Next, we classify the irreducible invariant algebraic curves of system \eqref{LD-mod1} by integrating the Gr\"obner basis elimination method and the triangular set method. This contribution presents an elimination-based derivation that corrects several erroneous entries in existing classifications. Consequently, it turns out that our algorithmic framework will work more conveniently and effectively when the underlying structure of the invariant algebraic curves is known.

\begin{theorem}\label{LD-main}
The Li\'enard system \eqref{LD-mod1} with $\zeta=3$, $e=0$, and $\varepsilon\neq0$ admits invariant algebraic curves if and only if restrictions on the parameters given in Table \ref{table:LD1} are satisfied. The expressions of irreducible algebraic invariants and their cofactors are of the form in Table \ref{table:LD1}.
\end{theorem}

\begin{table}
\centering
\caption{Irreducible invariant algebraic curves of system \eqref{LD-mod1} with $\zeta=3$, $e=0$, and $\varepsilon\neq0$.}
\vspace{5pt}
\footnotesize
\tabcolsep 4pt
\begin{tabular}{>{\raggedright\arraybackslash}p{7.2cm}c>{\raggedright\arraybackslash}p{5.4cm}}
\toprule
Invariant algebraic curves $F(x,y)$ & Cofactors $K(x,y)$ & Parameter conditions \\
\midrule
 ~$y+x^{3}+\left(\alpha -\varepsilon \right) x +\frac{\delta}{\varepsilon}$ & $-\varepsilon$ & $C_{1,1}=\{\beta=0,\,\, \sigma = \varepsilon  \left(\alpha -\varepsilon \right)\}$\\[6pt]
~$y+\frac{\varepsilon}{3} x-\frac{\beta\varepsilon}{9}$ & $-3x^{2}-\beta x-\alpha +\frac{\varepsilon}{3}$ & {$\begin{aligned}&C_{1,2}=\Big\{\delta =
-\frac{\beta  \varepsilon  \left(3 \alpha -\varepsilon \right)}{27},\\
&\quad \sigma =
\frac{\varepsilon  \left(-\beta^{2}+3 \alpha -\varepsilon \right)}{9}\Big\}\end{aligned}$}\\[6pt]
 ~$y^{2}+\big(x^{3}+\frac{4}{3} \varepsilon  x \big) y +\frac{ \varepsilon}{3}x^{4}+\frac{4 \varepsilon^{2}}{9} x^{2}$ & $-3 x^{2}-\frac{8 \varepsilon}{3}$ & {$\begin{aligned}&C_{2,1}=\Big\{\alpha=2\varepsilon,\,\,\beta=0,\\
&\quad \delta=0,\,\,\sigma = \frac{2}{9}\alpha^2\Big\}\end{aligned}$} \\[6pt] %[6pt]
 ~{$\begin{aligned}
&y^{2}+\big(x^{3}-\frac{4}{3} \varepsilon  x \big) y +\frac{\varepsilon}{3}x^{4} -\frac{8 \varepsilon^{2}}{9}x^{2}+\frac{16 \varepsilon^{3}}{27}
\end{aligned}$} & $-3x^2$ & {$\begin{aligned}&C_{2,2}=\Big\{\alpha = -\frac{2 \varepsilon}{3},\,\,\beta = 0,\\
&\quad \delta = 0,\,\,\sigma = -\frac{4 \varepsilon^{2}}{3}\Big\}\end{aligned}$} \\[6pt]
~{$\begin{aligned}
&y^{2}+\Big(x^{3}+\frac{28 \varepsilon}{27}x-\frac{160\varepsilon^{\frac{3}{2}}}{729}i\Big) y+\frac{\varepsilon}{3}x^{4} +\frac{8 \varepsilon^{2} }{27}x^{2}\\
&\quad-\frac{256 \varepsilon^{\frac{5}{2}}i}{2187}x-\frac{80 \varepsilon^{3}}{6561}
\end{aligned}$} & $-3 x^{2}-\frac{64 \varepsilon}{27}$ & {$\begin{aligned}&C_{2,3}=\Big\{\alpha = \frac{46 \varepsilon}{27},\,\, \beta=0,\\
&\quad \delta =
-\frac{32\varepsilon^{\frac{5}{2}}}{243}i,\,\,\sigma =
\frac{52 \varepsilon^{2}}{81}
\Big\}\end{aligned}$} \\[6pt]
~{$\begin{aligned}
&y^{2}+\Big(x^{3}+\frac{28 \varepsilon }{27}x+\frac{160\varepsilon^{\frac{3}{2}}}{729}i\Big) y+\frac{\varepsilon}{3}x^{4}+\frac{8 \varepsilon^{2} }{27}x^{2}\\
&\quad+\frac{256 \varepsilon^{\frac{5}{2}}i}{2187}x-\frac{80 \varepsilon^{3}}{6561}
\end{aligned}$} & $-3 x^{2}-\frac{64 \varepsilon}{27}$ & {$\begin{aligned}&C_{2,4}=\Big\{\alpha = \frac{46 \varepsilon}{27},\,\, \beta=0,\\
&\quad \delta =
\frac{32\varepsilon^{\frac{5}{2}}}{243}i, \,\,\sigma =
\frac{52 \varepsilon^{2}}{81}
\Big\}\end{aligned}$} \\[6pt]
~{$\begin{aligned}
&y^{3}+\Big(2 x^{3}+\frac{10}{3} \varepsilon  x \Big) y^{2}\\
&\quad +\Big(x^{6}+\frac{11}{3} \varepsilon  x^{4}+\frac{125}{36} \varepsilon^{2} x^{2}\Big) y \\
&\quad+\frac{\varepsilon}{3}x^{7} +\frac{11 \varepsilon^{2}}{9}x^{5}+\frac{125 \varepsilon^{3}}{108}x^{3}
\end{aligned}$} & $-3 x^{2}-\frac{25 \varepsilon}{6}$ & {$\begin{aligned}&C_{3,1}=\Big\{\beta = 0,\,\, \delta = 0,\\
&\quad \sigma = \frac{25 \varepsilon^{2}}{18},\,\,\alpha = \frac{5 \varepsilon}{2}\Big\}\end{aligned}$} \\[6pt]
~{$\begin{aligned}
&y^{3}+\Big(2 x^{3}+\frac{15}{7} \varepsilon  x \Big) y^{2}\\
&\quad +\Big(x^{6}+\frac{52}{21} \varepsilon  x^{4}+\frac{75}{49} \varepsilon^{2} x^{2}\Big) y \\
&\quad+\frac{\varepsilon}{3}x^{7} +\frac{314 \varepsilon^{2}}{441}x^{5} +\frac{125 \varepsilon^{3} }{343}x^{3}
\end{aligned}$} & $-3 x^{2}-\frac{25 \varepsilon}{7}$ & {$\begin{aligned}&C_{3,2}=\Big\{\beta = 0, \,\, \delta = 0,\\
&\quad \sigma = \frac{125 \varepsilon^{2}}{147},\,\, \alpha = \frac{40 \varepsilon}{21}\Big\}\end{aligned}$} \\[6pt]
~{$\begin{aligned}
&y^{3}+\Big(2 x^{3}-\frac{85}{27} \varepsilon  x \Big) y^{2}+\Big(x^{6}-\frac{76}{27} \varepsilon  x^{4}+\frac{1075}{729}  \varepsilon^{2}x^{2}\\
&\quad+\frac{400000}{531441} \varepsilon^{3}\Big) y+\frac{\varepsilon}{3}x^{7}-\frac{14 \varepsilon^{2}}{9}x^{5} +\frac{15875 \varepsilon^{3}}{6561} x^{3}\\
&\quad-\frac{2000000 \varepsilon^{4}}{1594323}x
\end{aligned}$} & $-3 x^{2}-\frac{25 \varepsilon}{27}$ & {$\begin{aligned}&C_{3,3}=\Big\{\beta = 0,\,\, \delta = 0,\\
&\quad \sigma = -\frac{125 \varepsilon^{2}}{81},\,\, \alpha = -\frac{20 \varepsilon}{27}\Big\}\end{aligned}$} \\[6pt]
\bottomrule
\end{tabular}\label{table:LD1}
\end{table}

\begin{proof}
Here we only present the detailed proof for the case $N=3$, the remaining two cases can be proved by using similar calculations and arguments. According to Proposition \ref{LD-Prop1}, we know that the expressions of the invariant algebraic curves and their cofactors of system \eqref{LD-mod2} are as follows:
\begin{equation}\label{LD-IACs}
\begin{split}
G(s,z)&=z^{3}+\sum_{i=0}^2\left(f_{i,3}s^{3}+f_{i,2}s^2+ f_{i,1}s+f_{i,0}\right) z^{i},\\
K(s,z)&=a_2s^2+a_1s+a_0,
\end{split}
\end{equation}
where $f_{i,j}$ and $a_i$ are undetermined coefficients. The associated polynomial set
\[\mathcal{P} := \{\operatorname{coef}(\mathcal{X}G - KG, \boldsymbol{t}) \mid \boldsymbol{t} \in \operatorname{term}(\mathcal{X}G - KG)\}\] contains 21 polynomials $\{P_1,\ldots,P_{21}\} $, which can be expressed as
\begin{align*}
P_1&=\delta  f_{1,0}+a_{0} f_{0,0},\\
P_2&=-\alpha  \varepsilon  f_{1,0}+\varepsilon^{2} f_{1,0}+\alpha  f_{0,1}+\delta  f_{1,1}+f_{1,0} \sigma -\varepsilon  f_{0,1}+a_{0} f_{0,1}+a_{1} f_{0,0}, \\
P_3&=2 f_{2,0} \delta +f_{1,0} \varepsilon +a_{0} f_{1,0}-f_{0,1}, \\
P_4&=-\frac{1}{2} \beta  \varepsilon  f_{1,0}-\alpha  \varepsilon  f_{1,1}+a_{1} f_{0,1}+\frac{1}{2} \beta  f_{0,1}+\varepsilon^{2} f_{1,1}+\sigma  f_{1,1}+\delta  f_{1,2}\\
&\quad{}+2 \alpha  f_{0,2}+a_{0} f_{0,2}-2 \varepsilon  f_{0,2}+a_{2} f_{0,0},\\
P_5&=-\frac{1}{2} \beta  \varepsilon  f_{1,1}-\alpha  \varepsilon  f_{1,2}+f_{0,1}+a_{1} f_{0,2}-3 \varepsilon  f_{0,3}+a_{2} f_{0,1}+\varepsilon^{2} f_{1,2}\\
&\quad{}+\sigma  f_{1,2}+\delta  f_{1,3}+\beta  f_{0,2}+3 \alpha  f_{0,3}+a_{0} f_{0,3},\\
P_6&=-f_{2,1}+3 \varepsilon +a_{0}, \\
P_7&=-\frac{1}{2} \beta  \varepsilon  f_{1,3}+3 f_{0,3}+a_{2} f_{0,3}, \\
P_8&=-\frac{1}{2} \beta  \varepsilon  f_{1,2}-\alpha  \varepsilon  f_{1,3}+2 f_{0,2}+a_{2} f_{0,2}+\varepsilon^{2} f_{1,3}+\sigma  f_{1,3}+\frac{3}{2} \beta  f_{0,3}+a_{1} f_{0,3},\\
P_9&=2 \varepsilon  f_{2,0}+a_{0} f_{2,0}+3 \delta -f_{1,1}, \\
P_{10}&=-3 f_{2,3}+a_{2}, \\
P_{11}&=a_{1}-2 f_{2,2}, \\
P_{12}&=f_{2,3} a_{2}+3 f_{2,3}, \\
P_{13}&=a_{1} f_{2,3}+a_{2} f_{2,2}+\frac{3}{2} \beta  f_{2,3}+2 f_{2,2}, \\
P_{14}&=-\beta  \varepsilon  f_{2,3}+a_{2} f_{1,3}+3 f_{1,3}, \\
P_{15}&=-\beta  \varepsilon  f_{2,2}-2 \alpha  \varepsilon  f_{2,3}+2 \varepsilon^{2} f_{2,3}+2 \sigma  f_{2,3}+a_{1} f_{1,3}+a_{2} f_{1,2}+\frac{3}{2} \beta  f_{1,3}+2 f_{1,2}, \\
P_{16}&=-2 \alpha  \varepsilon  f_{2,2}-\beta  \varepsilon  f_{2,1}+2 \varepsilon^{2} f_{2,2}+3 \alpha  f_{1,3}+\beta  f_{1,2}+2 \delta  f_{2,3}+2 \sigma  f_{2,2}\\
&\quad{}-2 \varepsilon  f_{1,3}+a_{0} f_{1,3}+a_{1} f_{1,2}+a_{2} f_{1,1}+f_{1,1}, \\
P_{17}&=-\beta  \varepsilon  f_{2,0}-2 \alpha  \varepsilon  f_{2,1}-\varepsilon  f_{1,2}+2 \varepsilon^{2} f_{2,1}+2 \sigma  f_{2,1}+2 \delta  f_{2,2}+a_{0} f_{1,2}\\
&\quad{}+a_{1} f_{1,1}+a_{2} f_{1,0}+\frac{1}{2} \beta  f_{1,1}+2 \alpha  f_{1,2}-3 f_{0,3}, \\
P_{18}&=-2 \alpha  \varepsilon  f_{2,0}+2 \varepsilon^{2} f_{2,0}+\alpha  f_{1,1}+2 \delta  f_{2,1}+2 \sigma  f_{2,0}+a_{0} f_{1,1}+a_{1} f_{1,0}-2 f_{0,2}, \\
P_{19}&=3 \alpha  f_{2,3}+\beta  f_{2,2}-\varepsilon  f_{2,3}+a_{0} f_{2,3}+a_{1} f_{2,2}+a_{2} f_{2,1}+f_{2,1},\\
P_{20}&=-3 \alpha  \varepsilon +\alpha  f_{2,1}+3 \varepsilon^{2}+\varepsilon  f_{2,1}+a_{0} f_{2,1}+a_{1} f_{2,0}+3 \sigma -2 f_{1,2},\\
P_{21}&=-\frac{3}{2} \beta  \varepsilon +a_{0} f_{2,2}+a_{1} f_{2,1}+a_{2} f_{2,0}+\frac{1}{2} \beta  f_{2,1}+2 \alpha  f_{2,2}-3 f_{1,3}.
\end{align*}
Note that $\varepsilon\neq0$, so we add one more constraint $1-\omega\varepsilon$ to the polynomial set $\mathcal{P}$, and obtain an ideal $\mathfrak{I}=\langle P_1,\ldots,P_{21},1-\omega\varepsilon \rangle$. Let $\mathfrak{I}_{16}$ be the 16th elimination ideal of the ideal $\mathfrak{I}$ in the ring $\mathbb{C}[\boldsymbol{\mu},\boldsymbol{f},\boldsymbol{a},\omega]$, i.e., eliminating the undetermined coefficients $\boldsymbol{f}=(f_{0,0},\ldots,f_{2,3}),\boldsymbol{a}=(a_0,a_1,a_2)$, and $\omega$. Computing the minimal prime decomposition of $\mathfrak{I}_{16}$, we obtain
$\sqrt{\mathfrak{I}_{16}}= \bigcap\limits_{i=1}^{4} \mathfrak{p}_i$, where
\begin{equation*}
\begin{split}
\mathfrak{p}_1&=\langle \beta, \delta, 2 \alpha -5 \varepsilon, -25 \varepsilon^{2}+18 \sigma \rangle,\\
\mathfrak{p}_2&=\langle \beta, \delta, 21 \alpha -40 \varepsilon,
-125 \varepsilon^{2}+147 \sigma
 \rangle,\\
\mathfrak{p}_3&=\langle \beta, \delta, 27 \alpha +20 \varepsilon,
125 \varepsilon^{2}+81 \sigma
 \rangle,\\
\mathfrak{p}_4&=\langle \beta, -\alpha  \varepsilon +\varepsilon^{2}+\sigma \rangle.
\end{split}
\end{equation*}
Thus, using Gr\"obner basis method we obtain necessary conditions for system \eqref{LD-mod2} to
have invariant algebraic surfaces of the form \eqref{LD-IACs}. Next, we will explain that by combining the triangular set method, we can not only verify whether the conditions are sufficient, but also derive the expressions of the invariant algebraic curves.

By adding the condition $\mathfrak{p}_1$ to the polynomial set $\mathcal{P}$, we compute the  triangular series under the constraint $\varepsilon\neq0$:
\[{\rm TriSer}([[\beta, \delta, 2 \alpha -5 \varepsilon, -25 \varepsilon^{2}+18 \sigma,P_1,\ldots,P_{21}],\{\varepsilon\}],[\boldsymbol{\mu},\boldsymbol{a},\boldsymbol{f}]);\]
the result contains the following one branch
\begin{equation*}
\begin{split}
[\mathcal{T}_0,\mathcal{U}_0]=&\,\,[[\beta, \delta, 2 \alpha^{2}-9 \sigma, 2 \alpha -5 \varepsilon,
5 \alpha +3 a_{0}, a_{1}, a_{2}+3, f_{0,0}, f_{0,1}, f_{0,2},
16 \alpha^{3}-3375 f_{0,3}, \\
&f_{1,0}, f_{1,1}, 8 \alpha^{2}-225 f_{1,2}
, f_{1,3}, f_{2,0}, 7 \alpha +15 f_{2,1}, f_{2,2}, 1+f_{2,3}], \{
\alpha\}].
\end{split}
\end{equation*}
Solving this branch, we obtain
\begin{equation*}
\begin{split}
a_{0} &= -\frac{5 \alpha}{3}, a_{1} = 0, a_{2} = -3, \\
f_{0,0} &= f_{0,1} = f_{0,2} = f_{1,0} = f_{1,1} = f_{1,3}= f_{2,0} = f_{2,2} = 0,\\
f_{0,3} &= \frac{16 \alpha^{3}}{3375},
f_{1,2} = \frac{8 \alpha^{2}}{225},
f_{2,1} = -\frac{7 \alpha}{15},
f_{2,3} = -1.
\end{split}
\end{equation*}
Substituting into \eqref{LD-IACs}, we obtain the expressions of the invariant algebraic curve and its cofactor of system \eqref{LD-mod2}:
\begin{equation*}
\begin{split}
G(s,z)=
z^{3}+\left(-s^{3}-\frac{7}{15} \alpha  s \right) z^{2}+\frac{8 \alpha^{2} s^{2} z}{225}+\frac{16 \alpha^{3} s^{3}}{3375},\quad K(s,z)=-3 s^{2}-\frac{25 \varepsilon}{6}.
\end{split}
\end{equation*}

Then, going back through the changes of variables \eqref{LD-change}, we have the following irreducible invariant algebraic curve and its cofactor of differential system \eqref{LD-mod1}:
\begin{equation*}
\begin{split}
F(x,y)&=y^{3}+\left(2 x^{3}+\frac{10}{3} \varepsilon  x \right) y^{2}+\left(x^{6}+\frac{11}{3} \varepsilon  x^{4}+\frac{125}{36} \varepsilon^{2} x^{2}\right) y +\frac{\varepsilon}{3}x^{7} +\frac{11 \varepsilon^{2}}{9}x^{5}+\frac{125 \varepsilon^{3}}{108}x^{3},\\
K(x,y)&=-3 x^{2}-\frac{25 \varepsilon}{6}.
\end{split}
\end{equation*}
This completes the proof of $C_{3,1}$ for the case $N=3$. In a similar manner, we can verify that the ideal $\mathfrak{p}_2$ corresponds to condition $C_{3,2}$ and the ideal  $\mathfrak{p}_3$ corresponds to condition $C_{3,3}$. For both cases, we further derive the associated irreducible invariant algebraic curves and their corresponding cofactors. However, for the ideal $\mathfrak{p}_4$, we find that the computation of the related triangular series contains the following three branches:
\begin{equation*}
\begin{split}
[\mathcal{T}_1,\mathcal{U}_1]=&\,\,[[\beta, \alpha  \varepsilon -\varepsilon^{2}-\sigma,
a_{0}+3 \varepsilon, a_{1}, a_{2}, -\varepsilon^{3} f_{0,0}+\delta^{3}
, f_{0,1}, f_{0,2}, f_{0,3}, \delta  f_{1,0}-3 \varepsilon  f_{0,0}, \\
&f_{1,1}, f_{1,2}, f_{1,3}, \delta  f_{2,0}-\varepsilon  f_{1,0},
f_{2,1}, f_{2,2}, f_{2,3}]
, \{\delta, \varepsilon\}],\\
[\mathcal{T}_2,\mathcal{U}_2]=&[[\beta, \delta, \alpha  \varepsilon -\varepsilon^{2}-\sigma,
a_{0}+3 \varepsilon, a_{1}, a_{2}, f_{0,0}, f_{0,1}, f_{0,2}, f_{0,3}
, f_{1,0}, f_{1,1}, f_{1,2},\\
&\quad f_{1,3}, f_{2,0}, f_{2,1}, f_{2,2},
f_{2,3}], \{\varepsilon\}],\\
[\mathcal{T}_3,\mathcal{U}_3]=&[[\beta, \delta, 3 \alpha^{2}-16 \sigma, 3 \alpha -4 \varepsilon,
4 a_{0}+9 \alpha, a_{1}, a_{2}+3, f_{0,0}, f_{0,1}, f_{0,2}, f_{0,3},
f_{1,0}, f_{1,1}, \\
&f_{1,2}, f_{1,3}, f_{2,0}, f_{2,1}, f_{2,2},
f_{2,3}+1], \{\alpha\}].
\end{split}
\end{equation*}
Solving the first branch, we obtain the associated invariant algebraic curve and the corresponding cofactor:
\[F(x,y)=\left(y+x^{3}+\left(\alpha -\varepsilon \right) x +\frac{\delta}{\varepsilon}\right)^3,\quad K(x,y)=-3\varepsilon.\]
So the irreducible invariant algebraic curve is $F(x,y)=y+x^{3}+\left(\alpha -\varepsilon \right) x +\frac{\delta}{\varepsilon}$ with cofactor $-\varepsilon$, which corresponds to condition $C_{1,1}=\{\beta=0,\,\,\sigma = \varepsilon  \left(\alpha -\varepsilon \right)\}$. Solving the second branch, we obtain the associated invariant algebraic curve and the corresponding cofactor:
\[F(x,y)=\left(y+x^{3}+\left(\alpha -\varepsilon \right) x\right)^3,\quad K(x,y)=-3\varepsilon.\]
So the irreducible invariant algebraic curve is $F(x,y)=y+x^{3}+\left(\alpha -\varepsilon \right) x$ with cofactor $-\varepsilon$, which corresponds to the condition $\{\beta=0,\,\, \delta=0,\,\,\sigma = \varepsilon  \left(\alpha -\varepsilon \right)\}\subset C_{1,1}$. Solving the third branch, we obtain the associated invariant algebraic curves and the corresponding cofactor:
\[F(x,y)=\left(y+\frac{\varepsilon}{3}x\right)\left(y+x^{3}+\frac{\varepsilon}{3}x\right)^2,\quad K(x,y)=-3 x^{2}-3\varepsilon.\]
So the irreducible invariant algebraic curves are $F(x,y)=y+\frac{\varepsilon}{3}x$ with cofactor $-3 x^{2}-\varepsilon$ and $F(x,y)=y+x^{3}+\frac{\varepsilon}{3}x$ with cofactor $-\varepsilon$, which correspond to the condition $C_{1,1} \wedge C_{1,2}=\{\beta=0,\,\,\delta=0,\,\,\alpha=\frac{4}{3}\varepsilon,\,\,\sigma = \frac{1}{3}\varepsilon^2\}$. This confirms that the ideal $\mathfrak{p}_4$ corresponds to the conditions given in the case $N=1$. We complete the proof for the case $N=3$.

\end{proof}

\begin{remark}\label{hand-error}
The condition $C_{1,2}$ is different from the one obtained in \cite[Table 1]{Demina2018}. One can check that, in the second row of \cite[Table 1]{Demina2018}, $\sigma$ is written incorrectly as $\sigma=-\frac{\varepsilon}{18}(\beta^2+2\varepsilon-6\alpha)$. In addition, there are a few typos in $F(x,y)$, for example, under the condition $C_{3,2}$, the term $\frac{314\varepsilon^2}{551}x^5$ (see the second to last expression in \cite[Table 1]{Demina2018}) should be $\frac{314\varepsilon^2}{441}x^5$. This indicates that our algorithmic approach allows us to detect hand-computation errors in published literature. Note that Theorem \ref{LD-main} gives the classification of irreducible invariant algebraic curves of system \eqref{LD-mod1}. What we would like to emphasize is that in some cases, the form of irreducible invariant algebraic curves is not unique. In fact, under the condition $C_{1,1} \wedge C_{1,2}$ (this corresponds to the ideal $\mathfrak{p}_4$), system \eqref{LD-mod1} has two types of irreducible invariant algebraic curves.

\end{remark}

\section{Hopf Bifurcation and Lyapunov Coefficients}\label{Hopf-sect4}

In this section, we present some basic notions about Hopf bifurcation and Lyapunov coefficients. The term Hopf bifurcation refers to the local birth or death of a periodic solution from an equilibrium as a parameter crosses a critical value. The Hopf bifurcation theorem of a conservative planar system was known to Andronov and his co-workers from around 1930 \cite{Jackson1989}. The occurrence of such bifurcations in general $n$-dimensional systems for $n>2$ was studied by Hopf in the 1940s \cite{MJMM1976,Jackson1989}. These bifurcations are now named after him as Hopf bifurcation. The theory of Lyapunov coefficients can be found in \cite[Chapters 3 and 10]{Kuznetsov2004}. The goal of this section is to develop necessary and sufficient conditions for Hopf bifurcation points and efficient algorithm to compute Lyapunov coefficients.

We use the notation in \cite{GJHP1983} with minor changes. Consider the following $n$-dimensional nonlinear polynomial differential system:
\begin{equation}\label{L-dot_x}
\dot{\boldsymbol{x}} = f(\boldsymbol{x},\boldsymbol{\mu}),\quad f:\mathbb{R}^n\times\mathbb{R}^m\rightarrow\mathbb{R}^n,
\end{equation}
where $\boldsymbol{x}=(x_1,\ldots,x_n)$ is a list of variables, $\boldsymbol{\mu}=(\mu_1,\ldots,\mu_m)$ is a list of real parameters, and $f(\boldsymbol{x},\boldsymbol{\mu})=(f_1,\ldots,f_n)$ is a list of polynomials in $\mathbb{R}[\boldsymbol{x},\boldsymbol{\mu}]$. According to the well-known Hopf bifurcation theorem \cite{MJMM1976,GJHP1983,Kuznetsov2004}, sufficient conditions for this type of bifurcation are:

(a) {\bf Existence of a Hopf bifurcation point (HBP)}. A Hopf bifurcation point $(\bar{\boldsymbol{x}},\bar{\boldsymbol{\mu}})$ is an equilibrium point of the vector field $f(\boldsymbol{x},\boldsymbol{\mu})$ where the Jacobian matrix $D_{\boldsymbol{x}}f(\bar{\boldsymbol{x}},\bar{\boldsymbol{\mu}})$ has a simple pair of purely imaginary eigenvalues $(\lambda(\bar{\boldsymbol{x}},\bar{\boldsymbol{\mu}})=i\omega_0,\bar{\lambda}(\bar{\boldsymbol{x}},\bar{\boldsymbol{\mu}})=-i\omega_0,\omega_0>0)$, and no other eigenvalue with zero real part.

(b) {\bf Transversality condition}. Assume that there is a smooth curve of equilibria $(\boldsymbol{x}_*(\boldsymbol{\mu}),\boldsymbol{\mu})$ around $\bar{\boldsymbol{\mu}}$, that is, $f(\boldsymbol{x}_*(\boldsymbol{\mu}),\boldsymbol{\mu})=0$ for all $\boldsymbol{\mu}$ close enough to $\bar{\boldsymbol{\mu}}$ such that $\boldsymbol{x}_*(\bar{\boldsymbol{\mu}})=\bar{\boldsymbol{x}}$. Let $\lambda(\boldsymbol{\mu})$ be the eigenvalue of $D_{\boldsymbol{x}} f(\boldsymbol{x}_*(\boldsymbol{\mu}),\boldsymbol{\mu})$ satisfying $\lambda(\bar{\boldsymbol{\mu}})=i\omega_0$. With the other parameters held fixed, transversality with respect to the chosen parameter $\mu_i$ means:
\[\frac{\partial {\rm Re} (\lambda(\boldsymbol{\mu}))}{\partial \mu_i}\Big|_{\mu=\bar{\mu}}\neq0.\]

(c) {\bf Genericity condition (or nondegeneracy condition)}. The first Lyapunov coefficient $\ell_1(\bar{\boldsymbol{x}})$ is nonvanishing.

Condition (a) gives an equilibrium $(\bar{\boldsymbol{x}},\bar{\boldsymbol{\mu}})$ whose Jacobian $J_f=D_{\boldsymbol{x}}f(\bar{\boldsymbol{x}},\bar{\boldsymbol{\mu}})$ has no zero eigenvalue. The Implicit Function
Theorem therefore gives a unique local function $\boldsymbol{x}_*(\boldsymbol{\mu})$ satisfying $f(\boldsymbol{x}_*(\boldsymbol{\mu}),\boldsymbol{\mu})=0$ and $\boldsymbol{x}_*(\bar{\boldsymbol{\mu}})=\bar{\boldsymbol{x}}$. Simplicity of $i\omega_0$ gives a locally differentiable eigenvalue $\lambda(\boldsymbol{\mu})$ of $D_xf(\boldsymbol{x}_*(\boldsymbol{\mu}),\boldsymbol{\mu})$ with $\lambda(\bar{\boldsymbol{\mu}})=i\omega_0$. When the three conditions are satisfied, a unique local branch of periodic orbits with periods tending to $2\pi/\omega_0$ is born at $(\bar{\boldsymbol{x}},\bar{\boldsymbol{\mu}})$, and the sign of $\ell_1(\bar{\boldsymbol{x}})$ gives stability of the small bifurcated limit cycle \textit{on the center manifold}. The parameter values that satisfy the above conditions correspond to the critical points for a Hopf bifurcation; at these critical points, a limit cycle has not yet emerged. A limit cycle can only be observed once the parameter value crosses this critical point (see \cite{GJHP1983,Kuznetsov2004}). We refer to the type of Hopf bifurcation satisfying (a) and (b) as a \textit{complete} Hopf bifurcation ({\bf CHB}) in order to distinguish it from the simple Hopf bifurcations introduced by Liu \cite{JMAA:1994}. Whether it is supercritical, subcritical, or degenerate depends on the higher order terms in system \eqref{L-dot_x}. It remains to compute the Lyapunov coefficients of the system near $\bar{\boldsymbol{x}}$. If the genericity condition is violated (i.e., $\ell_1(\bar{\boldsymbol{x}})=0$), then the system is degenerate (often leading to a Bautin bifurcation or generalized Hopf bifurcation of codimension 2). In this case, multiple families of periodic solutions or complex phenomena involving multiple limit cycles may bifurcate, requiring higher-order derivative analysis. Here we recall the projection method presented in \cite{Kuznetsov2004} for the calculation of the first Lyapunov coefficient associated to Hopf bifurcations.

Let $\dot{\boldsymbol{x}}=F(\boldsymbol{x})$ be a differential system having $\boldsymbol{x}_0=0$ as an equilibrium point (otherwise, we can translate the equilibrium to the origin). Consider the third-order Taylor approximation of $F$ around $\boldsymbol{x}_0=0$ given by
\begin{equation} \label{L-F-talor}
\begin{split}
F(\boldsymbol{x}) = &A\boldsymbol{x} + \frac{1}{2} B(\boldsymbol{x}, \boldsymbol{x}) + \frac{1}{6} C(\boldsymbol{x}, \boldsymbol{x}, \boldsymbol{x})+ \mathcal{O}(\|\boldsymbol{x}\|^4),
\end{split}
\end{equation}
where $A = F_{\boldsymbol{x}}(0)$ and $B(\boldsymbol{x},\boldsymbol{y})$ and $C(\boldsymbol{x},\boldsymbol{y},\boldsymbol{z})$ are multilinear functions with components
\begin{align}
&B_i(\boldsymbol{x},\boldsymbol{y}) = \sum_{j, k=1}^n \left. \frac{\partial^2 F_i(\xi)}{\partial \xi_j \partial \xi_k} \right|_{\xi=0} x_j y_k, \label{L-B}\\
&C_i(\boldsymbol{x},\boldsymbol{y},\boldsymbol{z}) = \sum_{j, k, \ell=1}^n \left. \frac{\partial^3 F_i(\xi)}{\partial \xi_j \partial \xi_k \partial \xi_{\ell}} \right|_{\xi=0} x_j y_k z_{\ell},\label{L-C}
\end{align}
where $i = 1,\,2,\,...\,,\,n$. The following result can be found in the book \cite[pp. 177--181]{Kuznetsov2004}.

\begin{theorem}\label{Hopf-th}
Suppose that $A$ has a simple pair of purely imaginary eigenvalues $\pm i \omega_0, \omega_0>0$, and no other eigenvalue with zero real part. Let $p,\,q \in \mathbb{C}^n$ be two complex eigenvectors such that
\begin{equation}\label{L_Aq}
A q = i \omega_0 q, \quad A^\top p = -i \omega_0 p, \quad \langle p, q \rangle = \sum_{i=1}^n \bar{p}_i q_i = 1,
\end{equation}
where $A^\top$ is the transposed matrix of $A$. If $I_n$ denotes the identity matrix, then we can define the first Lyapunov coefficient $\ell_1(\boldsymbol{x}_0)$ as
\begin{equation}\label{L-l1}
\ell_1(\boldsymbol{x}_0) = \frac{1}{2} \,{\rm{Re}} \, G_{21},
\end{equation}
where
\begin{equation*}% \label{L-G21}
G_{21} = \left\langle p,\, C(q, q, \bar{q}) + B(\bar{q}, (2i\omega_0 I_n - A)^{-1}B(q, q)) - 2B(q, A^{-1}B(q, \bar{q})) \right\rangle.
\end{equation*}

\end{theorem}

A {\bf HBP} that satisfies conditions (b) and (c) is called the Hopf point of codimension 1, the sign of $\ell_1(\boldsymbol{x}_0)$ gives stability of the small bifurcated periodic orbit on the
center manifold. When $\ell_1(\boldsymbol{x}_0)< 0$, the point $\boldsymbol{x}_0$ is a weak focus of the system restricted to the center manifold of $\boldsymbol{x}_0$ and the limit cycle that emerges from $\boldsymbol{x}_0$ is stable. In this case, we say that the Hopf bifurcation is supercritical. When $\ell_1(\boldsymbol{x}_0)>0$, the point $\boldsymbol{x}_0$ is also a weak focus of the system restricted to the center manifold of $\boldsymbol{x}_0$, but the limit cycle that bifurcates from $\boldsymbol{x}_0$ is unstable. In this second case, we say that the Hopf bifurcation is subcritical. More detailed discussions on Hopf bifurcation, including applications, can be found in \cite{Kuznetsov2004}.

The Algorithm \ref{alg:lyapunov}, presented below, is based on Theorem \ref{Hopf-th} that can be used to compute the first Lyapunov coefficient $\ell_1(\boldsymbol{x}_0)$ for a given differential system.

\begin{algorithm}
	\renewcommand{\algorithmicrequire}{\textbf{Input:}}
	\renewcommand{\algorithmicensure}{\textbf{Output:}}
\caption{Symbolic-numeric computation of the first Lyapunov coefficient}\label{alg:lyapunov}
\begin{algorithmic}[1]
\REQUIRE{A differential system \eqref{L-dot_x} with a {\bf HBP} $(\boldsymbol{x}_0,\boldsymbol{\mu}_0)$}
\ENSURE{The first Lyapunov coefficient $\ell_1(\boldsymbol{x}_0)$}

\STATE Evaluate the differential system using the parameter condition $\boldsymbol{\mu}_0$ and translate the point $\boldsymbol{x}_0$ to the origin.

%\STATE Evaluate $\mathbf{x}_0$ numerically using floating-point arithmetic with significant digits.

\STATE Evaluate the matrix $A$ in \eqref{L-F-talor} at the point $\boldsymbol{x}_0$ (using floating-point arithmetic with significant digits).

\STATE Compute the complex eigenvectors $p$ and $q$ using \eqref{L_Aq} and convert them to approximate rational expressions.

\STATE Compute the multilinear functions $B(\boldsymbol{x}, \boldsymbol{y})$ and $C(\boldsymbol{x},\boldsymbol{y},\boldsymbol{z})$ in \eqref{L-F-talor} using \eqref{L-B} and \eqref{L-C}.

\STATE Compute the first Lyapunov coefficient $\ell_1(\boldsymbol{x}_0)$ using \eqref{L-l1}.

\RETURN $\ell_1(\boldsymbol{x}_0)$

\end{algorithmic}
\end{algorithm}

\begin{remark}
The termination of the algorithm is obvious. The
correctness of Algorithm \ref{alg:lyapunov} follows directly from Theorem \ref{Hopf-th}. To enhance computational precision, Line 1 is the algorithm preprocessing stage, which restricts numerical evaluation exclusively to the parameters $\boldsymbol{\mu}$, while $\boldsymbol{x}_0$ is retained in full symbolic form. Line 2 then proceeds to evaluate the matrix $A$ at the point $\boldsymbol{x}_0$ using floating-point arithmetic. We implemented the algorithm in \textsc{Maple} (from line 2 to line 6). The \textsc{Maple} program is available at {\small{
{\textcolor{blue}{\underline{\url{https://github.com/Bo-Math/Computing-Lyapunov-Coefficient}}}}}}. Before running Algorithm \ref{alg:lyapunov}, we will first use the routine {\sf SamplePoints} (within the {\sf RegularChains} library in \textsc{Maple}) to derive several Hopf bifurcation points for a differential system to exhibit Hopf bifurcations. In Section \ref{sect-5}, we will apply our \textsc{Maple} program to analyze the nature of Hopf bifurcations for certain differential systems in order to show its feasibility.
\end{remark}

\subsection{Algebraic characterization of Hopf bifurcation points}\label{HBP-alg}

We are interested in detecting curves of {\bf HBP}s. This constitutes the location of the values $(\boldsymbol{x},\boldsymbol{\mu})$ at which $f(\boldsymbol{x},\boldsymbol{\mu})=0$ and $Df(\boldsymbol{x},\boldsymbol{\mu})$ has a simple pair of purely imaginary eigenvalues. We provide a complete algebraic characterization of {\bf HBP}s from different perspectives (one is based on the properties of Hurwitz determinants and subresultant sequences, and the other depends on the properties of the greatest common divisor of polynomials). The established results supplement and improve the conclusions about {\bf HBP}s in \cite{Weber2000,ChenMM-2011}.

\subsubsection{Approach 1: Characterization via Hurwitz determinants and subresultants}

The following result is a refined version of \cite[Theorem 10]{ChenMM-2011} for the case $k_1=0$ and $k_2=1$, which gives the necessary and sufficient conditions for a polynomial $p(\lambda)$ to have a simple pair of purely imaginary roots while no other roots have zero real parts. That is, the necessary and sufficient conditions for {\bf HBP}s of a differential system of the form \eqref{L-dot_x}.

\begin{theorem}\label{Hopf-th3}
Let $p(\lambda)=\lambda^n+a_1(\boldsymbol{x},\boldsymbol{\mu})\lambda^{n-1}+\cdots+a_n(\boldsymbol{x},\boldsymbol{\mu})$ be the characteristic polynomial of the Jacobian matrix $Df(\boldsymbol{x},\boldsymbol{\mu})$. We write $p(\lambda)$ in the form: $p(\lambda)=p_1(\lambda^2)+\lambda p_2(\lambda^2)$. Then $p(\lambda)$ has a simple pair of purely imaginary roots, $i\omega_0$ and $-i\omega_0$, with $\omega_0>0$, and no other root with zero real part if and only if one of the following statements holds:
\begin{itemize}
   \item[(a)]

   The condition $\bigcup\limits_{k\geq1}^{\lfloor \frac{n}{2} \rfloor}\mathcal{C}_k$ holds, where $\mathcal{C}_k$ is the condition on the coefficients of $p(\lambda)$ satisfy $\Delta_{n-1}=\Delta_{n-3}=\cdots=\Delta_{n-2k+1}=0$ and $a_n\Delta_{n-2k-1}\neq0$, and that $S_k(p_1,p_2,\lambda)$ has exactly one root in $(-\infty,0)$, which is simple.

  \item[(b)] The condition $\bigcup\limits_{k\geq1}^{\lfloor \frac{n}{2} \rfloor}\bar{\mathcal{C}}_k$ holds, where $\bar{\mathcal{C}}_k$ is the condition on the coefficients of $p(\lambda)$ such that
\begin{equation}\label{Sk-1}
\bar{S}_k=\left\{\begin{array}{l}
S_{k}(p_1,\lambda p_2,\lambda),\quad \text{when $n=2m+1$}\\
S_{k}(\lambda p_2,p_1,\lambda),\quad \text{when $n=2m$}\\
\end{array}\right.
\end{equation} has exactly one root in $(-\infty,0)$, and that root is simple; moreover, $\Delta_{n}=\Delta_{n-2}=\cdots=\Delta_{n-2k+2}=0$ and $a_n\Delta_{n-2k}\neq0$.
\end{itemize}

%The characteristic polynomial $p(\lambda)$ of the Jacobian matrix $Df(\mathbf{x},\mu)$ has a pair of distinct roots, $i\omega_0$ and $-i\omega_0$, on the imaginary axis, and no other root with zero real part if and only if the condition $\mathcal{C}=\bigcup\limits_{k\geq1}^{\lfloor \frac{n}{2} \rfloor}\bar{\mathcal{C}}_k$ holds, where $\bar{\mathcal{C}}_k$ are the conditions expressed in Theorem \ref{Hopf-th1} or Theorem \ref{Hopf-th2}.

\end{theorem}

\begin{proof}
Here we only present the proof for the statement (a), the proof of the statement (b) is similar. The only difference is
that during the proof we need to use Corollary \ref{lea-HS2} instead of Corollary \ref{lea-HS1}.

%The proof can be done using the arguments in \cite{ChenMM-2011}.

%By the assumption, given a polynomial $p(\lambda)$, we can construct the $i$th subresultant of $p_1$ and $p_2$ using \eqref{subequ1} for $i\in\{0,1,\ldots,\delta-1\}$.

\begin{itemize}
\item[``$\Rightarrow$'']
Assume that $p(\lambda)$ has a simple pair of purely imaginary roots, $i\omega_0$ and $-i\omega_0$, on the imaginary axis and no other root with zero real part. One then has $a_n\neq0$ and we can write $p(\lambda)=p_1^*(\lambda)(\lambda^2+\omega_0^2)$, where $\omega_0^2>0$ and $p_1^*(\lambda)$ does not have root with zero real part. Note that $p(\lambda)$ may have other symmetric roots, we denote them as $\pm z_1,\ldots,\pm z_j$ with ${\rm Re}(z_1)\cdots{\rm Re}(z_j)\neq0$ for $j\geq0$ (if $j=0$, we mean that $p(\lambda)$ does not have other symmetric roots). Therefore, $p(\lambda)$ can be written as
\[p(\lambda)=p_2^*(\lambda)(\lambda^2+\omega_0^2)(\lambda^2-z_1^2)\cdots(\lambda^2-z_j^2),\]
where $p_2^*(\lambda)$ does not have symmetric roots. When some of the $\{z_1^2,\ldots,z_j^2\}$ are equal, the multiplicity of these symmetric roots is greater than 1. By assumption we know that $p(\lambda)$ in total has exactly $1+j$ pairs of nonzero symmetric roots. By Corollary \ref{lea-HS1}, we have $\Delta_{n-1}=0,\ldots,\Delta_{n-2j-1}=0,\Delta_{n-2j-3}\neq0$. On the other hand, by Theorem \ref{HSres} we have $s_0(p_1,p_2,\lambda)=0,\ldots,s_j(p_1,p_2,\lambda)=0$, $s_{j+1}(p_1,p_2,\lambda)\neq0$, which implies that $S_{j+1}(p_1,p_2,\lambda)={\rm gcd}(p_1,p_2,\lambda)$. Since $\pm i\omega_0,\pm z_1,\ldots,\pm z_j$ are nonzero symmetric roots
of $p(\lambda)$, by Corollary \ref{lea-HS1}, $-\omega_0^2,z_1^2,\ldots,z_j^2$ are the common roots of $p_1(\lambda)$ and $p_2(\lambda)$, that is, they are the roots of $S_{j+1}(p_1,p_2,\lambda)$. In summary, take $k=j+1\in\{1,\ldots,\lfloor \frac{n}{2} \rfloor\}$, we have $S_k(p_1,p_2,\lambda)$ has exactly one negative real root ($-\omega_0^2$), and this root is simple; moreover,
\[\Delta_{n-1}=\Delta_{n-3}=\cdots=\Delta_{n-2k+1}=0,\quad a_n\Delta_{n-2k-1}\neq0.\]

\item[``$\Leftarrow$''] For any integer $k\in\{1,\ldots,\lfloor \frac{n}{2} \rfloor\}$, by the assumption, we have
$\Delta_{n-1}=\Delta_{n-3}=\cdots=\Delta_{n-2k+1}=0, a_n\Delta_{n-2k-1}\neq0$. It follows from Theorem \ref{HSres} that
\[s_0(p_1,p_2,\lambda)=s_1(p_1,p_2,\lambda)=\cdots=s_{k-1}(p_1,p_2,\lambda)=0, s_{k}(p_1,p_2,\lambda)\neq0.\]
Therefore the degree of $S_k(p_1,p_2,\lambda)$ is $k$ and $S_{k}(p_1,p_2,\lambda)={\rm gcd}(p_1,p_2,\lambda)$. Note that, by assumption, $p(\lambda)$ has exactly $k$ pairs of nonzero symmetric roots, denoted by $\pm z_1,\ldots,\pm z_k$. So $p(\lambda)$ can be written as $p(\lambda)=p_3^*(\lambda)(\lambda^2-z_1^2)\cdots(\lambda^2-z_k^2)$, where $p_3^*(\lambda)$ does not have symmetric roots. When some of the $\{z_1^2,\ldots,z_k^2\}$ are equal, the multiplicity of these symmetric roots is greater than 1. By Corollary \ref{lea-HS1}, $z_1^2,\ldots,z_k^2$ are common roots of $p_1(\lambda)$ and $p_2(\lambda)$. As a consequence, $z_1^2,\ldots,z_k^2$ are roots of $S_k(p_1,p_2,\lambda)$. On the other hand, $S_k(p_1,p_2,\lambda)$ has exactly one negative real root (which is simple), which implies that only one of the $k$ roots $z_1^2,\ldots,z_k^2$ is negative, say $z_1^2=-\omega_0^2<0$, and ${\rm Re}(z_j)\neq0$ for $j=2,\ldots,k$ (otherwise, $S_k(p_1,p_2,\lambda)$ will have another negative real root, which is a contradiction). Overall, $p(\lambda)$ has a simple pair of purely imaginary roots, $i\omega_0$ and $-i\omega_0$, and no other root with zero real part.

\end{itemize}

\end{proof}

The following result is a direct consequence of Theorem \ref{Hopf-th3}, which describes the equivalent conditions for the polynomial $p(\lambda)$ to have $\lfloor \frac{n}{2} \rfloor$ pairs of nonzero symmetric roots, of which exactly one simple pair of purely imaginary roots. This is precisely the situation that was not covered in \cite[Theorem 5]{ChenMM-2011}.
\begin{corollary}\label{coro-th-hp3}
Let $p(\lambda)=\lambda^n+a_1(\boldsymbol{x},\boldsymbol{\mu})\lambda^{n-1}+\cdots+a_n(\boldsymbol{x},\boldsymbol{\mu})$ be the characteristic polynomial of the Jacobian matrix $Df(\boldsymbol{x},\boldsymbol{\mu})$. When $n=2m$, the following statements are equivalent:
\begin{itemize}
    \item[(i)] $\bar{p}(y)=y^m+a_2y^{m-1}+\cdots+a_{2m}$ has exactly one root in $(-\infty,0)$, which is simple; moreover, $a_{2m}\neq0,\Delta_{2m-1}=\Delta_{2m-3}=\cdots=\Delta_{1}=0$.
   \item[(ii)]$\bar{p}(y)$ has exactly one root in $(-\infty,0)$, which is simple; moreover, $a_{2m}\neq0,\Delta_{2m}=\Delta_{2m-2}=\cdots=\Delta_{2}=0$.
 \item[(iii)] $\bar{p}(y)$ admits a unique root in $(-\infty,0)$, which is simple; moreover, $a_{2m}\neq0,a_1=a_3=\cdots=a_{2m-1}=0$.
\end{itemize}
Moreover, when $n=2m+1$, the following statements are equivalent:
\begin{itemize}
    \item[(a)] $p(\lambda)$ has $m$ pairs of nonzero symmetric roots (counted with multiplicity) and $p(0)\neq0$, of which exactly one simple pair of purely imaginary roots.
   \item[(b)]$\bar{p}(\lambda)$ has exactly one root in $(-\infty,0)$, which is simple; moreover, $a_1a_{2m}\neq0,a_3-a_1a_2=\cdots=a_{2m+1}-a_1a_{2m}=0$.
\end{itemize}

\end{corollary}

\begin{proof}
We will divide the proof into two parts. Part (A): when $n=2m$.

As stated in Remark \ref{rem-cor-HS2}, when $n=2m$,  neither $p_2(\lambda)$ nor $\lambda p_2(\lambda)$ can have $m$ pairs of nonzero symmetric roots $\{z_i,-z_i\}$ since $p_2(\lambda)$ is a polynomial of degree $m-1$. To make $p(\lambda)$ have $m$ pairs of symmetric roots, $p_2(\lambda)$ must be identically zero. In this case, we have $a_1=a_3=\cdots=a_{2m-1}=0$. At this point, we also have $S_{m}(p_1,p_2,\lambda)=S_{m}(\lambda p_2,p_1,\lambda)=p_1(\lambda)=\lambda^m+a_2\lambda^{m-1}+\cdots+a_{2m}=\bar{p}(\lambda)$. As a consequence, the equivalence of (i) and (ii) follows from the statements (a) and (b) of Theorem \ref{Hopf-th3} for $k=m$.

Now we prove the equivalence of (i) and (iii). Let $f(y)=y^n+a_1y^{n-1}+\cdots+a_{2m}$ be a polynomial of degree $2m$. We first establish the following identity: If $a_1=a_3=\cdots=a_{2k-3}=0$, then we have
\begin{equation}\label{Hur-claim}
\Delta_{2k-1}
=(-1)^{k(k-1)/2}a_{2k-1}^{k}.
\end{equation}

Note that the Hurwitz matrix defining $\Delta_{2k-1}$ has entries $a_{2j-i}$ (coefficients with indices out of range are treated as 0). Reorder its rows such that the even-numbered rows come first, followed by the odd-numbered rows, i.e.,
\[
2,4,\ldots,2k-2;1,3,\ldots,2k-1.
\]
Under the assumed vanishing of the odd-indexed coefficients, the resulting matrix has the block form
\[
\begin{pmatrix}
U & *\\
0 & V
\end{pmatrix},
\]
where $U$ is an upper triangular matrix of order $k-1$ with every diagonal entry equal to 1, and $V$ is an upper triangular matrix of order $k$ with every diagonal entry equal to $a_{2k-1}$. The row permutation has $k(k-1)/2$ inversions, which proves the identity \eqref{Hur-claim}. For $k=1$, this identity simply reads $\Delta_1=a_1$.

If $a_1=a_3=\cdots=a_{2m-1}=0$, applying \eqref{Hur-claim} for $k=1,\ldots,m$ gives $\Delta_1=\Delta_3=\cdots=\Delta_{2m-1}=0$. Conversely, suppose that $\Delta_1=\Delta_3=\cdots=\Delta_{2m-1}=0$. We proceed by induction. First, $a_1=\Delta_1=0$. If $a_1=a_3=\cdots=a_{2k-3}=0$, then \eqref{Hur-claim} gives
\[
0=\Delta_{2k-1}
=(-1)^{k(k-1)/2}a_{2k-1}^{k}.
\]
So we have $a_{2k-1}=0$. The conclusion follows by induction.

Part (B): when $n=2m+1$.
\begin{itemize}
\item[``$\Rightarrow$''] In this case, $p(\lambda)$ can be written in the form:
\begin{equation*}% \label{proof-HB3}
\begin{split}
p(\lambda)&=\lambda^{2m+1}+a_1\lambda^{2m}+a_2\lambda^{2m-1}+\cdots+a_{2m}\lambda+a_{2m+1}\\
&=(\lambda+\lambda_0)(\lambda^{2m}+b_1\lambda^{2m-1}+b_2\lambda^{2m-2}+\cdots+b_{2m-1}\lambda+b_{2m}).
\end{split}
\end{equation*}
By the assumption and Part (A), we know that $\lambda_0\neq0$ and $b_1=b_3=\cdots=b_{2m-1}=0,b_{2m}\neq0$. Moreover, we have $\bar{p}(y)=y^m+b_2y^{m-1}+\cdots+b_{2m}$ has one simple negative real root. By comparing the coefficients of $p(\lambda)$, we obtain
\begin{equation}\label{proof-HB4}
\begin{split}
a_1=\lambda_0,\,a_2=b_2,\,a_3=\lambda_0b_2,\,\ldots,\,a_{2k}=b_{2k},\,a_{2k+1}=\lambda_0b_{2k},\quad k=1,\ldots,m.
\end{split}
\end{equation}
The desired result follows directly from \eqref{proof-HB4}.

\item[``$\Leftarrow$''] Write $R(y)=y^m+a_2y^{m-1}+\cdots+a_{2m}$. By the assumption, the polynomial $p(\lambda)$ can be written as
\begin{equation*}% \label{proof-HB5}
\begin{split}
p(\lambda)&=\lambda^{2m+1}+a_1\lambda^{2m}+a_2\lambda^{2m-1}+a_1a_2\lambda^{2m-2}+\cdots+a_{2m}\lambda+a_1a_{2m}\\
&=\lambda^{2m}(\lambda+a_1)+a_2\lambda^{2m-2}(\lambda+a_1)+\cdots+a_{2m}(\lambda+a_1)\\
&=(\lambda+a_1)R(\lambda^2).
\end{split}
\end{equation*}
Since $\bar{p}(y)=R(\lambda^2)$ has one simple negative real root and $a_{2m}\neq0$, so by Part (A), we know that $p(\lambda)$ has $m$ pairs of nonzero roots $\{z_i,-z_i\}$ in $\mathbb{C}$, of which exactly one pair of pure imaginary roots. In addition, $p(\lambda)$ has one nonzero real root $-a_1$.

\end{itemize}

\end{proof}

We remark that, Corollary \ref{coro-th-hp3} provides an alternative way to deduce simple equivalent conditions for $p(\lambda)$ to have $\lfloor \frac{n}{2} \rfloor$ pairs of nonzero symmetric roots, of which exactly one simple pair of purely imaginary roots. Table \ref{tab:Hopf-1} provides a list of the necessary and sufficient conditions for {\bf HBP}s, as functions of polynomial coefficients, for vector fields of dimension two to five. As far as we know we are the first to provide these expressions.

\begin{table}[H]
	\centering
	\caption{Necessary and sufficient conditions for {\bf HBP}s for vector fields of dimension two through five.}
	\label{tab:Hopf-1}
\begin{tabular}{c|c}
  \hline \hline
  $n$ & Conditions for {\bf HBP}s \\
  \hline
 $2$ & $[\Delta_1=a_1=0,\quad a_2>0]$ \\
  \hline
 $3$ & $\left[\Delta_2=\left|
 \begin{matrix}
    a_1& a_3 \\
   1 & a_2 \\
  \end{matrix}
  \right|=0,\quad a_1a_3>0\right]$ \\
  \hline
 $4$ & $\left[\Delta_3=\left|
 \begin{matrix}
    a_1& a_3 & 0\\
   1 & a_2 & a_4 \\
   0 & a_1 & a_3 \\
  \end{matrix}
  \right|=0,\quad a_1a_4\Delta_2>0\right]$ $\cup$ $[a_1=a_3=0,\quad a_4<0]$\\
  \hline
 $5$ & $\left[\Delta_4=\left|
 \begin{matrix}
    a_1& a_3 & a_5 & 0\\
   1 & a_2 & a_4 &0\\
   0 & a_1 & a_3 & a_5\\
   0 & 1 & a_2 & a_4
  \end{matrix}
  \right|=0,\quad a_5\Delta_3\Delta_2>0\right]$ $\cup\phantom{[}$\strut\\
  & $[a_3-a_1a_2=a_5-a_1a_4=0,\quad a_4<0,\quad a_1\neq0]$\\
  \hline \hline
%  $6$ & $\Delta_5=\left|
% \begin{matrix}
%    a_1& a_3 & a_5 & 0 &0\\
%   1 & a_2 & a_4 &a_6 &0\\
%   0 & a_1 & a_3 & a_5 &0\\
%   0 & 1 & a_2 & a_4 &a_6\\
%   0 & 0 & a_1 & a_3 & a_5\\
%  \end{matrix}
%  \right|=0$,\quad $-\Delta_3\left|
% \begin{matrix}
%1& a_2 & a_6\\
%   a_1 & a_3 & 0 \\
%   0 & a_1 & a_5 \\
%  \end{matrix}
%  \right|>0$, \quad $a_6\neq0$ \\
%  \hline\hline
\end{tabular}
\end{table}

\begin{remark}\label{HP-2}
The following fact follows from Theorem \ref{Hopf-th3}, which establishes a semi-algebraic description of {\bf HBP}s of a differential system of the form \eqref{L-dot_x}. If there exists an integer $k\in\{1,\ldots,\lfloor \frac{n}{2} \rfloor\}$ such that
the following semi-algebraic system
\begin{equation}\label{Hopf-main-1}
\mathbb{S}_{n,k}:\left\{\begin{array}{l}
f_1(\bar{\boldsymbol{x}},\bar{\boldsymbol{\mu}})=\cdots=f_n(\bar{\boldsymbol{x}},\bar{\boldsymbol{\mu}})=0,\\
{\bf SAD}(\mathcal{C}_k(\bar{\boldsymbol{x}},\bar{\boldsymbol{\mu}}))_{\{=,\geq,>,\neq\}}\\
\end{array}\right.
\end{equation}
has a real root with respect to the variables $\bar{\boldsymbol{x}}=(x_1,\ldots,x_n)$, then the corresponding system \eqref{L-dot_x} has a {\bf HBP}. Here, the notation ${\bf SAD}(\mathcal{C}_k(\bar{\boldsymbol{x}},\bar{\boldsymbol{\mu}}))_{\{=,\geq,>\neq\}}$ is the semi-algebraic description of the characteristic polynomial $p(\lambda)$ of the Jacobian matrix $D_{\boldsymbol{x}}f(\bar{\boldsymbol{x}},\bar{\boldsymbol{\mu}})$ to have a simple pair of pure imaginary roots and no other root with zero real part. The semi-algebraic description of ${\bf SAD}(\mathcal{C}_k(\bar{\boldsymbol{x}},\bar{\boldsymbol{\mu}}))_{\{=,\geq,>,\neq\}}$ can be found in Table  \ref{tab:Hopf-1} for $n=2,\ldots,5$. For $n\geq6$, such a description can be obtained using techniques such as the Sturm--Habicht sequence \cite{Sturm-1989}, or the discriminant sequence \cite{Yang:1999}. In the following, we introduce the theory of triangular decomposition of semi-algebraic systems. In this theory, any semi-algebraic system of the form $\mathbb{S}_{n,k}$ can decompose into finitely many regular semi-algebraic systems such that the union of their zero sets is the zero set of $\mathbb{S}_{n,k}$. As a result, each regular semi-algebraic system describes a subset of the {\bf HBP}s.

\end{remark}
As a summary of this section, we provide a refined version of \cite[Theorem 10]{ChenMM-2011}, which gives an answer to determine when a univariate polynomial has zero as root of multiplicity $k_1$ and $k_2$ pairs of pure imaginary roots while no other roots have zero real parts.

\begin{theorem}\label{Th-10-ChenMM}
Let $f(x)=a_0x^n+a_1x^{n-1}+\cdots+a_n$ be a univariate polynomial of $\mathbb{R}[x]$. Then $f(x)$ has a root $0$ of multiplicity $k_1$ and has $k_2$ pairs of purely imaginary roots while no other roots have zero real parts if and only if the following holds:
\begin{itemize}
    \item[(a)] The coefficients of $f(x)$ satisfy $a_n=\cdots=a_{n-k_1+1}=0,a_{n-k_1}\neq0$.
   \item[(b)] Let $d=n-k_1$, $q(x)=a_0x^{n-k_1}+a_1x^{n-k_1-1}+\cdots+a_{n-k_1}=f_1(x^2)+xf_2(x^2)$. Then there exists an integer $k$ such that the condition $\bigcup\limits_{k\geq k_2}^{\lfloor \frac{d}{2} \rfloor}\mathcal{C}_k$ holds, where $\mathcal{C}_k$ is the condition on the coefficients of $f(x)$ such that $S_k(f_1,f_2,x)$ has $k_2$ negative real roots (counted with multiplicity) and
\[
 \Delta_{d-1-2j}(q)=0\ (0\le j<k),\qquad \Delta_{d-2k-1}(q)\ne0.
\]
\end{itemize}

\end{theorem}

\begin{proof}
Part (a) gives the exact multiplicity at zero. Applying Corollary \ref{lea-HS1} to $q$, the Hurwitz determinants in (b) determine $k=\deg\gcd(f_1,f_2)$. Each negative root $t$ of this GCD with multiplicity $j$ corresponds precisely to a pair of conjugate imaginary roots $\pm i\sqrt{-t}$ of $q$, each bearing multiplicity $j$. Since every nonzero imaginary root of $q$ arises in this manner, summing these multiplicities yields $k_{2}$, thereby establishing both implications.
\end{proof}

\subsubsection{Approach 2: Characterization via polynomial GCD}

The next result, presented below, provides an alternative method of determining the necessary and sufficient conditions for a polynomial $p(\lambda)$ to have a simple pair of pure imaginary roots while no other roots have zero real part.

\begin{theorem}\label{Hopf-gcd1}
Let $p(\lambda)=\lambda^n+a_1(\boldsymbol{x},\boldsymbol{\mu})\lambda^{n-1}+\cdots+a_n(\boldsymbol{x},\boldsymbol{\mu})$ be the characteristic polynomial of the Jacobian matrix $Df(\boldsymbol{x},\boldsymbol{\mu})$. Then $p(\lambda)$ has a simple pair of purely imaginary roots, $i\omega_0$ and $-i\omega_0$, with $\omega_0>0$, and no other root with zero real part, if and only if $a_n\neq0$ and $g(x)={\rm gcd}(u(x),v(x))$ has the two simple nonzero real roots $\omega_0$ and $-\omega_0$ and no other real roots, where $u(x),v(x)\in\mathbb{R}[x]$ are defined by $p(ix)=u(x)+iv(x)$.

\end{theorem}

\begin{proof}
We will first prove necessity and then prove sufficiency.
\begin{itemize}
\item[``$\Rightarrow$''] Suppose that the only roots of $p$ on the imaginary axis are $\pm i\omega_0$, where $\omega_0>0$, and both roots are simple. Since zero is not a root, $a_n=p(0)\ne0$. Moreover,
\[
 p(\pm i\omega_0)=0
 \quad\Longrightarrow\quad
 u(\pm\omega_0)=v(\pm\omega_0)=0.
\]
We have, $x-\omega_0$ and $x+\omega_0$ divide both $u$ and $v$, so both factors divide their greatest common divisor $g$. Thus $\pm\omega_0$ are real roots of $g$.

To show that these roots are simple, differentiate
$p(i x)=u(x)+i v(x)$ to obtain
\begin{equation}\label{eq:derivative}
i p'(i x)=u'(x)+i v'(x).
\end{equation}
If either $x_0=\omega_0$ or $x_0=-\omega_0$ were a multiple root of $g$,
then $(x-x_0)^2$ would divide both $u$ and $v$. Thus $u'(x_0)=v'(x_0)=0$, and \eqref{eq:derivative} would give
$p'(i x_0)=0$, contradicting the simplicity of $\pm i\omega_0$. Hence both real roots of $g$ are simple.

If $s_0\in\R$ is any root of $g$, then $g\mid u$ and $g\mid v$ give
$u(s_0)=v(s_0)=0$, and hence $p(i s_0)=0$. The assumption on the roots of $p$ therefore implies $s_0\in\{-\omega_0,\omega_0\}$. Thus $g$ has no other real roots.

\item[``$\Leftarrow$'']
Suppose that $a_n\ne0$ and the real roots of $g$ are exactly
$\pm\omega_0$, where $\omega_0>0$, and both roots are simple.
Since $g$ divides both $u$ and $v$,
\[
 u(\pm\omega_0)=v(\pm\omega_0)=0,
 \qquad p(\pm i\omega_0)=0.
\]
These two roots are distinct because $\omega_0>0$. Also, $p(0)=a_n\ne0$, so zero is not a root.

To prove simplicity, suppose that $p'(i x_0)=0$ for some
$x_0\in\{-\omega_0,\omega_0\}$. Since $u'(x_0)$ and $v'(x_0)$ are
real, \eqref{eq:derivative} implies $u'(x_0)=v'(x_0)=0$.
Together with $u(x_0)=v(x_0)=0$, this shows that $(x-x_0)^2$ divides
both $u$ and $v$, and hence also $g$, a contradiction. Therefore, $\pm i\omega_0$ are simple roots of $p$.

Now let $i s_0$, with $s_0\in\R\setminus\{0\}$, be any root of $p$
on the imaginary axis. Taking real and imaginary parts of
$p(i s_0)=0$ yields $u(s_0)=v(s_0)=0$. We know that $s-s_0$ divides both $u$ and $v$, and therefore divides $g$. Hence $g(s_0)=0$, which implies $s_0\in\{-\omega_0,\omega_0\}$. Consequently, $p$ has no other root with zero real part.

\end{itemize}

\end{proof}

The following result is a direct consequence of Theorem \ref{Hopf-gcd1}.

\begin{corollary}\label{REM-GCD}
Let ${\bf SAD}_n$ denote the conditions in Theorem \ref{Hopf-gcd1} that $p(\lambda)$ has a simple pair of purely imaginary roots while no other roots have zero real part. Then for $n = 2, 3, 4, 5$, these conditions can be explicitly expressed as follows:
\begin{align*}
{\bf SAD}_2 =& \{a_1=0,\ a_2>0\};\\
{\bf SAD}_3 =& \{a_1 a_2 = a_3,\ a_1a_3 > 0\};\\
{\bf SAD}_4=&\{a_1=a_3 = 0,\ a_4 < 0\} \, \vee \, \{a_1^2a_4-a_1a_2a_3+a_3^2=0,\ a_3a_4(a_1a_2-a_3)\neq0,\\
&\qquad a_1a_4(a_1a_2-a_3)>0\};\\
{\bf SAD}_5=&\{a_{1}^{2} a_{4}^{2}+a_{1} a_{2}^{2} a_{5}-a_{1} a_{2} a_{3} a_{4}-2 a_{1} a_{4} a_{5}\\
&\quad {}-a_{2} a_{3} a_{5}+a_{4} a_{3}^{2}+a_{5}^{2}=0,\\
&\,\, a_{5} a_{1} \left(a_{2} a_{5}-a_{3} a_{4}\right) \left(a_{1}^{2} a_{4}-a_{1} a_{2} a_{3}-a_{1} a_{5}+a_{3}^{2}\right)\neq0,\\
&\,\,-\left(a_{1}^{2} a_{4}-a_{1} a_{2} a_{3}-a_{1} a_{5}+a_{3}^{2}\right) a_{5} \left(a_{1} a_{2}-a_{3}\right)>0\}\, \vee \, \\
&\{a_1=0,\ a_2a_3a_5-a_3^2a_4-a_5^2=0,\,\left(a_{2} a_{5}-a_{3} a_{4}\right) a_{5} a_{3}\neq0,\ a_3a_5>0\}\, \vee \, \\
&\{a_1a_2-a_3=a_1a_4-a_5=a_2a_5-a_3a_4=0,\ a_1a_4a_5\neq0,\ a_1a_5<0\}.\\
\end{align*}

\end{corollary}

\begin{remark}\label{REM-GCD-1}
We remark that practical applications of Theorem \ref{Hopf-gcd1} involve computing the greatest common divisor of parametric univariate polynomials formed by the real and imaginary parts of a univariate polynomial, and determining whether $g(x)$ admits real roots other than a pair of nonzero symmetric roots. For the former problem, we refer to the algorithm of \cite{WangD,wang2wang2023} to compute the GCD of parametric univariate polynomials; for the latter, we may use techniques such as the Sturm--Habicht sequence \cite{Sturm-1989}, or the discriminant sequence \cite{Yang:1999}. Note that the above two procedures usually require detailed analysis of different cases determined by the polynomial coefficients, which is usually inefficient in practical applications. In contrast,  Theorem \ref{Hopf-th3} not only provides simple and explicit conditions, but also describes the structure of the roots of the characteristic polynomial. Therefore, in the subsequent analysis of Hopf bifurcations, we shall mainly rely on the algebraic criterion given by Theorem \ref{Hopf-th3}.
\end{remark}

\section{Semi-Algebraic Criteria for Hopf Bifurcation and Applications}\label{sect-5}

Here we discuss algebraic criterion for detecting Hopf bifurcations of system \eqref{L-dot_x}. Recall that a semi-algebraic description of {\bf HBP}s for a given differential system has been shown in Theorem \ref{Hopf-th3}. We need further to establish an algebraic criterion for detecting the transversality condition (b). Section \ref{sect-6-1} aims to provide a semi-algebraic description for detecting Hopf bifurcations in polynomial differential systems. In Section \ref{sect-6-2}, we propose a heuristic method based on triangular decomposition of semi-algebraic systems for automatic detection of Hopf bifurcations. In Sections \ref{sect-6-3} and \ref{sect-6-4}, we demonstrate the effectiveness of our algorithmic approach using results obtained for the Lorenz differential system, the Rabinovich--Fabrikant model and a four-dimensional business cycle model.

\subsection{Semi-algebraic description of Hopf bifurcation}\label{sect-6-1}

We first introduce some notations before presenting the semi-algebraic description of Hopf bifurcation. Consider a differential system of the form \eqref{L-dot_x}. Assume that there is a smooth curve of equilibria $\boldsymbol{x}_*(\boldsymbol{\mu})$ satisfying $f(\boldsymbol{x}_*(\boldsymbol{\mu}),\boldsymbol{\mu})=0$ and $\boldsymbol{x}_*(\bar{\boldsymbol{\mu}})=\bar{\boldsymbol{x}}$. Define the following equality
\begin{equation}\label{Hopf-semi-trans}
\widehat\Delta_{n-1}(\boldsymbol{\mu})=\Delta_{n-1}(\boldsymbol{x}_*(\boldsymbol{\mu}),\boldsymbol{\mu}).
\end{equation}

Let $\boldsymbol{e}_i\in\mathbb{R}^m$ be the $i$-th parameter coordinate vector, i.e., the vector  $\boldsymbol{e}_i$ has entry one at position $i$ and zero elsewhere. Define the parameter line and its equilibrium curve by
\begin{equation}\label{Hopf-notation-1}
\boldsymbol{\mu}(t)=\bar{\boldsymbol{\mu}}+t\boldsymbol{e}_i=(\bar{\mu}_1,\ldots,\bar{\mu}_{i-1},\bar{\mu}_i+t,\bar{\mu}_{i+1},\ldots,\bar{\mu}_m),\quad \boldsymbol{X}(t)=\boldsymbol{x}_*(\boldsymbol{\mu}(t)).
\end{equation}
Let  $\boldsymbol{\xi}_i=\boldsymbol{X}'(0)=\partial_{\mu_i}\boldsymbol{x}_*(\bar{\boldsymbol{\mu}})$. The curve $(\boldsymbol{X}(t),\boldsymbol{\mu}(t))$ has tangent $(\boldsymbol{\xi}_i,\boldsymbol{e}_i)$ in the state-parameter space. We define the following two polynomials:
\begin{equation}\label{Hopf-notation-2}
p(\lambda,\boldsymbol{\mu})=\det\{\lambda E_n-D_{\boldsymbol{x}}f(\boldsymbol{x}_*(\boldsymbol{\mu}),\boldsymbol{\mu})\},\quad p_t(\lambda)=p(\lambda,\boldsymbol{\mu}(t))=M_t(\lambda^2)+\lambda N_t(\lambda^2).
\end{equation}
The polynomials $M_t,N_t$ collect the coefficients of the even and odd powers of $\lambda$, respectively. Let $J_f=D_{\boldsymbol{x}}f(\bar{\boldsymbol{x}},\bar{\boldsymbol{\mu}})$ be the Jacobian matrix evaluated at the equilibrium point $(\bar{\boldsymbol{x}},\bar{\boldsymbol{\mu}})$, and $\det(J_f)$ denote the determinant of the matrix $J_f$.

The following result gives an algebraic representation of the transversality condition, i.e., condition (b) of Hopf bifurcations.

\begin{lemma}\label{Hopf-lemma-1}
Assume that there is a smooth curve of equilibria $(\boldsymbol{x}_*(\boldsymbol{\mu}),\boldsymbol{\mu})$ with $\boldsymbol{x}_*(\bar{\boldsymbol{\mu}})=\bar{\boldsymbol{x}}$ for system \eqref{L-dot_x}, and the differential system has a Hopf bifurcation point at $(\bar{\boldsymbol{x}},\bar{\boldsymbol{\mu}})$. Let $\lambda(\boldsymbol{\mu})$ be the unique differentiable function of $\boldsymbol{\mu}$ such that $\lambda(\boldsymbol{\mu})$ is a root of the characteristic polynomial $p(\lambda,\boldsymbol{\mu})$ and $\lambda(\bar{\boldsymbol{\mu}})=i\omega_0$. Then we have:
\begin{equation}\label{Hopf-main-2}
\frac{\partial {\rm Re} (\lambda(\boldsymbol{\mu}))}{\partial \mu_i}\Big|_{\mu=\bar{\mu}}\neq0\quad\Longleftrightarrow\quad
\left\{\begin{array}{l}
\frac{\partial\widehat\Delta_{n-1}(\boldsymbol{\mu})}{\partial\mu_i}\Big|_{\boldsymbol{\mu}=\bar{\boldsymbol{\mu}}}\neq0,\quad k=1\\
\gcd(G,H_i)\,\, \text{has no positive real root},\quad k\geq2
\end{array}\right.,
\end{equation}
where
\begin{equation*}
G(z)=\gcd\big(M_0(-z),N_0(-z)\big),\quad
H_i(z)=N_0'(-z)\dot{M}(-z)-M_0'(-z)\dot{N}(-z).
\end{equation*}
The subscript $i$ indexes the parameter direction, dots denote derivatives of the coefficients of $M$ and $N$ at fixed $z$, and primes $M'$, $N'$ denote derivatives with respect to $z$.

\end{lemma}

\begin{proof}
{\bf Case 1: $k=1$}. We first give the proof for case $n\geq4$ and $k=1$. The proof for $n\in\{2,3\}$ can be found in Remark \ref{Hopf-main-rem1}. Let $\lambda(\bar{\boldsymbol{\mu}})=i\omega_0$. According to \cite[Proposition 3.7]{Weber2000}, for $n\geq4$, we have
\begin{equation}\label{Hopf-re1}
\frac{\partial\widehat\Delta_{n-1}(\boldsymbol{\mu})}{\partial\mu_i}\Big|_{\boldsymbol{\mu}=\bar{\boldsymbol{\mu}}}=-2\Delta_{n-3}(\bar{\boldsymbol{x}},\bar{\boldsymbol{\mu}}){\rm res}(\frac{p(\lambda)}{\lambda^2+\omega_0^2},\lambda^2+\omega_0^2)\frac{\partial{\rm Re}(\lambda(\boldsymbol{\mu}))}{\partial \mu_i}\Big|_{\boldsymbol{\mu}=\bar{\boldsymbol{\mu}}}.
\end{equation}
It is obvious that $\frac{\partial\widehat\Delta_{n-1}(\boldsymbol{\mu})}{\partial \mu_i}(\bar{\boldsymbol{x}},\bar{\boldsymbol{\mu}})\neq0$ implies $\frac{\partial{\rm Re}(\lambda)}{\partial \mu_i}(\bar{\boldsymbol{x}},\bar{\boldsymbol{\mu}})\neq0$ and $\Delta_{n-3}(\bar{\boldsymbol{x}},\bar{\boldsymbol{\mu}})\neq0$. Simplicity of ($i\omega_0$) gives the resultant nonzero. By Theorem \ref{Hopf-th3}(a), $\Delta_{n-3}\neq0$ when $k=1$, and $\Delta_{n-3}=0$ when $k\geq2$. This proves the both directions for $k=1$.

{\bf Case 2: $k\geq2$}. Write $\lambda(t)=\lambda(\boldsymbol{\mu}(t))$, the definition of a partial derivative gives
\[\frac{\partial {\rm Re}\lambda(\bar{\boldsymbol{\mu}})}{\partial \mu_i}=\lim_{h\to0}\frac{{\rm Re}\lambda(\bar{\boldsymbol{\mu}}+h \boldsymbol{e}_i)-{\rm Re}\lambda(\bar{\boldsymbol{\mu}})}{h}=\frac{d}{dt}{\rm Re}\lambda(t)\Big|_{t=0}.\]
In the difference quotient, $h$ is exactly the change in $\mu_i$. For $z>0$, it follows from \eqref{Hopf-notation-2} that
\[p_0(\lambda)\equiv M_0(-z)+\lambda N_0(-z) \pmod{\lambda^2+z}.\]
It vanishes exactly when $p_0$ has the root pair $\pm i\sqrt{z}$. From the hypothesis, we see that the polynomial $G$ has a unique positive root $z=\omega_0^2$ and this root is simple. Note that $p_t(\lambda)$ can be written as
\begin{equation}\label{Hopf-conb-lem-eq1}
p_t(\lambda)=h_t(\lambda)q_t(\lambda),\quad h_t(\lambda)=\lambda^2+\eta(t)\lambda+\zeta(t),
\end{equation}
where $\eta(t)=-2{\rm Re}\lambda(t)$ and $\zeta(t)=|\lambda(t)|^2$. At $t=0$, we have $h_0=\lambda^2+z=\lambda^2+\omega_0^2$, $\eta(0)=0$ and $\zeta(0)=\omega_0^2$. Since $M_0(-z)=N_0(-z)=0$, we can write $M_0(y)=(y+z)m_0(y)$, $N_0(y)=(y+z)n_0(y)$. Then
\[q_0(\lambda)=m_0(\lambda^2)+\lambda n_0(\lambda^2)\equiv u+\lambda v\pmod{\lambda^2+z},\]
where
\[u=m_0(-z)=M_0'(-z),\quad v=n_0(-z)=N_0'(-z).\]
Note that $i\omega_0$ is the simple root of $p_0$, so $q_0(i\omega_0)\neq0$, which implies that $d=u^2+zv^2=|q_0(i\omega_0)|^2>0$.

Differentiate $p_t=h_tq_t$ at $t=0$, keeping $\lambda$ fixed, and reduce the result modulo the polynomial $h_0=\lambda^2+z$, we have
\[\dot{M}(-z)+\lambda\dot{N}(-z)\equiv(\dot\zeta+\lambda\dot\eta)(u+\lambda v)\pmod{\lambda^2+z}.\]
Matching the constant and linear coefficients yields $\dot{M}(-z)=u\dot{\zeta}-zv\dot{\eta}$ and $\dot{N}(-z)=v\dot{\zeta}+u\dot{\eta}$. By multiply the first identity by $v$ and subtract $u$ times the second, we cancel the term $\dot{\zeta}$. Note that $\dot{\eta}=-2\partial_{\mu_i}{\rm Re}\lambda(\bar{\boldsymbol{\mu}})$. Substitution gives
\[H_i(z)=v\dot{M}(-z)-u\dot{N}(-z)=-d\cdot\dot{\eta}=2d\cdot\frac{\partial {\rm Re}\lambda(\bar{\boldsymbol{\mu}})}{\partial \mu_i}.\]

\begin{itemize}
\item[``$\Rightarrow$''] Since $i\omega_0$ is simple $d>0$. The identity above gives $H_i(z)\ne0$. Suppose that $\gcd(G,H_i)$ had a positive real root $z_0$. As $\gcd(G,H_i)$ divides $G$, we would have $G(z_0)=0$. On the other hand, we know that $z=\omega_0^2$ is the unique positive root of $G$, so $z_0=z$. But $\gcd(G,H_i)$ also divides $H_i$, which would give $H_i(z)=0$, a contradiction. Therefore $\gcd(G,H_i)$ has no positive real root.

\item[``$\Leftarrow$''] Assume $\gcd(G,H_i)$ has no positive real root. Note that $G(z)$ has a unique positive root $z=\omega_0^2>0$. If $H_i(z)=0$, then we have \(z-\omega_0^2\) a common divisor of $G$ and $H_i$ in $\mathbb R[z]$. It would then divide $\gcd(G,H_i)$, contrary to the assumption. Thus $H_i(z)\ne0$, which implies $\partial_{\mu_i}{\rm Re}\lambda(\bar{\mu})\neq0$.

\end{itemize}

\end{proof}

\begin{remark}\label{Hopf-main-rem1}
When $k=1$, the statement of the first equation in \eqref{Hopf-main-2} holds for any $n\geq2$. In fact, let $\lambda(\boldsymbol{\mu})$ be the unique differentiable function of $\boldsymbol{\mu}$ such that $\lambda(\boldsymbol{\mu})$ is a root of the characteristic polynomial $p(\lambda)$ and $\lambda(\bar{\boldsymbol{\mu}})=i\omega_0$. Let us factor the polynomial $p(\lambda)$ into the form
\begin{equation}\label{Hopf-main-eq1}
p(\lambda)=(\lambda^2-2{\rm Re}(\lambda(\boldsymbol{\mu}))\lambda+\lambda(\boldsymbol{\mu})\bar{\lambda}(\boldsymbol{\mu}))p^*(\lambda).
\end{equation}
When $n=2$, we have $p^*(\lambda)=1$. By comparing the coefficients of the right-hand side of \eqref{Hopf-main-eq1} with $\lambda^2+a_1\lambda+a_2$, we have $a_1=\widehat\Delta_{1}(\boldsymbol{\mu})=-2{\rm Re}(\lambda(\boldsymbol{\mu}))$. So, $\frac{\partial\widehat\Delta_{1}(\boldsymbol{\mu})}{\partial \mu_i}\Big|_{\boldsymbol{\bar{\boldsymbol{\mu}}}}=-2\frac{\partial{\rm Re}(\lambda(\boldsymbol{\mu}))}{\partial \mu_i}\Big|_{\boldsymbol{\bar{\boldsymbol{\mu}}}}$. When $n=3$, we have $p^*(\lambda)=\lambda-\lambda_3$ with $\lambda_3\neq0$. Then, using Orlando’s formula (see \cite{Gantmakher:1959}), one obtains
\begin{equation}\label{Hopf-main-eq2}
\begin{split}
\widehat\Delta_{2}(\boldsymbol{\mu})&=-(\lambda(\boldsymbol{\mu})+\bar{\lambda}(\boldsymbol{\mu}))\cdot (\lambda(\boldsymbol{\mu})+\lambda_3)\cdot(\bar{\lambda}(\boldsymbol{\mu})+\lambda_3),\\
&=-2{\rm Re}(\lambda(\boldsymbol{\mu}))\cdot (\lambda(\boldsymbol{\mu})+\lambda_3)\cdot(\bar{\lambda}(\boldsymbol{\mu})+\lambda_3).
\end{split}
\end{equation}
Differentiating both sides of \eqref{Hopf-main-eq2} with respect to $\mu_i$ and noting that ${\rm Re}(\lambda(\bar{\boldsymbol{\mu}}))=0$, one obtains
\begin{equation*}
\frac{\partial\widehat\Delta_{2}(\boldsymbol{\mu})}{\partial \mu_i}\Big|_{\boldsymbol{\bar{\boldsymbol{\mu}}}}=-2(\lambda_3^2+\omega_0^2)\frac{\partial{\rm Re}(\lambda(\boldsymbol{\mu}))}{\partial \mu_i}\Big|_{\boldsymbol{\bar{\boldsymbol{\mu}}}}.
\end{equation*}
This gives the desired result.
\end{remark}

Let ${\rm adj}(J_f)$ denote the adjugate matrix of the Jacobian matrix $J_f=D_{\boldsymbol{x}}f(\bar{\boldsymbol{x}},\bar{\boldsymbol{\mu}})$. The following result provides an equivalent polynomial characterization of \eqref{Hopf-main-2} when $k=1$.

\begin{lemma}\label{lem-transversality}
We have the following equivalence relation
\begin{equation}\label{equ-trans-1}
\frac{\partial\widehat\Delta_{n-1}(\boldsymbol{\mu})}{\partial\mu_i}\Big|_{\boldsymbol{\mu}=\bar{\boldsymbol{\mu}}}\neq0\quad \Longleftrightarrow\quad  (\det J_f)\,\partial_{\mu_i}\Delta_{n-1}
 -(\nabla_{\boldsymbol{x}}\Delta_{n-1})^{\mathsf T}{\rm adj}(J_f) f_{\mu_i}\neq0.
\end{equation}
\end{lemma}

\begin{proof}
Taking partial derivative of the equilibrium equation $f(\boldsymbol{x}_*(\boldsymbol{\mu}),\boldsymbol{\mu})=0$ with respect to $\mu_i$, we have
\begin{equation}\label{equ-trans-2}
J_f(\partial_{\mu_i}\boldsymbol{x}_*)(\bar{\boldsymbol{\mu}})=-f_{\mu_i}.
\end{equation}
Differentiating both sides of \eqref{Hopf-semi-trans} with respect to $\mu_i$ and using the relation \eqref{equ-trans-2}, we have
\begin{equation}\label{equ-trans-3}
\frac{\partial\widehat\Delta_{n-1}}{\partial\mu_i}\Big|_{\bar{\boldsymbol{\mu}}}
 =\left[\partial_{\mu_i}\Delta_{n-1}
       -(\nabla_{\boldsymbol{x}}\Delta_{n-1})^{\mathsf T}(J_f)^{-1}f_{\mu_i}\right].
\end{equation}
The right-hand partial derivatives both in \eqref{equ-trans-2} and \eqref{equ-trans-3} act on the original polynomials in independent variables $(\boldsymbol{x},\boldsymbol{\mu})$ and are evaluated at $(\bar{\boldsymbol{x}},\bar{\boldsymbol{\mu}})$. Note that $\det J_f=\det D_{\boldsymbol{x}}f(\bar{\boldsymbol{x}},\bar{\boldsymbol{\mu}})\neq0$. Multiplying both sides of equation \eqref{equ-trans-3} by $\det J_f$, we obtain
\begin{equation}\label{equ-trans-4}
(\det J_f)\frac{\partial\widehat\Delta_{n-1}}{\partial\mu_i}\Big|_{\bar{\boldsymbol{\mu}}}
 =\left[(\det J_f)\partial_{\mu_i}\Delta_{n-1}
       -(\nabla_{\boldsymbol{x}}\Delta_{n-1})^{\mathsf T}{\rm adj}(J_f)f_{\mu_i}\right].
\end{equation}
The above equality holds in both directions. This completes the proof.

\end{proof}

We are now in a position to derive semi-algebraic descriptions of Hopf bifurcation. In particular, the following result, based on Lemmas \ref{Hopf-lemma-1} and \ref{lem-transversality},  characterizes \textit{complete} Hopf bifurcation (\textbf{CHB}) in a semi-algebraic form.

\begin{proposition}\label{Hopf-prop1}
Assume that there is a smooth curve of equilibria $(\boldsymbol{x}_*(\boldsymbol{\mu}),\boldsymbol{\mu})$ with $\boldsymbol{x}_*(\bar{\boldsymbol{\mu}})=\bar{\boldsymbol{x}}$ for system \eqref{L-dot_x}. If at least one of the semi-algebraic systems in the disjunction
\[\bigvee\limits_{k=1}^{\lfloor \frac{n}{2} \rfloor}\mathbb{H}_{n,k}(\bar{\boldsymbol{x}},\bar{\boldsymbol{\mu}})\]
admits a real solution with respect to $\bar{\boldsymbol{x}}=(x_1,\ldots,x_n)$, then the differential system exhibits a {\bf CHB}. The expressions of $\mathbb{H}_{n,k}$ are the following:
\begin{equation}\label{Hopf-main-2-1}
\mathbb{H}_{n,1}:\left\{\begin{array}{l}
\bigwedge\limits_{j=1}^n
\Big(f_j(\bar{\boldsymbol{x}},\bar{\boldsymbol{\mu}})=0\Big)\\
\bigwedge\Big({\bf SAD}(\mathcal{C}_1(\bar{\boldsymbol{x}},\bar{\boldsymbol{\mu}}))_{\{=,\geq,>,\neq\}}\Big)\\
\bigwedge\Big((\det J_f)\,\partial_{\mu_i}\Delta_{n-1}
 -(\nabla_{\boldsymbol{x}}\Delta_{n-1})^{\mathsf T}{\rm adj}(J_f) f_{\mu_i}\neq0\Big)
\end{array}\right.,
\end{equation}
and
\begin{equation}\label{Hopf-main-2-2}
\quad \mathbb{H}_{n,k}:\left\{\begin{array}{l}
\bigwedge\limits_{j=1}^n
\Big(f_j(\bar{\boldsymbol{x}},\bar{\boldsymbol{\mu}})=0\Big)\\
\bigwedge\Big({\bf SAD}(\mathcal{C}_k(\bar{\boldsymbol{x}},\bar{\boldsymbol{\mu}}))_{\{=,\geq,>,\neq\}}\Big)\\
\bigwedge\Big(T_i(\bar{\boldsymbol{x}},\bar{\boldsymbol{\mu}})\Big)
\end{array}\right.\quad k\geq2,
\end{equation}
where $\mu_i$ is the chosen parameter, $T_i$ is the semi-algebraic description of the condition that $\gcd(G,\allowbreak H_i)$ has no positive real root.

\end{proposition}

\begin{proof}
The desired result is a direct consequence of Theorem \ref{Hopf-th3}, Lemmas \ref{Hopf-lemma-1} and \ref{lem-transversality}.
\end{proof}

\begin{remark}
For a given differential system of the form \eqref{L-dot_x}, the expression of $T_i$ may be derived by computing $\gcd$ of two polynomials. Below, we present the detailed derivation steps for the four-dimensional case. Write
\begin{equation}\label{4d-example-eq1}
\begin{split}
p_t(\lambda)&=\det\{\lambda E_4-D_{\boldsymbol{x}}f(\boldsymbol{X}(t),\boldsymbol{\mu}(t))\}\\
&=\lambda^4+a_1(t)\lambda^3+a_2(t)\lambda^2+a_3(t)\lambda+a_4(t)\\
&=M_t(\lambda^2)+\lambda N_t(\lambda^2),
\end{split}
\end{equation}
where $\boldsymbol{X}(t)$ and $\boldsymbol{\mu}(t)$ are defined in \eqref{Hopf-notation-1}. Hence, $M_0(y)=y^2+a_2(0)y+a_4(0)$ and $N_0(y)=a_1(0)y+a_3(0)$. Note that, at $t=0$, we have $a_j(0)=a_j(\bar{\boldsymbol{x}},\bar{\boldsymbol{\mu}})$. When $k=2$, by Theorem \ref{Hopf-th3} we have $a_1(0)=a_3(0)=0$ and $a_4(0)<0$. In this case, the polynomial $G(z)=\gcd(M_0(-z),N_0(-z))=z^2-a_2(0)z+a_4(0)$. The polynomial $H_i(z)$ simplifies to
\[H_i(z)=(2z-a_2(0))(\dot{a}_3-\dot{a}_1z).\]
A dot on $a_j$ denotes the derivative of $a_j(\boldsymbol{X}(t),\boldsymbol{\mu}(t))$ at $t=0$. The factor $2z-a_2(0)$ is nonzero at the roots of $G$ since $a_2^2(0)-4a_4>0$. So $T_i$ is equivalent to the condition that $\gcd(G,\dot{a}_3-\dot{a}_1z)$ has no positive real root. Assume that $\dot{a}_3=A$, $\dot{a}_1=B$. We obtain the semi-algebraic $T_i$:
\begin{equation}\label{4d-example-eq2}
\begin{split}
T_i(\bar{\boldsymbol{x}},\bar{\boldsymbol{\mu}})=\Big(-AB>0\Big)\,\,\vee\,\,\Big(A^2-a_2(0)AB+a_4(0)B^2\neq0\Big).
\end{split}
\end{equation}
The above elimination can be derived using a procedure similar to that in Corollary \ref{REM-GCD}. The expressions for $A$ and $B$ are evaluated at $a_j(0)=a_j(\bar{\boldsymbol{x}},\bar{\boldsymbol{\mu}})$. Indeed, by differentiating the associated equilibrium equation $f(\boldsymbol{X}(t),\boldsymbol{\mu}(t))=0$ and a coefficient $a_j(\boldsymbol{X}(t),\boldsymbol{\mu}(t))$ at $t=0$, applying the chain rule we have
\begin{equation}\label{4d-example-eq3}
\begin{split}
J_f\boldsymbol{\xi}_i+\partial_{\mu_i}f=0,\quad \dot{a}_j=\partial_{\mu_i}a_j(\bar{\boldsymbol{x}},\bar{\boldsymbol{\mu}})+\left(\nabla_{\boldsymbol{x}}a_j(\bar{\boldsymbol{x}},\bar{\boldsymbol{\mu}})\right)^T\boldsymbol{\xi}_i.
\end{split}
\end{equation}
Combining the two equations above and performing an operation similar to that in \eqref{equ-trans-4}, we have
\begin{equation}\label{4d-example-eq4}
\begin{split}
\det(J_f)\dot{a}_j=\det(J_f)\partial_{\mu_i}a_j(\bar{\boldsymbol{x}},\bar{\boldsymbol{\mu}})-\left(\nabla_{\boldsymbol{x}}a_j(\bar{\boldsymbol{x}},\bar{\boldsymbol{\mu}})\right)^T{\rm adj}(J_f)\partial_{\mu_i}f.
\end{split}
\end{equation}

In practice, we multiply $H_i$ by $\det(J_f)$ to clear denominators and facilitate symbolic polynomial manipulation.

\end{remark}

%$J_f(\boldsymbol{x},\boldsymbol{\mu})=D_{\boldsymbol{x}}f(\boldsymbol{x},\boldsymbol{\mu})$ be the Jacobian matrix of system \eqref{L-dot_x}, $A(t)=J_f(\boldsymbol{X}(t),\boldsymbol{\mu}(t))$.
%

%For $n\in\{2,3\}$ and $k=1$, or for $n\geq4$ and an integer $k\in\{1,\ldots,\lfloor \frac{n}{2} \rfloor-1\}$,

In many applications, researchers are typically concerned with a specific type of Hopf bifurcation where a pair of complex conjugate eigenvalues of the Jacobian matrix crosses the imaginary axis while all other eigenvalues have negative real parts. This type of bifurcation, termed a simple Hopf bifurcation ({\bf SHB}), was introduced by Liu \cite{JMAA:1994}. The following criterion can be used to detect {\bf SHB} for differential systems of the form \eqref{L-dot_x}.

\begin{proposition}\label{Hopf-S}
Assume that there is a smooth curve of equilibria $(\boldsymbol{x}_*(\boldsymbol{\mu}),\boldsymbol{\mu})$ with $\boldsymbol{x}_*(\bar{\boldsymbol{\mu}})=\bar{\boldsymbol{x}}$ for system \eqref{L-dot_x}. If the following semi-algebraic system
\begin{equation}\label{Hopf-main-L}
\mathbb{H}_{n}:\left\{\begin{array}{l}
\bigwedge\limits_{j=1}^n
\Big(f_j(\bar{\boldsymbol{x}},\bar{\boldsymbol{\mu}})=0\Big)\ \bigwedge \
\Big(\Delta_{n-1}(\bar{\boldsymbol{x}},\bar{\boldsymbol{\mu}})=0\Big)\\
\bigwedge \Big(a_n(\bar{\boldsymbol{x}},\bar{\boldsymbol{\mu}})>0\Big)\ \bigwedge\limits_{\ell=1}^{n-2} \Big(\Delta_{\ell}(\bar{\boldsymbol{x}},\bar{\boldsymbol{\mu}})>0\Big)\\
\bigwedge\Big((\det J_f)\,\partial_{\mu_i}\Delta_{n-1}
 -(\nabla_{\boldsymbol{x}}\Delta_{n-1})^{\mathsf T}{\rm adj}(J_f) f_{\mu_i}\neq0\Big)
\end{array}\right.
\end{equation}
has real roots with respect to $\bar{\boldsymbol{x}}=(x_1,\ldots,x_n)$. The differential system then undergoes a {\bf SHB}.

\end{proposition}

\begin{proof}
By Theorem \ref{Hopf-th3}, we have $k=1$. So the equality \eqref{Hopf-re1} holds. Thus the desired result follows from Lemma \ref{lem-transversality}. In fact, this result is a consequence of \cite[Theorem 3.6 and Proposition 3.7]{Weber2000}, as well as a slightly modified version of the main result given in \cite[Section 3]{JMAA:1994}.
\end{proof}

\begin{remark}
We will use the method of triangular decomposition to solve the obtained semi-algebraic systems. An algorithm for computing triangular decomposition of a semi-algebraic system is called
\textsf{RealTriangularize}, which has been implemented as part of the \textsf{RegularChains} library in \textsc{Maple}. In practice, we will take $\mu_i$ a suitable parameter in the variable ordering when running the \textsf{RealTriangularize} command. In the experimental part of this paper, we will mark the required parameter in red to illustrate the choice we have made.
\end{remark}

\subsection{Heuristic real quantifier elimination via {\sf RealTriangularize}}\label{sect-6-2}

In this section, we establish when a real triangular decomposition of a semi-algebraic system suffices to yield the required bifurcation conditions. Specifically, Algorithm \ref{algorithm:rtd} implements a heuristic real quantifier elimination procedure via {\sf RealTriangularize}.

\begin{algorithm}
\renewcommand{\algorithmicrequire}{\textbf{Input:}}
\renewcommand{\algorithmicensure}{\textbf{Output:}}
\renewcommand{\algorithmiccomment}[1]{\hfill // #1}
\caption{Extract Parameter Conditions from Regular Semi-algebraic Systems.}
\label{algorithm:rtd}
\begin{algorithmic}[1]
\REQUIRE
{A regular semi-algebraic system $rsa:=[{\cal Q}, {\cal T}, {\cal P}]$
of $\Q[\bmu,\x]$; free variables $\bfu$ and main variables $\bfy$ of ${\cal T}$.}
\ENSURE {The projection of $Z_{\R}(rsa)$ onto the parameter space.}
\IF{no variables in $\x$ appear in $rsa$}
\STATE return $rsa$ \COMMENT{this includes the case that $\x=\emptyset$}
\ELSIF{$\x=\bfy$}
\STATE return ${\cal Q}$
\ELSIF{$\bfy\subset\x$}
\STATE return the output of the QE problem: $\exists (\x\setminus\bfy), {\cal Q}$ \COMMENT{if $\bfy=\emptyset$, then ${\cal Q}\Leftrightarrow rsa$}
\ELSIF{$\x\subset \bfy$}
\STATE let $\bfz=\bfy\setminus\x$\;
\STATE let $rsa' := rsa\cap \Q[\bfu,\bfz]=[{\cal Q}', {\cal T}', {\cal P}']$ and $rsa'' := [{\cal Q}'', {\cal T}'', {\cal P}'']= [true, {\cal T}\setminus {\cal T}', {\cal P}\setminus {\cal P}']$\;
\IF{$rsa'=[{\cal Q}', {\cal T}', {\cal P}']$ has only one real solution in $\bfz$ for any $u\in Z_{\R}({\cal Q})$}
\STATE \COMMENT{This condition holds if the main degrees of all polynomials
in ${\cal T}'$ are all one.}\;
\STATE return $rsa'$
\ELSIF{$rsa''$ has real solutions in $\x$ for any real zero of $rsa'$}
\STATE \COMMENT{It holds if the main degrees of all polynomials
in ${\cal T}''$ are all one and ${\cal P}''=\emptyset$.}\;
%, which is satisfied in high probability due to shape lemma
\STATE return $rsa'$
\ELSE
 \STATE return the output of the QE problem: $\exists \x, rsa$\;
\ENDIF
\ELSIF{$\x\cap\bfy\neq\emptyset$}
\STATE let $\x' := \x\setminus \bfy$; $\x'' := \x\cap\bfy$ and $\boldsymbol{z} := \boldsymbol{y}\setminus \x$ \COMMENT{we must have $\x''>\bfz$ as $\x>\bmu$}
\STATE let $rsa' := rsa\cap \Q[\bfu,\bfz]=[{\cal Q}', {\cal T}', {\cal P}']$ and $rsa'' := [{\cal Q}'', {\cal T}'', {\cal P}'']= [true, {\cal T}\setminus {\cal T}', {\cal P}\setminus {\cal P}']$\;
\IF{$rsa'=[{\cal Q}', {\cal T}', {\cal P}']$ has only one real solution in $\bfz$ for any $u\in Z_{\R}({\cal Q})$}
\STATE return $\exists \x', rsa'$\;
\ELSIF{$rsa''=[true, {\cal T}'', {\cal P}'']$ has real solutions in $\x''$ for any real zero of $rsa'$}
\STATE return $\exists \x', rsa'$\;
\ELSE
\STATE return the output of the QE problem: $\exists \x, rsa$\;
\ENDIF
\ELSE
\STATE return the output of the QE problem: $\exists \x, rsa$\;
\ENDIF
% \FOR{$i = 1$ \TO $s$}
% \STATE \;
% \ENDFOR
\end{algorithmic}
\end{algorithm}

\begin{proposition}
Algorithm \ref{algorithm:rtd} is correct.
\end{proposition}
\begin{proof}
If no variables in $\x$ appear in $rsa$, the algorithm obviously returns the correct output at line $2$.

By the definition of regular semi-algebraic system, we have
${\cal Q}\Leftrightarrow\exists \bfy, rsa$.
To get the output, we need to solve $\exists \x, rsa$.
Thus, if $\x=\bfy$, the algorithm returns correct output at line $4$.
If $\bfy\subset \x$, then $\exists \x, rsa$ is equivalent to $\exists (\x\setminus\bfy, \bfy), rsa$, which implies that the algorithm correctly returns at line $6$.

Note that starting from line $7$, we must have $\x\neq\emptyset$
and $\bfy\neq\emptyset$.
If $\x\subset\bfy$, then $\bmu=(\bfu, \bfz)$ and ${\cal Q}$ only involves variables in $\bmu$,
which implies that ${\cal Q}={\cal Q}'$ and ${\cal Q}''=\emptyset$.
Clearly we have $\pi_{\bmu}(Z_{\R}(rsa))\subseteq Z_{\R}(rsa')$.

Given any point $u$ of $Z_{\R}({\cal Q})$, if the condition at line $10$ holds,
there exists a unique $z$ such that $(u,z)\in Z_{\R}(rsa')$.
Thus the projection $\pi: Z_{\R}(rsa')\rightarrow Z_{\R}({\cal Q})$ is a bijection.
On the other hand, as $rsa$ is a regular semi-algebraic system,
there exists $x$ such that $rsa(u,z,x)$ is satisfied.
Thus the algorithm returns the correct output at line $12$.
If the condition at line $11$ holds, the condition at line $10$ trivially holds.

Given any point $\mu=(u,z)$ of $Z_{\R}(rsa')$, if the condition at line $13$
is satisfied, then there exists $x$ such that $rsa(\mu,x)$ is satisfied,
and thus the algorithm returns the correct output at line $15$.
Next we show that the condition at line $14$ holds implies that
the condition at line $13$ holds.
As $rsa$ is a regular semi-algebraic system, the regular system $[{\cal T}, {\cal P}]$
specializes well at any point $u$ of $Z_{\R}({\cal Q})$
and the specialized regular chain ${\cal T}(u)$ is zero-dimensional.
This implies that neither the initials of polynomials in ${\cal T}''$ nor the polynomials in ${\cal P}''$ vanish at any point $(u, z)$ of $Z_{\R}(rsa')$.
As the main degrees of all polynomials in ${\cal T}''$ are linear and ${\cal P}''=\emptyset$, there must exist $x$ such that $rsa''(u,z,x)$ has real solutions.

At line $18$, we assume that $\x\cap\bfy\neq\emptyset$.
Now $\x''\subseteq\x$ and $\bfz\subset\bmu$ are both main variables of ${\cal T}$,
which implies that $\x''>\bfz$ as $\x>\bmu$.
Therefore both $rsa'$ and $rsa''$ are well defined.
Following similar arguments as before, the algorithm returns
the correct output at line $22$ and line $24$.
For line $17$, $26$ and $28$, the structure of regular semi-algebraic system itself does not provide much useful property
to be exploited, and a general quantifier procedure has to be applied to get the projection.

\end{proof}

\subsection{Illustrative example}\label{sect-6-3}
We first illustrate the complete workflow of our algorithmic tests on the Lorenz system, followed by a targeted comparison with existing methods in the literature to validate our approach's performance in detecting Hopf bifurcations. The famous Lorenz system is given by (see \cite{Lorenz1963})
\begin{equation}\label{Lorenz-1}
\left\{\begin{array}{l}
\dot{x}=\sigma(y-x),\\
\dot{y}=\rho x-y-xz,\\
\dot{z}=-bz+xy,
\end{array}\right.
\end{equation}
where $\boldsymbol{x}=(x,y,z)$ are variables, and $\boldsymbol{\mu}=(\rho,\sigma,b)$ are real parameters. The characteristic polynomial $p(\lambda)$ of the Jacobian matrix $Df(\boldsymbol{x},\boldsymbol{\mu})$ is
\begin{equation}\label{Lorenz-Cp}
p(\lambda)=\lambda^3+a_1\lambda^2+a_2\lambda+a_3,
\end{equation}
where
\begin{equation*}
\begin{split}
a_1&=b +1+\sigma,\\
a_2&=b \sigma -\rho  \sigma +\sigma  z +x^{2}+b +\sigma,\\
a_3&=-b \rho  \sigma +b \sigma  z +x^{2} \sigma +y \sigma  x +b \sigma.
\end{split}
\end{equation*}

Using Theorem \ref{Hopf-th3}, Propositions \ref{Hopf-prop1} and \ref{Hopf-S}, we formulate the following
semi-algebraic set $\mathbb{S}_{3,1}$, $\mathbb{H}_{3,1}$ and $\mathbb{H}_{3}$:
\begin{equation}\label{Hopf-main-L1}
\mathbb{S}_{3,1}:\left\{\begin{array}{l}
\sigma(y-x)=0,\quad \rho x-y-xz=0,\quad -bz+xy=0,\\
\Delta_2=\left|
 \begin{matrix}
    a_1& a_3 \\
   1 & a_2 \\
  \end{matrix}
  \right|=0,\quad a_1a_3>0,\\
\end{array}\right.
\end{equation}

\begin{equation}\label{Hopf-main-L2}
\mathbb{H}_{3,1}:\left\{\begin{array}{l}
\sigma(y-x)=0,\quad \rho x-y-xz=0,\quad -bz+xy=0,\\
\Delta_2=\left|
 \begin{matrix}
    a_1& a_3 \\
   1 & a_2 \\
  \end{matrix}
  \right|=0,\quad a_1a_3>0,\\
\frac{\partial\widehat\Delta_2}{\partial b}\neq0
\end{array}\right.
\end{equation}
and
\begin{equation}\label{Hopf-main-L3}
\mathbb{H}_{3}:\left\{\begin{array}{l}
\sigma(y-x)=0,\quad \rho x-y-xz=0,\quad -bz+xy=0,\\
\Delta_2=\left|
 \begin{matrix}
    a_1& a_3 \\
   1 & a_2 \\
  \end{matrix}
  \right|=0,\quad a_3>0,\quad a_1>0,\\
\frac{\partial\widehat\Delta_2}{\partial b}\neq0
\end{array}\right.
\end{equation}
for detecting {\bf HBP}, {\bf CHB}, and {\bf SHB} of the Lorenz system \eqref{Lorenz-1}. In practice, the last equation in \eqref{Hopf-main-L2} or \eqref{Hopf-main-L3} will be replaced by using Lemma \ref{lem-transversality}. By running the {\sf RealTriangularize} command on the semi-algebraic set $\mathbb{S}_{3,1}$ for the variable ordering $x>y>z>\rho>\sigma>{\color{red}b}$, we obtain the three following regular semi-algebraic systems for {\bf HBP} (unspecified $R_i^P$ are empty):
\begin{equation}\label{Lorenz-R1}
\begin{split}
\mathcal{C}_{{\rm HBP}}=\big(R_1^T\,\,\wedge \,\,R_1^{\mathcal{Q}}\big)\,\,\vee\,\,\big(R_2^T\,\,\wedge \,\,R_2^{\mathcal{Q}}\big)\,\,\vee\,\,
\big(R_{3}^T\,\,\wedge \,\,R_{3}^{\mathcal{Q}}\big),
\end{split}
\end{equation}
where
\begin{equation*}
\begin{split}
R_1^T&=\left\{\begin{array}{l}
x=0\\
y=0\\
z=0\\
\sigma+1=0
\end{array}\right.,\quad
R_1^{\mathcal{Q}}=\left\{\begin{array}{l}
\rho-1>0\\
b\neq0
\end{array}\right.,\quad
R_2^T=\left\{\begin{array}{l}
x-y=0\\
y^2+2\rho-2=0\\
z-\rho+1=0\\
b+2=0\\
\sigma+1=0
\end{array}\right.,\\
R_2^{\mathcal{Q}}&=\left\{\begin{array}{l}
1-\rho>0
\end{array}\right.,\quad
R_3^T=\left\{\begin{array}{l}
x-y=0\\
y^2+(-\rho+1)b=0\\
z-\rho+1=0\\
\left(-b +\sigma -1\right) \rho -\sigma^{2}+\left(-3-b \right) \sigma=0
\end{array}\right.,\\
R_3^{\mathcal{Q}}&=
\left\{\begin{array}{l}
\big(\sigma+b+1>0 \wedge b-\sigma+1<0 \wedge b>0 \wedge \sigma-1>0\big) \\
\vee (\sigma+b+1>0 \wedge b-\sigma+1>0 \wedge b<0 \wedge \sigma-1<0 \wedge \sigma>0\big) \\
\vee (\sigma+b+1<0 \wedge b-\sigma+1<0 \wedge b<0 \wedge \sigma<0 \wedge \sigma+1>0\big)\\
\vee (\sigma+b+1<0 \wedge b-\sigma+1>0 \wedge b<0 \wedge \sigma+1<0\big)
\end{array}\right..
\end{split}
\end{equation*}

%\begin{equation}\label{Lorenz-R}
%\begin{split}
%R_1^T&=\left\{\begin{array}{l}
%x\\
%y\\
%z\\
%\sigma+1
%\end{array}\right.,\quad
%R_1^{\mathcal{Q}}=\left\{\begin{array}{l}
%b\neq0
%\end{array}\right.,
%\quad R_1^P=\left\{\begin{array}{l}
%\rho-1>0
%\end{array}\right.,\\
%R_2^T&=\left\{\begin{array}{l}
%x-y\\
%y^2+2\rho-2\\
%z-\rho+1\\
%b+2\\
%\sigma+1
%\end{array}\right.,\quad R_2^P=\left\{\begin{array}{l}
%1-\rho>0
%\end{array}\right.,\quad
%R_3^T=\left\{\begin{array}{l}
%x-y\\
%y^2+(-\rho+1)b\\
%z-\rho+1\\
%\left(-b +\sigma -1\right) \rho -\sigma^{2}+\left(-3-b \right) \sigma
%\end{array}\right.,\\
%R_3^{\mathcal{Q}}&=\left\{\begin{array}{l}
%\big(\sigma+b+1>0 \wedge b-\sigma+1<0 \wedge b>0 \wedge \sigma-1>0\big) \\
%\vee (\sigma+b+1>0 \wedge b-\sigma+1>0 \wedge b<0 \wedge \sigma-1<0 \wedge \sigma>0\big) \\
%\vee (\sigma+b+1<0 \wedge b-\sigma+1<0 \wedge b<0 \wedge \sigma<0 \wedge \sigma+1>0\big)\\
%\vee (\sigma+b+1<0 \wedge b-\sigma+1>0 \wedge b<0 \wedge \sigma+1<0\big)
%\end{array}\right..\nonumber
%\end{split}
%\end{equation}

It follows directly from $R_1$ that the origin is a Hopf bifurcation point when the conditions $\sigma+1=0$, $\rho-1>0$ and $b\neq0$ hold. This result coincides with the result given in \cite[Section 2]{AlgabaA2005}. We remark that, by using analysis techniques, the authors in \cite{AlgabaA2005} obtained conditions on $(\rho,\sigma,b)$ for the nontrivial equilibria of the Lorenz system to undergo a Hopf bifurcation, see \cite[Proposition 2]{AlgabaA2005}. Here, using triangular decomposition of semi-algebraic systems, we have found a simple closed form description involving three free parameters (see $R_2$ and $R_3$), which coincides (after some elementary transformation) with the result of a manual computation given in \cite[Proposition 2]{AlgabaA2005}.

Similarly, by solving the semi-algebraic set $\mathbb{H}_{3,1}$ and $\mathbb{H}_{3}$ for the variable ordering $x>y>z>\rho>\sigma>{\color{red}b}$, we obtain the corresponding regular semi-algebraic systems for the Lorenz system to have {\bf CHB} and {\bf SHB}. The output results are the following:

{\bf CHB:}\quad $\longleftarrow$\quad $\mathbb{H}_{3,1}$
\begin{equation}\label{Lorenz-R2}
\begin{split}
\mathcal{C}_{{\rm CHB}}=\big(R_2^T\,\,\wedge \,\,R_2^{\mathcal{Q}}\big)\,\,\vee\,\,
\big(R_{3}^T\,\,\wedge \,\,R_{3}^{\mathcal{Q}}\big).
\end{split}
\end{equation}

{\bf SHB:}\quad $\longleftarrow$\quad $\mathbb{H}_{3}$
\begin{equation}\label{Lorenz-R3}
\begin{split}
\mathcal{C}_{{\rm SHB}}=R_3^T\,\,\wedge \,\,R_3^{\mathcal{Q}}\,\,\wedge \,\,(\sigma+b+1>0).
\end{split}
\end{equation}

\begin{remark}\label{Hopf-Lcomp}
Here we take the Lorenz system as an example to compare our approach with three other methods in the literature for studying Hopf bifurcation points.
\begin{itemize}
  \item {\bf Quantifier elimination procedure.} The algebraic criterion for detecting Hopf bifurcations points has been introduced in \cite{JMAA:1994,Weber2000}. More concretely, the semi-algebraic description of the set of parameter values for which a Hopf bifurcation point occurs for the Lorenz system can be expressed by the following first-order formula:
\begin{equation}\label{Hopf-QE-L1}
\begin{split}
\exists x \exists y \exists z \,\big(&
\sigma(y-x)=0\,\wedge\,\rho x-y-xz=0\, \wedge \, -bz+xy=0\\
&\wedge \,
\Delta_2=\left|
 \begin{matrix}
    a_1& a_3 \\
   1 & a_2 \\
  \end{matrix}
  \right|=0\, \wedge \, a_3\cdot a_1>0\big),
\end{split}
\end{equation}
where $a_1$, $a_2$ and $a_3$ are the coefficients in \eqref{Lorenz-Cp}. In order to decide or to reduce the generated first-order formulae we can use tools from computational logic, such as \textsf{QEPCAD B} \cite{qepcadb} (an improvement of \textsf{QEPCAD} \cite{Collins_Hong:91}),
functions \textsf{Resolve} and \textsf{Reduce} in \textsc{Mathematica} (as described in \cite{AS:2010}), the \textsf{REDLOG} package in \textsc{Reduce} (implemented on the basis of the methods of cylindrical algebraic decomposition \cite{Collins:75b,Collins_Hong:91} and virtual term substitution \cite{VW:1994,VW:1997}), and the \textsc{Maple} package \textsf{RegularChains} (a quantifier elimination algorithm by cylindrical algebraic decomposition based on regular chains \cite{che2016qua}).

By imposing extra positivity conditions on the parameters $\rho,\sigma,b$ in \eqref{Hopf-QE-L1},
the following result is obtained by using \textsf{QuantifierElimination} command of the \textsf{RegularChains} library in \textsc{Maple} for the test of a Hopf bifurcation point:
\begin{equation*}% \label{Hopf-QE-L2}
\begin{split}
\rho  b -\rho  \sigma +b \sigma +\sigma^{2}+\rho +3 \sigma=0\, \wedge\, b>0\, \wedge\, \sigma-1>0\, \wedge\, -1+\sigma-b>0.
\end{split}
\end{equation*}
We remark that using the quantifier elimination procedure, one can detect the parameter conditions for the existence of a Hopf bifurcation point of a given differential system. However, this procedure does not reveal any information about the equilibria of the differential system, that is, which equilibrium point of the system can eventually become a Hopf bifurcation point. In contrast, our method of triangular decomposition of semi-algebraic systems can provide such information on Hopf bifurcation fixed point. As a byproduct, we can directly read the Hopf curve of the corresponding equilibria of the system.

    \item {\bf Elimination ideal approach.} The techniques used here follow from the ideas in \cite{LiVR2019,CPVR2020}. Such a method for studying Hopf bifurcations can also be found in \cite{NKSW:2020}. Assume that the characteristic polynomial $p(\lambda)$ in \eqref{Lorenz-Cp} has a pair of imaginary roots $\lambda_{1,2}=\pm i\omega$ and one real root $\lambda_3=b_0$, then $p(\lambda)$ can be written in the form
\begin{equation}\label{Lorenz-Cp1}
\bar{p}(\lambda)=(\lambda-b_0)(\lambda^2+\omega^2).
\end{equation}
Equating the coefficients of $\lambda$ on both sides of the equality $p(\lambda)=\bar{p}(\lambda)$ yields the system of algebraic equations
\begin{equation}\label{Lorenz-P1}
\begin{split}
g_1&=-b \rho  \sigma +b \sigma  z +\omega^{2} b_{0}+x^{2} \sigma +y \sigma  x +b \sigma=0,\\
g_2&=b \sigma -\omega^{2}-\rho  \sigma +\sigma  z +x^{2}+b +\sigma=0,\\
g_3&=b +1+\sigma +b_{0}=0.
\end{split}
\end{equation}

The characteristic polynomial $p(\lambda)$ can be represented in the form \eqref{Lorenz-Cp1} only for
those values of parameters and equilibria of system \eqref{Lorenz-1} for which system \eqref{Lorenz-P1} has a solution. Note that $b_0$ and $\omega$ should be different from zero. To take this condition into account we add to system \eqref{Lorenz-P1} the equations
\begin{equation}\label{Lorenz-P2}
\begin{split}
h_1=1-\nu_1b_0=0,\quad h_2=1-\nu_2\omega=0,
\end{split}
\end{equation}
where $\nu_1$ and $\nu_2$ are new variables. We denote by $\mathcal{I}$ the ideal generated by the polynomials:
\[\mathcal{I}=\langle \sigma(y-x),\rho x-y-xz,-bz+xy,g_1,g_2,g_3,h_1,h_2\rangle.\]
We then have to eliminate from $\mathcal{I}$ the variables $b_0$, $\omega$, $\nu_1$ and $\nu_2$. By computing the fourth elimination ideal of $\mathcal{I}$ in the polynomial ring $\mathbb{C}[b_0,\omega,\nu_1,\nu_2,b, \rho, \sigma, x, y, z]$, we have
$$\mathcal{I}_4= \mathcal{I} \cap \mathbb{C}[b, \rho, \sigma, x, y, z].$$
Computing a prime decomposition of the elimination ideal $\mathcal{I}_4$, we obtain a sequence of irredundant prime ideals:
\begin{equation*}% \label{Lorenz-IP1}
\begin{split}
\mathcal{H}_1&=\langle x, y, z, \sigma +1\rangle,\\
\mathcal{H}_2&=\langle x, y, z, -b^{2}-b \sigma +\rho  \sigma -b -\sigma\rangle,\\
\mathcal{H}_3&=\langle y -x, \rho -z -1, b z -x^{2}, \sigma^{2} z +x^{2} \sigma -\sigma  \,z^{2}+x^{2} z +2 \sigma  z +x^{2}+z^{2}+z
,\\
&\quad\,\,b \sigma +\sigma^{2}-\sigma  z +x^{2}+b +2 \sigma +z +1
\rangle.
\end{split}
\end{equation*}
By substituting the generators of $\mathcal{H}_1$ into the polynomial $p(\lambda)$, we have
\[p(\lambda)=\left(\lambda^{2}+\rho -1\right) \left(\lambda+b \right).\]
So, the origin is a Hopf bifurcation point when the conditions $\sigma+1=0$, $\rho-1>0$ and $b\neq0$ hold. However, using a sample point $[b=0,\rho=1]$ in $\mathcal{H}_2$ we obtain the polynomial $p(\lambda)=\lambda^{2} \left(\lambda +1+\sigma \right)$. As a consequence, $p(\lambda)$ can not have a simple pair of purely imaginary roots though we impose the equation $h_2=0$ in the ideal $\mathcal{I}$.

In fact, by the Closure Theorem for ``almost all'' values of $b, \rho, \sigma, x, y, z$ satisfying the condition $\mathcal{H}_i=0$ the corresponding variety $V(\mathcal{I})$ has a solution. However it can happen that for some values of parameters and variables it does not hold. Using Extension Theorem, we may obtain necessary and sufficient conditions for the existence of a simple pair of purely imaginary roots of $p(\lambda)$. For more details about the theorems, we refer the reader to \cite[Chapter 3]{cox2013ide}. In contrast, using our method of triangular decomposition of semi-algebraic systems, we can directly obtain the desired results for the analysis of Hopf bifurcations.

\item {\bf Hopf algorithms using polynomial resultants.} The techniques used here follow from the ideas in \cite{Guckenheimer:1996,Guckenheimer:1997}. Guckenheimer and his co-workers discussed classical algebraic constructions for determining matrices of arbitrary dimension with pure imaginary eigenvalues. A single equality condition and an inequality in the coefficients of the characteristic polynomial are formulated in terms of resultants, determinants of special matrices whose elements are functions of the characteristic polynomial coefficients. See \cite[Section 2]{Guckenheimer:1997} for more details about the method for detecting curves of Hopf bifurcation points. We remark that the criteria developed in \cite{Guckenheimer:1996,Guckenheimer:1997} can only provide necessary conditions for the Hopf bifurcation points as the real part of the eigenvalues other than a pair of pure imaginary eigenvalues maybe zero. For example, under the condition $b=0$, $\sigma=-1$ and $\rho-1>0$, it is easy to check that the characteristic polynomial $p(\lambda)$ in \eqref{Lorenz-Cp} at the origin is expressed as
\begin{equation*}
p(\lambda)=\lambda^{3}+(\rho -1)\lambda.
\end{equation*}
In this case, $p(\lambda)$ has a pair of pure imaginary eigenvalues $\lambda_{1,2}=\pm i\sqrt{\rho-1}$ and one eigenvalue $\lambda_3=0$. The coefficients of $p(\lambda)$ satisfy the conditions of \cite[Theorem 2.1]{Guckenheimer:1997}, but do not satisfy the conditions for Hopf bifurcation. In contrast, using our symbolic method, we can derive necessary and sufficient conditions for Hopf bifurcation points.

\end{itemize}

\end{remark}

\begin{figure}
\centering
\begin{minipage}[t]{0.48\linewidth}
  \vspace{0pt}
		\centering
		\includegraphics[width=0.99\linewidth]{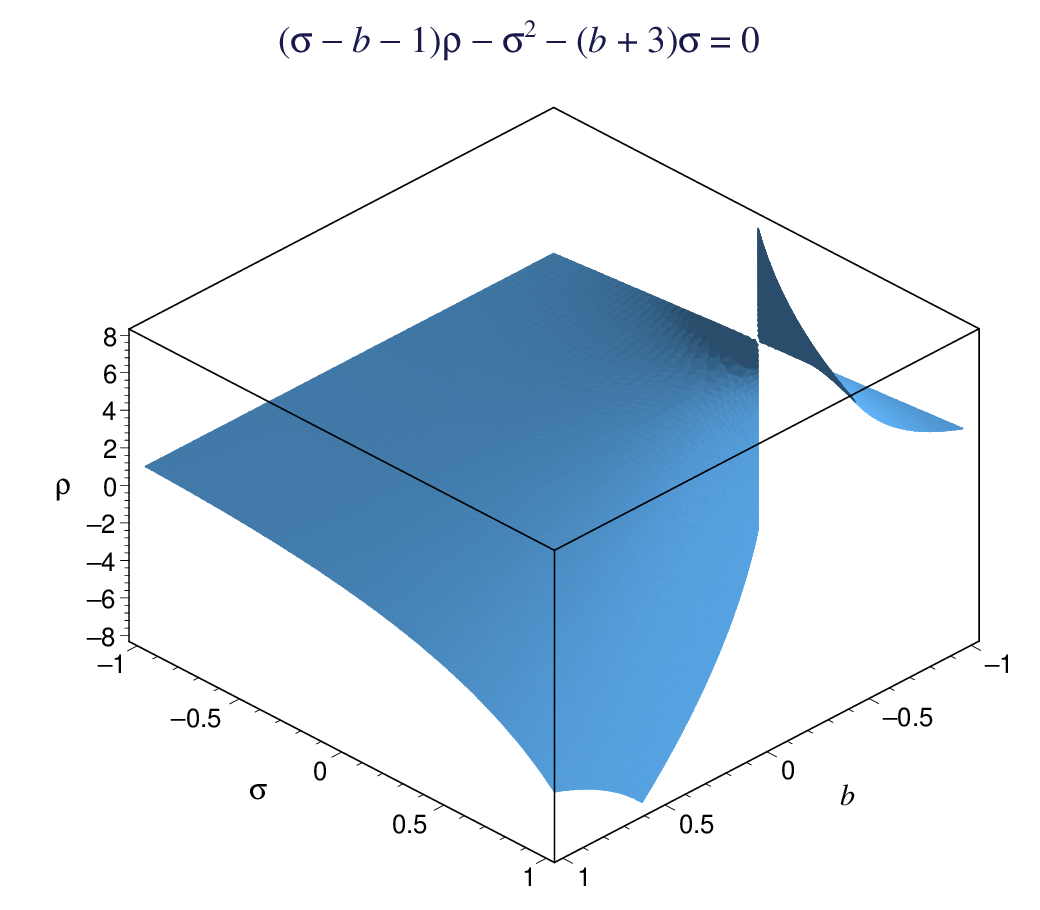}
		\parbox{\linewidth}{\centering (a) Hopf surface: certain patches of $R_3^T$\par}
	\end{minipage}\hfill
\begin{minipage}[t]{0.48\linewidth}
  \vspace{0pt}
		\centering
		\includegraphics[width=0.99\linewidth]{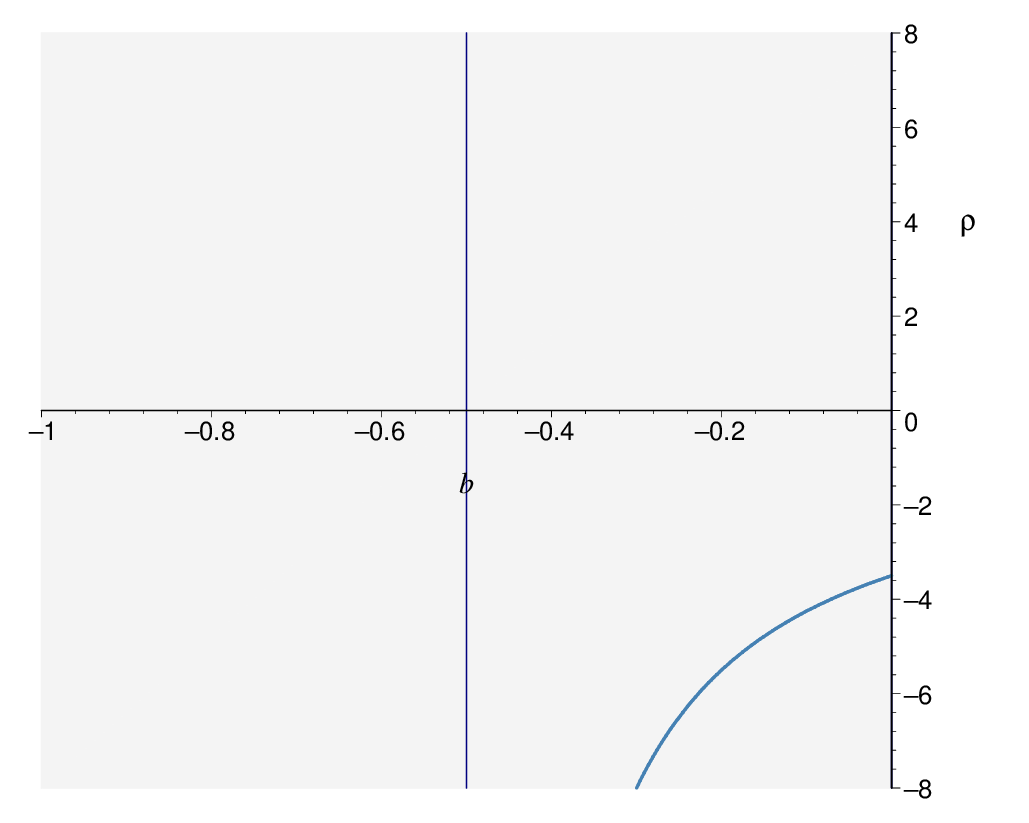}
		\parbox{\linewidth}{\centering (b) Hopf curve when $\sigma=\frac{1}{2}$\par}
	\end{minipage}\par\medskip
\begin{minipage}[t]{0.48\linewidth}
  \vspace{0pt}
\centering
\includegraphics[width=1\linewidth]{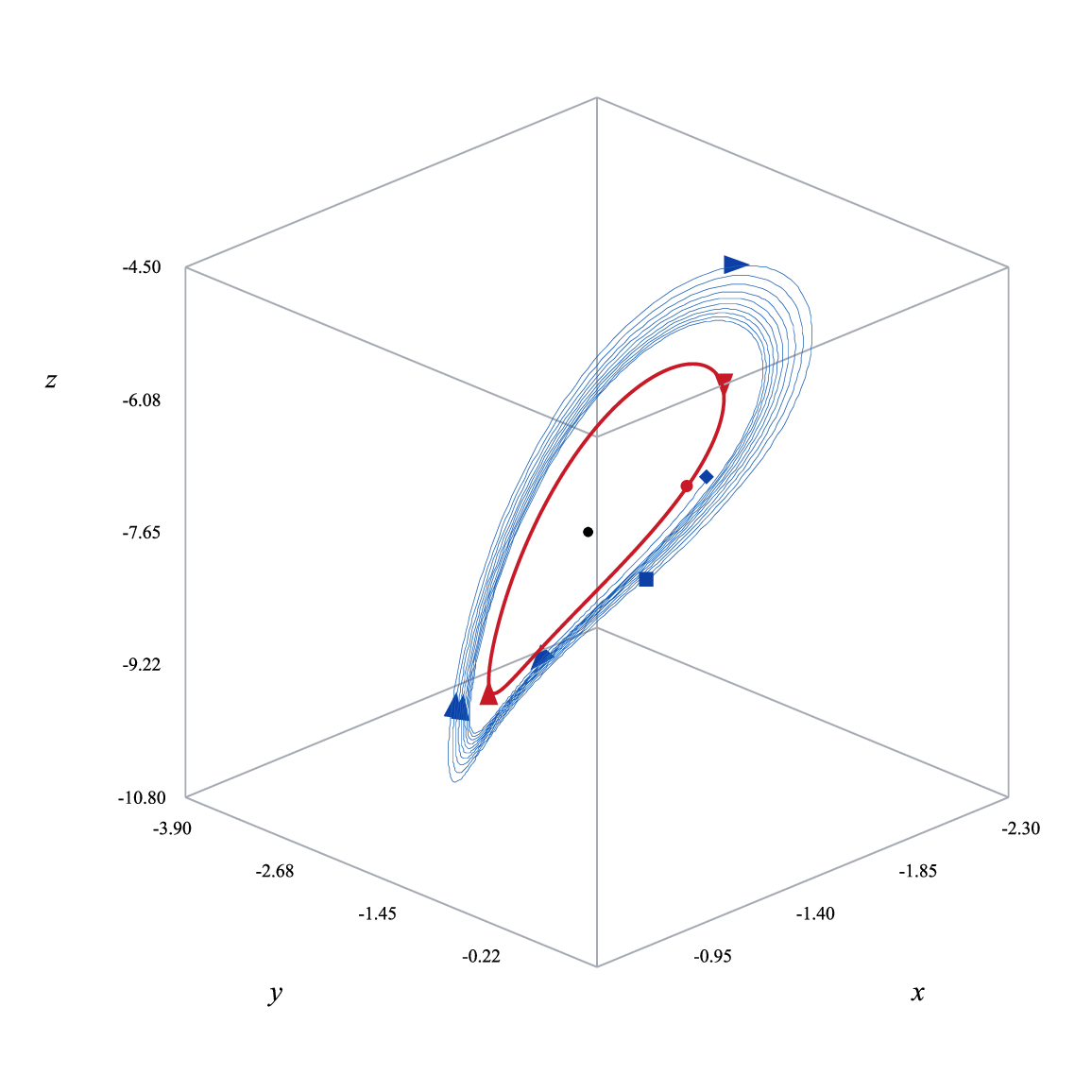}
\parbox{\linewidth}{\centering (c) {\bf SHB} of system \eqref{Lorenz-1} at the point $\boldsymbol{x}_0=(-{\frac{11762252359}{8589934592}},-{\frac{376392075483}{274877906944}},-\frac{15}{2})$ under the condition $\boldsymbol{\mu}_0=(-\frac{13}{2},\frac{1}{2},-\frac{1}{4})$. The parameters chosen for the numerical simulation are: $\rho=-6.5$, $\sigma=0.5$ and $b=-0.23$. The initial values $(x(0), y(0), z(0))$ in blue curve are (-1.31339255365637, 0.095533439591827, -6.26429818819526); and in red are (-1.31339255365637, -0.139287559282873, -6.47024849016272). The bifurcated limit cycle is unstable.\par}
\end{minipage}\hfill
\begin{minipage}[t]{0.48\linewidth}
  \vspace{0pt}
\centering
\includegraphics[width=1\linewidth]{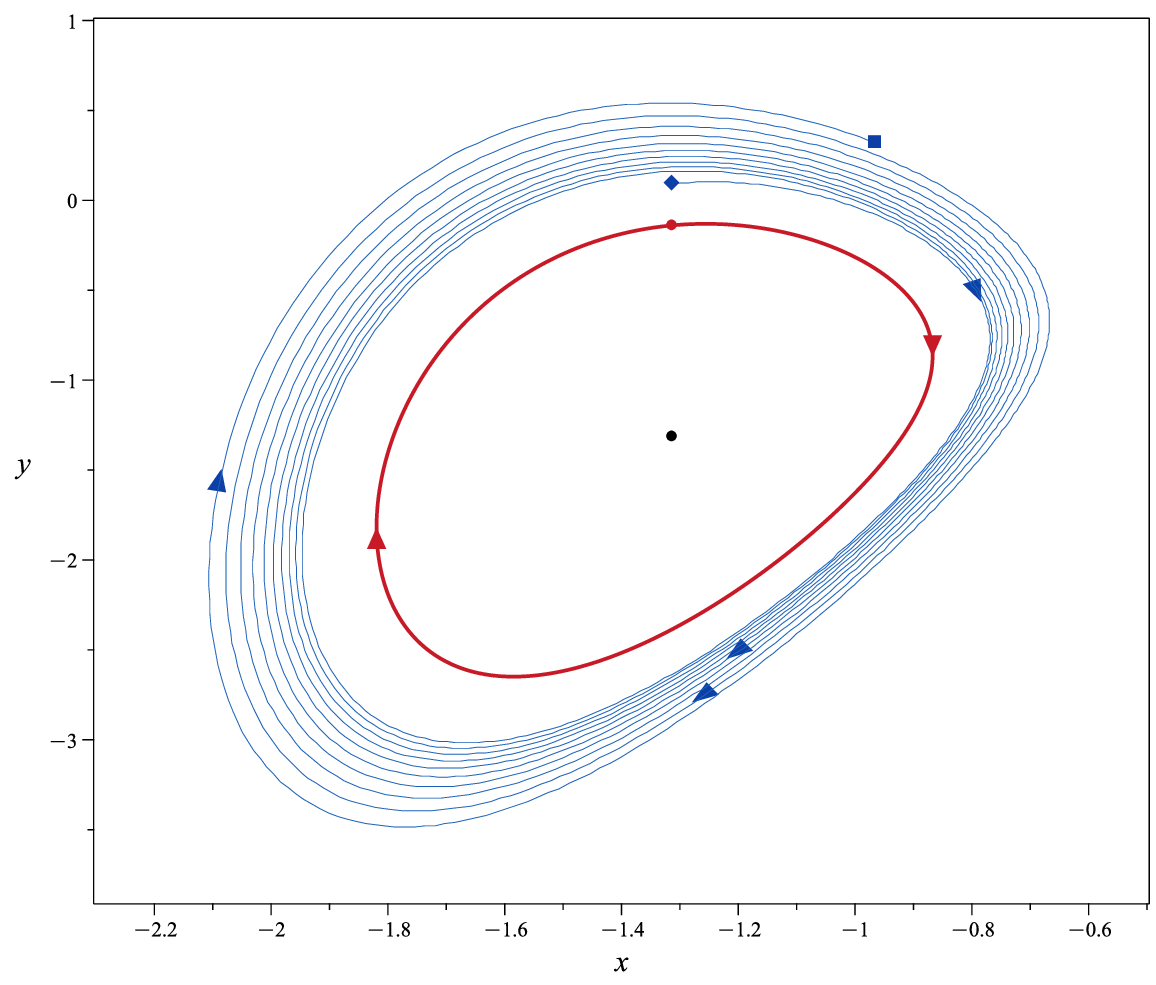}
\parbox{\linewidth}{\centering (d) Projection of the bifurcated limit cycle onto the $x-y$ plane.\par}
\end{minipage}

\caption{Hopf curve and Hopf bifurcation of system \eqref{Lorenz-1}.}\label{fig:L-1}
\end{figure}

In what follows, we characterize the nature of Hopf bifurcations in the differential system \eqref{Lorenz-1}. To this end, we employ the {\sf SamplePoints} routine, a built-in function within the {\sf RegularChains} library in \textsc{Maple}, to solve the semi-algebraic system $\mathbb{H}_{3}$ formulated in \eqref{Hopf-main-L3}. This computation yields four distinct sample points $(\boldsymbol{x}_0,\boldsymbol{\mu}_0)$, as summarized in Table \ref{tab:ex1-sample}, at which the differential system \eqref{Lorenz-1} admits simple Hopf bifurcations. For each sample point $(\boldsymbol{x}_0,\boldsymbol{\mu}_0)$, we execute Algorithm \ref{alg:lyapunov} to numerically compute the first Lyapunov coefficient $\ell_1(\mathbf{x}_0)$, whose corresponding values are presented in Table \ref{tab:ex1-sample}. Note that the sign of $\ell_1(\boldsymbol{x}_0)$ gives stability of the small bifurcated limit cycle on the center manifold. The numerical simulation results for the Hopf curve and simple Hopf bifurcation of system \eqref{Lorenz-1} are illustrated in Figure \ref{fig:L-1}.

\begin{table}[H]
	\centering
	\caption{Sample points of $(\boldsymbol{x}_0,\boldsymbol{\mu}_0)$ corresponding to {\bf SHB} of system \eqref{Lorenz-1}.}
	\label{tab:ex1-sample}
\begin{tabular}{c|c|c}
  \hline
  $\boldsymbol{x}_0$ & $\boldsymbol{\mu}_0$ &  $\ell_1(\boldsymbol{x}_0)$ \\
  \hline
  \noalign{\vskip 4pt}
  $(-{\frac{11762252359}{8589934592}},-{\frac{376392075483}{274877906944}},-\frac{15}{2})$ & $(-\frac{13}{2},\frac{1}{2},-\frac{1}{4})$ &  $0.007408529250>0$ \\[3pt]
  $({\frac{23524504717}{17179869184}},{\frac{188196037741}{137438953472}},-\frac{15}{2})$ & $(-\frac{13}{2},\frac{1}{2},-\frac{1}{4})$ &  $0.007408529250>0$ \\[3pt]
    $({\frac{22522998513}{8589934592}},{\frac{720735952423}{274877906944}},\frac{55}{2})$ & $(\frac{57}{2},\frac{3}{2},\frac{1}{4})$ &  $0.007366251735>0$ \\[3pt]
   $(-{\frac{45045997027}{17179869184}},-{\frac{90091994053}{34359738368}},\frac{55}{2})$ & $(\frac{57}{2},\frac{3}{2},\frac{1}{4})$ &  $0.007366251735>0$ \\[3pt]
    \hline
\end{tabular}
\end{table}

\subsection{Other models and remarks}\label{sect-6-4}
\subsubsection{The Rabinovich--Fabrikant model}

We study the following Rabinovich--Fabrikant model proposed in \cite{RF1979}:
\begin{equation}\label{RFM-1}
\left\{\begin{array}{l}
\dot{x}=y \left(x^{2}+z -1\right)+\beta  x,\\
\dot{y}=x \left(-x^{2}+3 z +1\right)+\beta  y,\\
\dot{z}=-2 z \left(x y +\alpha \right),
\end{array}\right.
\end{equation}
where
$\boldsymbol{x}=(x,y,z)$ are variables, and $\boldsymbol{\mu}=(\alpha,\beta)$ are real parameters with $\alpha,\beta>0$. This system models the stochasticity arising from the modulation instability in a nonequilibrium dissipative medium.

In a similar way, by running the {\sf RealTriangularize} command on the resulting semi-algebraic set $\mathbb{S}_{3,1}$, $\mathbb{H}_{3,1}$ and $\mathbb{H}_{3}$ for the variable ordering $x>y>z>\beta>{\color{red}\alpha}$, we obtain the following regular semi-algebraic systems for {\bf HBP}, {\bf CHB} and {\bf SHB}.

{\bf HBP:}\quad $\longleftarrow$\quad $\mathbb{S}_{3,1}$
\begin{equation}\label{RF-HBP}
\begin{split}
\mathcal{C}_{{\rm HBP}}=\big(\bar{R}_{1}^T\,\,\wedge\,\,\bar{R}_{1}^{Q}\big)\,\,\vee\,\,\big(\bar{R}_{2}^T\,\,\wedge\,\,\bar{R}_{2}^{Q}\big),
\end{split}
\end{equation}
where
\begin{equation*}
\begin{split}
\bar{R}_1^T&=\left\{\begin{array}{l}
yx+\alpha=0\\
y^2-\alpha=0\\
z-\alpha+1=0\\
\beta-2\alpha+2=0\\
\end{array}\right., \quad
\bar{R}_1^Q=\left\{\begin{array}{l}
\alpha-1>0\\
\alpha-2\neq0
\end{array}\right.,\\
\bar{R}_2^T&=\left\{\begin{array}{l}
yx+\alpha=0\\
\left(z -1\right) y^{2}+\alpha^{2}-\alpha  \beta=0\\
\left(\beta +\alpha \right) z -\beta=0\\
\beta^{4}-2 \alpha^{2} \beta^{2}+4 \alpha  \beta+\alpha^{4}-3 \alpha^{2}=0\\
\end{array}\right., \\
\bar{R}_2^Q&=\left\{\begin{array}{l}
-3 \alpha^{8}+4 \beta  \,\alpha^{7}+3 \beta^{2} \alpha^{6}-4 \beta^{3} \alpha^{5}+9 \alpha^{6}-18 \beta  \,\alpha^{5}+8 \beta^{2} \alpha^{4}>0\\
\beta>0\,\wedge\, \alpha^2-3<0\, \wedge\, \alpha>0
\end{array}\right..
\end{split}
\end{equation*}

{\bf CHB:}\quad $\longleftarrow$\quad $\mathbb{H}_{3,1}$
\begin{equation}\label{RF-CHB}
\begin{split}
\mathcal{C}_{{\rm CHB}}=&\big(\bar{R}_{1}^T\,\,\wedge\,\,\bar{R}_1^Q \,\,\wedge\,\,(3\alpha-4\neq0)\big)\,\,\vee\,\,\\
&\big(\bar{R}_{2}^T\,\,\wedge\,\,\bar{R}_2^Q\,\,\wedge\,\,\big(3\alpha-4\neq0\,\,\wedge\,\,112 \alpha^{6}-264 \alpha^{4}+183 \alpha^{2}-48\neq0\big)\big).
\end{split}
\end{equation}

{\bf SHB:}\quad $\longleftarrow$\quad $\mathbb{H}_{3}$
\begin{equation}\label{RF-SHB}
\begin{split}
\mathcal{C}_{{\rm SHB}}=\big(\bar{R}_{1}^T\,\,\wedge\,\,\bar{R}_{1}^{\bar{Q}}\big)\,\,\vee\,\,\big(\bar{R}_{2}^T\,\,\wedge\,\,\bar{R}_{2}^{\bar{Q}}\big),
\end{split}
\end{equation}
where
\begin{equation*}
\begin{split}
&\bar{R}_{1}^{\bar{Q}}=\left\{\begin{array}{l}
\alpha>1\\
\alpha<2\\
3\alpha-4\neq0
\end{array}\right.,\\
&\bar{R}_{2}^{\bar{Q}}
=\left\{\begin{array}{l}
-2 \alpha^{7}+2 \beta  \,\alpha^{6}+3 \alpha^{5} \beta^{2}-\beta^{3} \alpha^{4}+6 \alpha^{5}-8 \beta  \,\alpha^{4}>0\\
\alpha-\beta>0\\
\beta>0\\
\alpha^2-3<0\\
\alpha>0\\
3\alpha-4\neq0\\
112 \alpha^{6}-264 \alpha^{4}+183 \alpha^{2}-48\neq0
\end{array}\right..
\end{split}
\end{equation*}

{
As $\bar{R}_{1}$ satisfies the condition at line $11$ of Algorithm~\ref{algorithm:rtd}, the parameter conditions for
it to have real solutions are exactly the conjunctions
of its constraints involving the parameters $\alpha$ and $\beta$.
For $\bar{R}_{2}$, the condition at line $14$ does not hold.
To check if the condition at line $13$ still holds, it is equivalent
to verify if the discriminant of the only quadratic (in main variables) polynomial
$\left(z -1\right) y^{2}+\alpha^{2}-\alpha  \beta$ is nonnegative,
which turns out to be {\em true} by using \textsf{QuantifierElimination} command of the \textsf{RegularChains} library in \textsc{Maple}.
Thus, $\bar{R}_{2}$, the parameter conditions for
it to have real solutions are also exactly the conjunctions
of its constraints involving the parameters $\alpha$ and $\beta$.
}

Applying the {\sf SamplePoints} routine to the resulting semi-algebraic system $\mathbb{H}_{3}$ yields ten distinct sample points $(\boldsymbol{x}_0,\boldsymbol{\mu}_0)$. Selecting four of these points, we compute the desired Lyapunov coefficient $\ell_1(\boldsymbol{x}_0)$ for each via Algorithm \ref{alg:lyapunov}. The computational results are summarized in Table \ref{tab:ex2-sample}. For readability, all large rational numbers presented here represent approximate values.

\begin{table}[H]
	\centering
	\caption{Sample points of $(\boldsymbol{x}_0,\boldsymbol{\mu}_0)$ corresponding to {\bf SHB} of system \eqref{RFM-1}.}
	\label{tab:ex2-sample}
\begin{tabular}{c|c|c}
  \hline
  $\boldsymbol{x}_0$ & $\boldsymbol{\mu}_0$ &  $\ell_1(\boldsymbol{x}_0)$ \\
  \hline
  \noalign{\vskip 4pt}
$(\frac{37071315507}{34359738368}, -\frac{9267828877}{8589934592}, \frac{21}{128})$ & $(\frac{149}{128}, \frac{21}{64})$ &  $-0.1358788031<0$ \\[3pt]
  $(\frac{433355034582832159}{335544320000000000}, -{\frac{22187777773}{17179869184}}, {\frac{171}{256}})$ & $({\frac{427}{256}}, {\frac{171}{128}})$ &  $-2.188145314<0$ \\[3pt]
  $({\frac{245816212559970767}{167772160000000000}}, -{\frac{4249256123014107}{10485760000000000}}, {\frac{232138666295}{549755813888}})$ & $({\frac{19}{32}}, {\frac{1908565127415}{4398046511104}})$ &  $0.08969396095>0$ \\[3pt]
  $({\frac{20846503569}{17179869184}}, -{\frac{285817634399}{274877906944}}, {\frac{101782471623479}{281474976710656}})$ & $({\frac{323}{256}}, {\frac{1571576075581}{2199023255552}})$ &  $0.2048564397>0$ \\[3pt]
    \hline
\end{tabular}
\end{table}

\subsubsection{A four-dimensional business cycle model}

We consider the following four-dimensional business cycle model suggested by Puu \cite{Puu2003}:
\begin{equation}\label{4dim-equ-1}
\begin{split}
\dot{Y}_1 &= W_1, \\
\dot{Y}_2 &= W_2, \\
\dot{W}_1 &= -(s_1 + m_1)Y_1 + m_2 Y_2 + (v_1 - 1 - s_1)W_1 - \frac{v_1}{3}W_1^3, \\
\dot{W}_2 &= m_1 Y_1 - (s_2 + m_2)Y_2 + (v_2 - 1 - s_2)W_2 - \frac{v_2}{3}W_2^3,
\end{split}
\end{equation}
where the subscripts $i\in\{1,2\}$ denote the economic regions (or the
countries), $\boldsymbol{x}=(Y_1,Y_2,W_1,W_2)$ are variables and $\boldsymbol{\mu}=(m_1,m_2,s_1,s_2,v_1,v_2)$ are real parameters. More concretely, $Y_i$=regional income, $W_i$=rate of change of regional income, $m_i$=marginal propensity to import $(0<m_i<1)$, $s_i$=propensity to save $(0<s_i<1)$ and $\nu_i$=accelerator
coefficient of investment $(v_i>0)$.

To better display the output results, in practice we set $[m_1=m_2=s_2=s_1]$ and $v_2=2v_1$ with $s_1\in(0,1)$ and $v_1>0$. By running the {\sf RealTriangularize} command on the resulting semi-algebraic set $\mathbb{S}_{4,1}$, $\mathbb{H}_{4,1}$ and $\mathbb{H}_{4}$ for the variable ordering $Y_1>Y_2>W_1>W_2>s_1>{\color{red}v_1}$, we obtain the following regular semi-algebraic systems for {\bf HBP}, {\bf CHB} and {\bf SHB}.

{\bf HBP:}\quad $\longleftarrow$\quad $\mathbb{S}_{4,1}$

\begin{equation}\label{4D-HBP}
\begin{split}
\mathcal{C}_{{\rm HBP}}=\big(\widetilde{R}_{1}^T\,\,\wedge\,\,\widetilde{R}_{1}^{Q}\big)\,\,\vee\,\,\big(\widetilde{R}_{2}^T\,\,\wedge\,\,\widetilde{R}_{2}^{Q}\big),
\end{split}
\end{equation}
where
\begin{equation*}
\begin{split}
\widetilde{R}_1^T&=\left\{\begin{array}{l}
Y_1=0,\quad Y_2=0\\
W_1=0,\quad W_2=0\\
2s_1-2v_1+1=0
\end{array}\right., \quad
\widetilde{R}_1^Q=\left\{\begin{array}{l}
3-2v_1>0\\
2v_1-1>0\\
v_1-1\neq0
\end{array}\right.,\\
\widetilde{R}_2^T&=\left\{\begin{array}{l}
Y_1=0,\quad Y_2=0\\
W_1=0,\quad W_2=0\\
s_1-2v_1+2=0
\end{array}\right., \quad
\widetilde{R}_2^Q=\left\{\begin{array}{l}
3-2v_1>0\\
v_1-1>0
\end{array}\right..
\end{split}
\end{equation*}

{\bf CHB:}\quad $\longleftarrow$\quad $\mathbb{H}_{4,1}$

\begin{equation}\label{4D-CHB}
\begin{split}
\mathcal{C}_{{\rm CHB}}=\big(\widetilde{R}_{1}^T\,\,\wedge\,\,\widetilde{R}_{1}^{Q}\big)\,\,\vee\,\,\big(\widetilde{R}_{2}^T\,\,\wedge\,\,\widetilde{R}_{2}^{Q}\big).
\end{split}
\end{equation}

{\bf SHB:}\quad $\longleftarrow$\quad $\mathbb{H}_{4}$

\begin{equation}\label{4D-SHB}
\begin{split}
\mathcal{C}_{{\rm SHB}}=\widetilde{R}_{1}^T\,\,\wedge\,\,\big(1-v_1>0\,\,\wedge\,\,2v_1-1>0\big).
\end{split}
\end{equation}

\begin{remark}\label{4D-remark}
Applying {\sf RealTriangularize} to the resulting semi-algebraic set $\mathbb{S}_{4,2}$ yields an empty output, implying that system \eqref{4dim-equ-1} admits {\bf HBP}s only within $\mathbb{S}_{4,1}$. Indeed, even with six free parameters $\boldsymbol{\mu}$ in the economically meaningful domain $[0<m_i<1,0<s_i<1,v_i>0]$, we have $a_4=s_{2} m_{1}+m_{2} s_{1}+s_{1} s_{2}>0$, which confirms via Theorem \ref{Hopf-th3} that $\mathbb{S}_{4,2}$ is empty. On the other hand, if $s_1, v_1$, and $v_2$ are kept as free parameters, {\sf RealTriangularize} decomposes the system for {\bf HBP}s into 28 regular semi-algebraic branches (taking about 2.484 seconds under the variable ordering $Y_1>Y_2>W_1>W_2>s_1>v_2>{\color{red}v_1}$), which are omitted here for brevity. As observed, the required symbolic computation easily becomes prohibitively expensive and intractable as the number of free parameters increases.

\end{remark}

The experimental results of Algorithm \ref{alg:lyapunov} related to the semi-algebraic system $\mathbb{H}_{4}$ are summarized in Table \ref{tab:ex3-sample}. The computational times of {\bf HBP}, {\bf CHB} and {\bf SHB} for the three differential models are presented in Table \ref{Tab-examples}.

\begin{table}[H]
	\centering
	\caption{Sample points of $(\boldsymbol{x}_0,\boldsymbol{\mu}_0)$ corresponding to {\bf SHB} of system \eqref{4dim-equ-1}.}
	\label{tab:ex3-sample}
\begin{tabular}{c|c|c}
  \hline
  $\boldsymbol{x}_0$ & $\boldsymbol{\mu}_0$ &  $\ell_1(\boldsymbol{x}_0)$ \\
  \hline
  \noalign{\vskip 4pt}
  $(0,0,0,0)$ & $(\frac{1}{6},\frac{1}{6},\frac{1}{6},\frac{1}{6},\frac{2}{3},\frac{4}{3})$ &  $-0.7116035710<0$ \\[3pt]
    \hline
  $(0,0,0,0)$ & $(\frac{1}{4},\frac{1}{4},\frac{1}{4},\frac{1}{4},\frac{3}{4},\frac{3}{2})$ &  $-1.421874998<0$ \\[3pt]
    \hline
  $(0,0,0,0)$ & $(\frac{3}{10},\frac{3}{10},\frac{3}{10},\frac{3}{10},\frac{4}{5},\frac{8}{5})$ &  $-1.911342351<0$ \\[3pt]
  %  \hline
  $(0,0,0,0)$ & $(\frac{1}{3},\frac{1}{3},\frac{1}{3},\frac{1}{3},\frac{5}{6},\frac{5}{3})$ &  $-2.138888890<0$ \\[3pt]
    \hline
\end{tabular}
\end{table}

\begin{table}[H]
\caption{Computational times (in sec.) of {\bf HBP}, {\bf CHB} and {\bf SHB} for differential models.}\label{Tab-examples}
\smallskip
\begin{center}
\begin{tabular}{cccc}
\toprule
  Models & {\bf HBP} & {\bf CHB} & {\bf SHB}\\ \hline
   Lorenz & 0.219 & 0.281 & 0.250\\ \hline
   Rabinovich--Fabrikant model & 0.531 & 1.594 & 1.016\\ \hline
   Business cycle model & 0.110 & 0.125 & 0.109\\
   \bottomrule
\end{tabular}
\end{center}
\end{table}

\section{Concluding Remarks}\label{sect-6}

In this work, we introduce two algorithmic approaches for determining the parameter conditions under which a parametric polynomial differential system admits a prescribed invariant algebraic surface of bounded degree. Specifically, we combine Gr\"obner bases with saturation and prime decomposition to obtain a concise representation of the Zariski closure of the existence locus (necessary conditions), whereas we project the nonconstant solutions via triangular decomposition to recover the exact constructible locus (necessary and sufficient conditions). We report the performance of the two methods in Table \ref{table:dsys-exp} using the BZ reaction model. It turns out that, triangular decomposition method is more efficient and flexible, making it particularly suitable for computing higher-degree invariant algebraic surfaces or imposing additional constraints on system parameters. In addition, an alternative method of combining Gr\"obner bases with triangular decomposition is discussed. As an application, certain types of invariant algebraic surfaces of a class of Li\'enard differential systems are characterized. It is shown that our algorithmic approach will work more conveniently and effectively once we know the structure of invariant algebraic curves of the differential system. The obtained results indicate that the proposed approach allows us to detect hand-computation errors in published literature.

We provide a systematic approach for studying Hopf bifurcations by using method of triangular decomposition of semi-algebraic systems. We give a complete, necessary and sufficient characterization of Hopf bifurcation points via two equivalent formulations: one based on the properties of Hurwitz determinants and subresultant sequences, and the other governed by the greatest common divisor of the characteristic polynomial's real and imaginary parts. We derive a semi-algebraic description of the transversality condition for Hopf bifurcation via the polynomial GCD. This description naturally yields an algorithmic decision procedure for detecting Hopf bifurcations in the considered systems. Moreover, we present a symbolic-numeric projection algorithm for evaluating the first Lyapunov coefficient. Experimental results demonstrate the effectiveness of the proposed algorithmic approach.

Several promising directions remain for future research. First, extending our symbolic-numeric framework to certify the sign of the first Lyapunov coefficient using validated interval arithmetic would provide rigorous computer-assisted proofs of supercriticality and subcriticality. Second, when the genericity condition fails (i.e., $\ell_1(\boldsymbol{x}_0)=0$), generalizing our semi-algebraic and projection techniques to higher-order Lyapunov coefficients offers a systematic pathway toward the algorithmic detection of codimension-two Bautin (generalized Hopf) bifurcations. Finally, exploiting structural sparsity and symmetry in parametric polynomial systems could further reduce the complexity of branch-wise triangular decomposition in higher dimensions.

\section*{Acknowledgements}
The second author would like to thank Professor Dingkang Wang for the valuable discussions during the ``Symposium on Symbolic Computation and Beyond'' held in Hangzhou on June 27--28, 2026. Special thanks are also due to him for sharing his software package for computing the greatest common divisor of parametric polynomials. The central mathematical ideas and overall approach are due to the authors. The authors acknowledge the use of GPT-6-Astra for assistance with the \textsc{Maple} code used in the numerical simulations shown in Figure \ref{fig:L-1}, as well as for proofreading and improving the clarity of the exposition (specifically, the simplicity statements in Theorem \ref{Hopf-th3}, Corollary \ref{coro-th-hp3}, and Theorem \ref{Hopf-gcd1}). The authors take full responsibility for the paper's content and the accuracy of its references.

\bibliographystyle{abbrv}
\bibliography{references}

\end{document}